\documentclass[twoside,10pt]{article}
\usepackage{amssymb,amsthm,amscd,amsmath}
\usepackage{curves}
\renewcommand{\Sp}{\operatorname{Sp}}
\newcommand{\bK}{\bar{K}}
\newcommand{\tK}{\tilde{K}}
\allowdisplaybreaks[1]
\theoremstyle{plain}
\newtheorem{Theorem}{Theorem}[section]
\newtheorem{Lemma}[Theorem]{Lemma}

\newtheorem{Claim}[Theorem]{Claim}
\newtheorem{Proposition}[Theorem]{Proposition}
\newtheorem{Corollary}[Theorem]{Corollary}

\newtheorem*{Theorem*}{Theorem}
\newtheorem*{Lemma*}{Lemma}

\theoremstyle{definition}
\newtheorem*{Definition*}{Definition}
\newtheorem*{Remark*}{Remark}
\newtheorem*{Example*}{Example}
\newtheorem{Definition}[Theorem]{Definition}
\newtheorem{Remark}[Theorem]{Remark}
\newtheorem{Example}[Theorem]{Example}

\newcommand{\isom}{\cong}

\newcommand{\n}{\mathbin{\triangleleft}}

\newcommand{\sn}{\mathbin{\triangleleft\triangleleft}}

\newcommand{\gen}[1]{\mathopen{<}#1\mathclose{>}}
\newcommand{\abs}[1]{\lvert#1\rvert}

\DeclareMathOperator{\core}{core}
\DeclareMathOperator{\Aut}{Aut}
\DeclareMathOperator{\Out}{Out}

\DeclareMathOperator{\Gal}{Gal}
\DeclareMathOperator{\Isom}{Isom}
\newcommand{\Nor}{\mathbf{N}}
\newcommand{\Cen}{\mathbf{C}}

\newcommand{\gFit}{\mathbf{F}^*}
\newcommand{\Op}{\mathbf{O}_p}
\renewcommand{\O}[1]{\mathbf{O}_#1}

\DeclareMathOperator{\ind}{ind}

\newcommand{\CCC}{\mathbb C}
\newcommand{\FFF}{\mathbb F}            
\newcommand{\NNN}{\mathbb N}
\newcommand{\PPP}{\mathbb P}            
\newcommand{\QQQ}{\mathbb Q}
\newcommand{\RRR}{\mathbb R}
\newcommand{\ZZZ}{\mathbb Z}
\newcommand{\PP}{\mathbb P^1}           

\newcommand{\hK}{\hat K}

\newcommand{\bQ}{\overline\QQQ}
\newcommand{\hQ}{\hat\QQQ}
\newcommand{\ok}{{\mathcal O}_k}        
\newcommand{\ol}{{\mathcal O}_l}        
\newcommand{\oK}{{\mathcal O}_K}        
\newcommand{\cB}{\mathcal B}            
\newcommand{\cC}{\mathcal C}

\newcommand{\cG}{\mathcal G}
\newcommand{\cH}{\mathcal H}

\newcommand{\cM}{\mathcal M}

\newcommand{\cT}{\mathcal T}

\newcommand{\frakm}{\mathfrak m}        
\newcommand{\frakp}{\mathfrak p}
\newcommand{\frakP}{\mathfrak P}

\newcommand{\bn}[1]{\mathbf #1}         
\newcommand{\s}[1]{\sigma_#1}

\DeclareMathOperator{\gL}{\Gamma L}
\DeclareMathOperator{\GL}{GL}

\DeclareMathOperator{\SL}{SL}

\DeclareMathOperator{\PgL}{P\Gamma L}
\DeclareMathOperator{\PGL}{PGL}
\DeclareMathOperator{\PsL}{P\Sigma L}
\DeclareMathOperator{\PSL}{PSL}

\DeclareMathOperator{\PSp}{PSp}
\DeclareMathOperator{\Sz}{Sz}
\DeclareMathOperator{\AGL}{AGL}
\DeclareMathOperator{\ASL}{ASL}

\newcommand{\M}[2]{\mathrm{M}_{#1#2}}             

\begin{document}
\title{The rational function analogue of a question of Schur and
exceptionality of permutation representations}
\author{
Robert M.~Guralnick
\thanks{Partially supported by the NSF and NATO.}
\and
Peter M\"uller
\thanks{Supported by the DFG.}
\and
Jan Saxl
\thanks{Partially supported by the NATO.}
}
\maketitle
\begin{abstract}In 1923 Schur considered the following problem. Let
  $f\in\mathbb Z[X]$ be a polynomial that induces a bijection on the
  residue fields $\mathbb Z/p\mathbb Z$ for infinitely many primes
  $p$. His conjecture, that such polynomials are compositions of
  linear and Dickson polynomials, was proved by M.~Fried in 1970. Here
  we investigate the analogous question for rational functions, and
  also we allow the base field to be any number field. As a result,
  there are many more rational functions for which the analogous
  property holds. The new infinite series come from rational isogenies
  or endomorphisms of elliptic curves. Besides them, there are
  finitely many sporadic examples which do not fit in any of the
  series we obtain.
  
  The Galois theoretic translation, based on Chebotar\"ev's density
  theorem, leads to a certain property of permutation groups, called
  exceptionality. One can reduce to primitive exceptional groups.
  While it is impossible to explicitly describe all primitive
  exceptional permutation groups, we provide certain reduction
  results, and obtain a classification in the almost simple case.
  
  The fact that these permutation groups arise as monodromy groups of
  covers of Riemann spheres $f:\PP\to\PP$, where $f$ is the rational
  function we investigate, provides genus $0$ systems. These are
  generating systems of permutation groups with a certain
  combinatorial property. This condition, combined with the
  classification and reduction results of exceptional permutation
  groups, eventually gives a precise geometric classification of
  possible candidates of rational functions which satisfy the
  arithmetic property from above. Up to this point, we make frequent
  use of the classification of the finite simple groups.
  
  We do not directly obtain information about the number fields over
  which these examples exist, or if they exist at all. However, except
  for finitely many cases, these remaining candidates are connected to
  isogenies or endomorphisms of elliptic curves. Thus we use results
  about elliptic curves, modular curves, complex multiplication, and
  the techniques used in the inverse regular Galois problem to settle
  these finer arithmetic questions.
\end{abstract}

\section{Introduction} In 1923 Schur \cite{Schur} posed the following
question. Let $f(X)\in\ZZZ[X]$ be a polynomial which induces a
bijection on $\ZZZ/p\ZZZ$ for infinitely many primes $p$. He proved
that if $f$ has prime degree, then $f$ is, up to linear changes over
$\bQ$, a cyclic polynomial $X^k$, or a Chebychev polynomial $T_k(X)$
(defined implicitly by $T_k(Z+1/Z)=Z^k+1/Z^k$). He conjectured that
without the degree assumption, the polynomial is a composition of such
polynomials. This was proved almost 50 years later by Fried
\cite{Fried:Schur}.
 
The obvious extension of Schur's original question to rational
functions over number fields then leads to the notion of {\em
  arithmetical exceptionality} as follows.

\begin{Definition}\label{exdef} Let $K$ be a number field, and $f\in
K(X)$ be a non--constant rational function. Write $f$ as a quotient of
two relatively prime polynomials in $K[X]$. For a place $\frakp$ of
$K$ denote by $K_\frakp$ the residue field. Except for finitely many
places $\frakp$, we can apply $\frakp$ to the coefficients of $f$ to
obtain $f_\frakp\in K_\frakp(X)$.

Regard $f_\frakp$ as a function from the projective line
$\PP(K_\frakp)=K_\frakp\cup\{\infty\}$ to itself via $x\mapsto
f_\frakp(x)$, with the usual rules to handle vanishing denominator or
$x=\infty$.

We say that $f$ is \emph{arithmetically exceptional} if there are
infinitely many places $\frakp$ of $K$ such that $f_\frakp$ is
bijective on $\PP(K_\frakp)$.\end{Definition}

A slight extension (see \cite{Turnwald:Schur}) of Fried's above result
is the following. Here $D_m(a,X)$ denotes the Dickson polynomial of
degree $m$ belonging to $a\in K$, which is most conveniently defined
by $D_m(a,Z+a/Z)=Z^m+(a/Z)^m$.

\begin{Theorem}[Fried]\label{polyex} Let $K$ be a number field, and
$f\in K[X]$ be an arithmetically exceptional polynomial. Then $f$ is
the composition of linear polynomials in $K[X]$ and Dickson
polynomials $D_m(a,X)$.\end{Theorem}

In this paper we generalize Schur's question to rational functions
over number fields. Under the assumption that the degree is a prime,
this has already been investigated by Fried \cite{Fried:CM}, showing
that there are non--polynomial examples over suitable number fields.

In contrast to the polynomial case, the classification result for
rational functions is quite complicated --- in terms of the associated
monodromy groups as well as from the arithmetic point of view. It is
immediate from the definition that if an arithmetically exceptional
rational function is a composition of rational functions over $K$,
then each of the composition factors is arithmetically exceptional.
Thus we are interested in indecomposable arithmetically exceptional
functions. Note however that the composition of arithmetically
exceptional functions need not be arithmetically exceptional, see
Theorem \ref{converse}.

In order to state the classification of arithmetically exceptional
functions, we need to introduce some notation. Let $K$ be a field of
characteristic $0$, and $f\in K(X)$ be a non--constant rational
function. Write $f=R/S$ with coprime polynomials $R,S\in K[X]$. Then
the degree $\deg(f)$ of $f$ is defined to be the maximum of $\deg(R)$
and $\deg(S)$. Note that $\deg(f)$ is also the degree of the field
extension $K(X)|K(f(X))$.

Let $t$ be a transcendental over $K$, and denote by $L$ a splitting
field of $R(X)-tS(X)$ over $K(t)$. Let $\hK$ be the algebraic closure
of $K$ in $L$, and set $A=\Gal(L|K(t))$, $G=\Gal(L|\hK(t))$.
Throughout the paper, $A$ and $G$ are called the \emph{arithmetic} and
\emph{geometric monodromy group} of $f$, respectively. Note that $G$
is a normal subgroup of $A$, with $A/G\cong\Gal(\hK|K)$. Furthermore,
$A$ and $G$ act transitively on the roots of $R(X)-tS(X)$.

Let $\bK$ be an algebraic closure of $K$, $\frakp_1,\dots,\frakp_r$ be
the places of $\bK(t)$ which ramify in $\bK L$, and let $e_i$ be the
ramification index of $\frakP_i|\frakp_i$, where $\frakP_i$ is a place
of $\bK L$ lying above $\frakp_i$. The unordered tuple
$(e_1,\dots,e_r)$ is a natural invariant associated to $f$, we call it
the \emph{ramification type} of $f$. Denote by $K_{\text{min}}$ a
minimal (with respect to inclusion) field, such that $f$ exists over
$K=K_{\text{min}}$ with given geometric monodromy group $G$ and given
ramification type. Of course the field $K_{\text{min}}$ need not be
unique, but it turns out that in most cases in the theorem below, it
is.

If we fix $G$ and the ramification type, then there are usually
several possibilities for $A$, and one can also ask for a minimal
field of definition if one fixes $A$ too. To keep the theorem below
reasonably short, we have not included that. Results can be found in
the sections where the corresponding cases are dealt with.

A rough distinction between the various classes of indecomposable
arithmetically exceptional functions is the genus $g$ of $L$. Note
that (see also Section \ref{AriPrep}, equation \eqref{RH})
\begin{equation}\label{gL}
2(\abs{G}-1+g)=\abs{G}\sum(1-1/e_i).
\end{equation}

The cases $g=0$ and $g=1$ belong to certain well--understood series,
to be investigated in Sections \ref{g=0} and \ref{g=1}, respectively.
If $g>1$, then $f$ belongs to one of finitely many possibilities
(finite in terms of $G$ and ramification type).  These sporadic cases
are investigated in Section \ref{g>1}.

\begin{Definition} Let $K$ be a field, and $f(X)\in K(X)$ be a
  rational function. Then we say that $f$ is \emph{equivalent} to
  $h(X)\in K(X)$, if $f$ and $h$ differ only by composition with
  linear fractional functions over $K$. (Of course, the arithmetic and
  geometric monodromy group is preserved under this equivalence, and
  so is the property of being arithmetically
  exceptional.)\end{Definition}

The rough classification of indecomposable arithmetically exceptional
rational functions is given in the following theorem. See Section
\ref{nota} for notation, and Sections \ref{g=0}, \ref{g=1}, and
\ref{g>1} for more details.

\begin{Theorem}\label{main} Let $K$ be a number field, $f\in K(X)$ be
  an indecomposable arithmetically exceptional function of degree
  $n>1$, $G$ the geometric monodromy group of $f$, $T=(e_1,\dots,e_r)$
  the ramification type, and $K_{\text{min}}$ be a minimal field of
  definition as defined above. Let $g$ be as defined in equation
  \eqref{gL}.\begin{itemize}
\item[(a)] If $g=0$, then either
\begin{itemize}
\item[(i)] $n\ge3$ is a prime, $G=C_n$, $T=(n,n)$,
$K_{\text{min}}=\QQQ$; or
\item[(ii)] $n\ge5$ is a prime, $G=C_n\rtimes C_2$, $T=(2,2,n)$,
$K_{\text{min}}=\QQQ$; or
\item[(iii)] $n=4$, $G=C_2\times C_2$, $T=(2,2,2)$,
$K_{\text{min}}=\QQQ$.
\end{itemize}
In case (i), $f$ is equivalent to $X^n$ or to a R\'edei function, in
case (ii), $f$ is equivalent to a Dickson polynomial (see Section
\ref{g=0}).
\item[(b)] If $g=1$, then $n=p$ or $p^2$ for an odd prime $p$. In the
former case set $N=C_p$, and $N=C_p\times C_p$ in the latter
case. Then one of the following holds.
\begin{itemize}
\item[(i)] $n\ge5$, $G=N\rtimes C_2$, $T=(2,2,2,2)$.
  $K_{\text{min}}=\QQQ$ if and only if $n=p^2$ or
  $n\in\{5,7,11,13,17,19,37,43,67,163\}$; or
\item[(ii)] $n\equiv1\pmod{6}$, $G=N\rtimes C_6$, $T=(2,3,6)$,
  $K_{\text{min}}=\QQQ$ if $n=p^2$, and
  $K_{\text{min}}=\QQQ(\sqrt{-3})$ if $n=p$; or
\item[(iii)] $n\equiv1\pmod{6}$, $G=N\rtimes C_3$, $T=(3,3,3)$,
  $K_{\text{min}}=\QQQ$ if $n=p^2$, and
  $K_{\text{min}}=\QQQ(\sqrt{-3})$ if $n=p$; or
\item[(iv)] $n\equiv1\pmod{4}$, $G=N\rtimes C_4$, $T=(2,4,4)$,
  $K_{\text{min}}=\QQQ$ if $n=p^2$, and
  $K_{\text{min}}=\QQQ(\sqrt{-1})$ if $n=p$.
\end{itemize}
The functions $f$ arise from isogenies or endomorphisms of elliptic
curves as described in Section \ref{g=1}.
\item[(c)] If $g>1$, then
\begin{itemize}
\item[(i)] $n=11^2$, $G=C_{11}^2\rtimes\GL_2(3)$, $T=(2,3,8)$,
  $g=122$, $K_{\text{min}}=\QQQ(\sqrt{-2})$; or
\item[(ii)] $n=5^2$, $G=C_5^2\rtimes S_3$, $T=(2,3,10)$, $g=6$,
  $K_{\text{min}}=\QQQ$; or
\item[(iii)] $n=5^2$, $G=C_5^2\rtimes H$ with $H$ a Sylow $2$-subgroup
  of the subgroup of index $2$ in $\GL_2(5)$, $T=(2,2,2,4)$, $g=51$,
  $K_{\text{min}}=\QQQ(\sqrt{-1})$; or
\item[(iv)] $n=5^2$, $G=C_5^2\rtimes(C_6\rtimes C_2)$, $T=(2,2,2,3)$,
  $g=26$, $K_{\text{min}}=\QQQ$; or
\item[(v)] $n=3^2$, $G=C_3^2\rtimes(C_4\rtimes C_2)$, $T=(2,2,2,4)$,
  $g=10$ or $T=(2,2,2,2,2)$, $g=19$, and $K_{\text{min}}=\QQQ$ in both
  cases; or
\item[(vi)] $n=28$, $G=\PSL_2(8)$, $T=(2,3,7)$, $g=7$, $T=(2,3,9)$,
  $g=15$, or $T=(2,2,2,3)$, $g=43$, and $K_{\text{min}}=\QQQ$ in all
  three cases; or
\item[(vii)] $n=45$, $G=\PSL_2(9)$, $T=(2,4,5)$, $g=10$, and
  $K_{\text{min}}=\QQQ$.
\end{itemize}
\end{itemize}
\end{Theorem}

The main part of the proof is group theoretic. Arithmetic
exceptionality translates equivalently to a property of finite
permutation groups, see Section \ref{AriPrep}. Let $A$ and $G$ be
finite groups acting transitively on $\Omega$; with $G$ a normal
subgroup of $A$. Then the triple $(A,G,\Omega)$ is said to be
\emph{exceptional} if the diagonal of $\Omega\times\Omega$ is the only
common orbit of $A$ and $G$. We do not completely classify such
triples -- we may return to this in a future paper. Instead, we use
further properties coming from the arithmetic context. Namely we get a
group $B$ with $G<B\le A$ with $B/G$ cyclic such that $(B,G,\Omega)$
is still exceptional. We say that the pair $(A,G)$ is
\emph{arithmetically exceptional} if such a group $B$ exists.

The notion of exceptionality has arisen in different context. It
typically comes in the following situation: Let $K$ be a field, and
$\phi:X\to Y$ be a separable branched cover of projective curves which
are absolutely irreducible over $K$. Let $A$ and $G$ denote the
arithmetic and geometric monodromy group of this cover, respectively.
Then $(A,G)$ is exceptional if and only if the diagonal is the only
absolutely irreducible component of the fiber product $X\times_\phi
X$. The arithmetic relevance lies in the well known fact that an
irreducible, but not absolutely irreducible curve, has only finitely
many rational points. Thus if the fiber product has the exceptionality
property, then there are only finitely many $P,Q\in X(K)$ with $P\ne
Q$ and $\phi(P)=\phi(Q)$, so $\phi$ is essentially injective on
$K$-rational points. If $K$ has a procyclic Galois group, for instance
if $K$ is a finite field, then $A/G$ is cyclic, and exceptionality
automatically implies arithmetic exceptionality.

For these reasons we have considered it worth the effort to obtain
reasonably complete reduction or classification results for
(arithmetically) exceptional permutation groups. For our original
question about the analog of Schur's question, we could have used the
genus $0$ condition on $G$ (see Section \ref{genusg}) at earlier
stages to remove many groups from consideration before studying
exceptionality.

A special case of exceptionality has arisen recently in some work on
graphs (see \cite{GLPS} and \cite{LiPr}). The question there involves
investigating those exceptional groups in which the quotient $A/G$ is
cyclic of prime order.

If $G\trianglelefteq A$ are transitive permutation groups on $n$
points, with $U<A$ the stabilizer of a point, then $A$ is primitive if
and only if $U$ is a maximal subgroup of $A$. Suppose that $A$ is
imprimitive, so there is a group $M$ with $U<M<A$. If $(A,G)$ is
arithmetically exceptional, then so are the pairs $(M,M\cap G)$ (in
the action on $M/U$) and $(A,G)$ in the action on $A/M$. The converse
need not hold (but it holds if $A/G$ is cyclic, see Lemma
\ref{excomp}).

Thus our main focus is on primitive permutation groups. We use the
Aschbacher--O'Nan--Scott classification of primitive groups, which
divides the primitive permutation groups into five classes (cf.\ 
Section \ref{PrimGroups}). The investigation of these five classes
requires quite different arguments. Exceptionality arises quite often
in affine permutation groups, one cannot expect to obtain a reasonable 
classification result.

Let $(A,G)$ be arithmetically exceptional, and $A$ primitive but not
affine. We show that $A$ does not have a regular nonabelian normal
subgroup, which completely removes one of the five types of the
Aschbacher--O'Nan--Scott classification. If $A$ preserves a product
structure (that means $A$ is a subgroup of a wreath product in product
action), then we can reduce to the almost simple case. If $A/G$ is
cyclic, does not preserve a product structure, but acts in diagonal
action, then we essentially classify the possibilities (up to the
question of existence of outer automorphisms of simple groups with a
certain technical condition). So the main case to investigate is the
almost simple action. We obtain the following classification, where
all the listed cases indeed do give examples. All these examples are
groups of Lie type. There are many possibilities if the Lie rank is
one. For higher Lie rank, the stablizers are suitable subfield
subgroups.

\begin{Theorem}\label{classificationexcept} Let $A$ be a primitive
  permutation group of almost simple type, so $L\trianglelefteq
  A\le\Aut(L)$ with $L$ a simple nonabelian group. Suppose there are
  subgroups $B$ and $G$ of $A$ with $G\trianglelefteq A$ and $B/G$
  cyclic, such that $(B,G)$ is exceptional. Let $M$ be a point
  stabilizer in $A$. Then $L$ is of Lie type, and one of the following
  holds.
\begin{enumerate}
\item[(a)] $L=\PSL_2(2^a)$, $G=L$ and either
\begin{enumerate}
  \item[(i)] $M\cap L=\PSL_2(2^{b})$ with $3 < a/b=r$ a prime and $B/G$
    generated by a field automorphism $x$ such that $\Cen_L(x)$ is not
    divisible by $r$; or
  \item[(ii)] $a>1$ is odd and $M\cap G$ is dihedral of order
    $2(q+1)$.
\end{enumerate}
\item[(b)] $L=\PSL_2(q)$ with $q=p^a$ odd, $G=L$ or $\PGL_2(q)$ and
  either
\begin{enumerate}
  \item[(i)] $M \cap L= \PSL_2(p^b)$ with $r:=a/b$ prime and $r$ not a
    divisor of $p(p^2-1)$, and $B$ any subgroup containing an
    automorphism $x$ which is either a field automorphism of order $e$
    or the product of a field automorphism of order $e$ and a diagonal
    automorphism, such that $r$ does not divide $p^{2a/e}-1$; or
  \item[(ii)] $M \cap L=L_2(p^{a/p})$ and $B$ is a subgroup of $\Aut(G)$
    such that $p|[B:G]$ and $B$ contains a diagonal automorphism; or
  \item[(iii)] $a$ is even, $M \cap L$ is the dihedral group of order
    $p^a-1$ and $B/G$ is generated by a field-diagonal automorphism; or
  \item[(iv)] $p=3$, $a$ is odd, $M \cap L$ is dihedral of order $p^a+1$
    and $B=\Aut(G)$.
\end{enumerate}
\item[(c)] $L=Sz(2^a)$ with $a>1$ odd, and $M \cap L = Sz(2^{a/b})$
  with $b$ a prime not $5$ and $a \ne b$ or $M \cap L$ is the
  normalizer of a Sylow $5$-subgroup (which is the normalizer of a
  torus); or
\item[(d)] $L=Re(3^a)$ with $a > 1$ odd and $M \cap L=Re(3^{a/b})$
  with $b$ a divisor of $a$ other than $3$ or $7$; or
\item[(e)] $L=U_3(p^a)$ and $L \cap M=U_3(p^{a/b})$ with $b$ a prime
  not dividing $p(p^2-1)$; or
\item[(f)] $L=U_3(2^a)$ with $a > 1$ odd and $M \cap L$ is the
  subgroup preserving a subspace decomposition into the direct sum of
  $3$ orthogonal nonsingular $1$-spaces; or
\item[(g)] $L$ has Lie rank $\ge2$, $L\ne\Sp_4(2)'$ ($\isom\PSL_2(9)$,
  a case covered in [b]).
  Then $M \cap L$ is a subfield group, the centralizer in $L$ of a
  field automorphism of odd prime order $r$.  Moreover,
  \begin{enumerate}
  \item[(i)]$r \ne p$,
  \item[(ii)] if $r =3$, then $L$ is of type $\Sp_4$ with $q$ even,
    and
  \item[(iii)] there are no $\Aut(L)$-stable $L$-conjugacy classes of
    $r$-elements.
\end{enumerate}
\end{enumerate}
\end{Theorem}

If we now consider the primitive arithmetically exceptional groups
which are not affine and add the genus $0$ condition on $G$, then only
two groups survive, accounting for the non--prime--power degrees $28$
and $45$ in Theorem \ref{main}.

In the affine case, we first use the genus $0$ condition (except for
noting that exceptionality implies that $G$ is not doubly transitive),
and then check which ones give arithmetically exceptional
configurations. We end up with six infinite series (and a few
sporadics). Two of the series are well known and correspond to R\'edei
functions and Dickson polynomials, respectively.  The other four
series have an interesting connection with elliptic curves, to be
investigated in Section \ref{g=1}. Section \ref{g>1} takes care of the
sporadic cases. There we will also see that there is a big variety of
arithmetically exceptional functions of degree $4$ over number fields.
This will have the following surprising consequence.

\begin{Theorem}\label{converse} Let $K$ be an arbitrary number
  field. There are four arithmetically exceptional rational functions
  of degree $4$ over $K$, such that their composition is not
  arithmetically exceptional.\end{Theorem}

\newpage

\tableofcontents

\newpage

\section{Arithmetic--geometric preparation}\label{AriPrep}

\subsection{Arithmetic and geometric monodromy groups}\label{AGMG}

Let $K$ be a field of characteristic $0$ and $t$ be a transcendental.
Fix a regular extension $E$ of $K(t)$ of degree $n$. (Regular means
that $K$ is algebraically closed in $E$.) Denote by $L$ a Galois
closure of $E|K(t)$. Then $A:=\Gal(L|K(t))$ is called the
\emph{arithmetic monodromy group} of $E|K(t)$, where we regard $A$ as
a permutation group acting transitively on the $n$ conjugates of a
primitive element of $E|K(t)$.

Denote by $\hat K$ the algebraic closure of $K$ in $L$. Then
$G:=\Gal(L|\hat K(t))$ is called the \emph{geometric monodromy group}
of $E|K(t)$. Note that $A/G\cong\Gal(\hat K|K)$, and that $G$ still
permutes the $n$ conjugates of $E$ transitively, as $E$ and $\hat
K(t)$ are linearly disjoint over $K(t)$.

If $f(X)\in K(X)$ is a non--constant rational function, then we use
the terms arithmetic and geometric monodromy group for the extension
$E=K(x)$ over $K(t)$, where $x$ is a solution of $f(X)=t$ in some
algebraic closure of $K(t)$. In this case, we call $\hK$ the
\emph{field of constants} of $f$.

\subsection{Distinguished conjugacy classes of inertia generators}

We continue to use the above notation. Let $\bK$ be an algebraic
closure of $K$. We identify the places $\frakp:t\mapsto b$ (or
$1/t\mapsto0$ if $b=\infty$) of $\bK(t)$ with $\bK\cup\{\infty\}$. Let
$\frakp$ be a ramified place of $L|K(t)$. Set $y:=t-\frakp$ (or
$y:=1/t$ if $\frakp=\infty$). There is a minimal integer $e$ such that
$L$ embeds into the power series field $\bK((y^{1/e}))$. For such an
embedding let $\sigma$ be the restriction to $L$ of the automorphism
of $\bK((y^{1/e}))$ which is the identity on the coefficients and maps
$y^{1/e}$ to $\zeta y^{1/e}$, where $\zeta$ is a primitive $e$th root
of unity. The embedding of $L$ is unique only up to a
$\bK(t)$--automorphism of $\bK L$, that is $\sigma$ is unique only up
to conjugacy in $G=\Gal(L|\hK(t))\cong\Gal(\bK L|\bK(t))$. We call the
$G$--conjugacy class $\cC$ of $\sigma$ the \emph{distinguished
  conjugacy class} associated to $\frakp$. Note that each $\sigma$ in
$\cC$ is the generator of an inertia group of a place of $L$ lying
above $\frakp$, and that $e$ is the ramification index of $\frakp$.

Let the tuple $\cB:=(\frakp_1,\dots,\frakp_r)$ consist of the places
of $\bK(t)$ which are ramified in $\bK L$, and let
$\cH=(\cC_1,\dots,\cC_r)$ be the tuple of the associated distinguished
conjugacy classes of $G$.  We call the pair $(\cB,\cH)$ the
\emph{ramification structure} of $L|K(t)$.

\subsection{Branch cycle descriptions}

A fundamental fact, following from Riemann's existence theorem and a
topological interpretation of the geometric monodromy group as
described below is the following. See \cite{H:Buch} for a
selfcontained proof.

\begin{Proposition}\label{RET1} With the notation from above, let
  $\frakp_1,\dots,\frakp_r$ be the places of $K(t)$ which ramify in
  $L$, and let $\cC_i$ be the distinguished conjugacy class associated
  to $\frakp_i$. Then there are $\s i\in\cC_i$, such that the $\s i$
  generate $G$, and the product relation $\s1\s2\dots\s r=1$
  holds.\end{Proposition}

We call such a tuple $(\s1,\dots,\s r)$ a \emph{branch cycle
description} of $L|K(t)$.

The branch cycle description allows us to compute the genus of the
field $E$. First note that this genus is the same as the genus of $\bK
E$. For $\sigma\in G$, let $\ind(\sigma)$ be $[E:K(t)]$ minus the
number of cycles of $\sigma$ in the given permutation representation.
Then the Riemann--Hurwitz genus formula yields
\begin{equation}\label{RH}
\sum_{i}\ind(\s i)=2([E:K(t)]-1+\text{genus}(E)).
\end{equation}

If the genus of $E$ from above is $0$, then we also call
$(\s1,\dots,\s r)$ a \emph{genus $0$ system}.

A sort of converse to Proposition \ref{RET1} is the following
algebraic version of Riemann's existence theorem, see \cite{H:Buch}.

\begin{Proposition}\label{RET2} Let $K$ be an algebraically closed
  field of characteristic $0$, $\s1,\dots,\s r$ be a generating system
  of a finite group $G$ with $\s1\s2\dots\s r=1$, and
  $\frakp_1,\dots,\frakp_r$ distinct elements from
  $K\cup\{\infty\}$.  Then there exists a Galois extension $L|K(t)$
  with group $G$, ramified only in $\frakp_1,\dots,\frakp_r$, such
  that the distinguished conjugacy class associated with $\frakp_i$ is
  the conjugacy class $\s i^G$.\end{Proposition}

\subsection{The branch cycle argument} We have seen that we can
realize any given ramification structure over $\bK(t)$. A situation we
will encounter frequently is whether we can get a regular extension
$E|K(t)$ with given ramification structure, and given pair $(G,A)$ of
geometric and arithmetic monodromy group. A powerful tool to either
exclude possible pairs, or to determine $A$ if $G$ is known, is
provided by Fried's \emph{branch cycle argument} (see
\cite[Lemma~2.8]{H:Buch}, \cite[Section~I.2.3]{MM},
\cite[Part~1]{Shih:PSL}), of which the following proposition is an
immmediate consequence.

\begin{Proposition}\label{BCA} Let $K$ be a field of characteristic $0$,
  $L|K(t)$ be a finite Galois extension with group $A$, $\hK$ the
  algebraic closure of $K$ in $L$, and $G:=\Gal(L|\hK(t))$. Let
  $(\cB,\cH)$ be the ramification structure of $L|\hK(t)$. Let
  $\zeta_n$ be a primitive $n$--th root of unity, where $n=\abs{G}$.
  Let $\gamma\in\Gal(\bK|K)=:\Gamma$ be arbitrary, and $m$ with with
  $\gamma^{-1}(\zeta_n)=\zeta_n^m$. Let $a\in A$ be the restriction to
  $L$ of an extension of $\gamma$ to $\Gal(\bK L|K(t))$. Then $\cB$ is
  $\Gamma$--invariant and $\cC_{\gamma(\frakp)}=(\cC_\frakp^{m})^a$
  for each $\frakp\in\cB$.
\end{Proposition}

Often, we will use the following special case.

\begin{Corollary}\label{BCAC}
  Let $L|\QQQ(t)$ be a finite Galois extension with group $A$, $\hQ$
  the algebraic closure of $\QQQ$ in $L$, and $G:=\Gal(L|\hQ(t))$. Let
  $\frakp$ be a rational place of $\QQQ(t)$, and $\sigma$ be in the
  associated distinguished conjugacy class of $G$. Then for each $m$
  prime to the order of $\sigma$, the element $\sigma^m$ is conjugate
  in $A$ to $\sigma$.
\end{Corollary}

Often the branch cycle argument does not help in excluding certain
pairs $(A,G)$ and given ramification structure. Then the following
observation is sometimes helpful.

\begin{Lemma}\label{AG} Let $K$ be a field of characteristic $0$, $t$
  a transcendental, and $L|K(t)$ be a finite Galois extension with
  group $A$. Let $\hK$ be the algebraic closure of $K$ in $L$, and set
  $G=\Gal(L|\hK(t))$. Let $\frakp$ be a rational place of $K(t)$, and
  $\frakP$ be a place of $L$ lying above $\frakp$. Denote by $D$ and
  $I$ the decomposition and inertia group of $\frakP$, respectively.
  Then $I\le D\cap G$ and $A=GD$. In particular,
  $A=G\Nor_A(I)$.\end{Lemma}

\begin{proof} Without loss assume that $\frakp$ is a finite place, so
  $K[t]$ is contained in the valuation ring. Let $\ol$ be the integral
  closure of $K[t]$ in $L$. Then $D/I$ is the Galois group of
  $(\ol/\frakP)|(K[t]/\frakp)$, see \cite[Chapter I, \S7,
  Prop.~20]{Serre:LF}. On the other hand, $\hK$ embeds into
  $\ol/\frakP$, so $D/I$ surjects naturally to $A/G=\Gal(\hK|K)$.
  Furthermore, if $\phi\in I$, then $u-u^\phi\in\frakP$ for all
  $u\in\hK$, hence $\phi$ is trivial on $\hK$, so $I\le D\cap
  G$.\end{proof}

\subsection{Weak rigidity}

We use the weak rigidity criterion in order to prove existence of
certain regular extensions $E|K(t)$ with given ramification structure
and geometric monodromy group. The main reference of the material in
this section is \cite{MM} and \cite{H:Buch}.

\begin{Definition}\label{rigid} Let $G$ be a finite group, and
  $\cH=(\cC_1,\cC_2,\dots,\cC_r)$ be an $r$--tuple of conjugacy
  classes of $G$. Consider the set of $r$--tuples $(\s1,\s2,\dots,\s
  r)$ with $\s i\in\cC_i$ which generate $G$ and $\s1\s2\dots\s r=1$.
  We say that $\cH$ is \emph{weakly rigid} if this set is not empty,
  and if any two $r$--tuples in this set are conjugate under
  $\Aut(G)$.\end{Definition}

\begin{Definition} Let $K$ be a subfield of $\CCC$, and
  $\cB\subset\bK\cup\{\infty\}$ be an unordered tuple of distinct
  elements.  Associate to each $\frakp\in\cB$ a conjugacy class
  $\cC_\frakp$ of a finite group $G$ of order $n$, and let $\cH$ be
  the tuple of these classes. Set $\Gamma:=\Gal(\bK|K)$, and let
  $\zeta_n$ be a primitive $n$--th root of unity. For each
  $\gamma\in\Gamma$ let $m(\gamma)$ be an integer with
  $\gamma^{-1}(\zeta_n)=\zeta_n^{m(\gamma)}$.
\begin{itemize}
\item[(a)] The pair $(\cB,\cH)$ is called $K$--\emph{rational}, if
  $\cB$ is $\Gamma$--invariant and
  $\cC_{\gamma(\frakp)}=\cC_\frakp^{m(\gamma)}$ for each
  $\gamma\in\Gamma$ and $\frakp\in\cB$.
\item[(b)] The pair $(\cB,\cH)$ is called \emph{weakly} $K$--{\em
    rational}, if $\cB$ is $\Gamma$--invariant and
  $\cC_{\gamma(\frakp)}=\alpha(\cC_\frakp^{m(\gamma)})$ for each
  $\gamma\in\Gamma$ and $\frakp\in\cB$, where $\alpha\in\Aut(G)$
  depends on $\gamma$, but not on $\frakp$.
\end{itemize}
\end{Definition}

\begin{Proposition}\label{WeakRigCrit} Let $K$ be a field of
characteristic $0$, $G$ be a finite group which is generated by
$\s1,\dots,\s r$ with $\s1\s2\dots\s r=1$. Let
$\cB=\{\frakp_1,\frakp_2,\dots,\frakp_r\}$ $\subset\PP(\bK)$, and let
$\cC_{\frakp_i}$ be the conjugacy class of $\s i$. Suppose that
\[
\cH:=(\cC_{\frakp_1},\dots,\cC_{\frakp_r})
\]
is weakly rigid and
\[
(\cB,\cH)
\]
is weakly $K$--rational.

Furthermore, let $H<G$ be a subgroup which is self--normalizing and
whose $G$--conjugacy class is $\Aut(G)$--invariant. Suppose that one
of the $\frakp_i$ is rational, and that $\cC_{\frakp_i}$ intersects
$H$ non--trivially for this index $i$. Also, suppose that the $\s i$
are a genus $0$ system with respect to the action on $G/H$.

Then there exists $f\in K(X)$ with geometric monodromy group $G$
acting on $G/H$ and ramification structure $(\cB,\cH)$.
\end{Proposition}

\begin{proof} By Proposition \ref{RET2}, there exists a Galois
  extension $L|\bK(t)$ with group $G$ and ramification structure
  $(\cB,\cH)$. Weak rigidity implies that $L$ is Galois over $K(t)$,
  see e.g.~\cite[Remark 3.9(a)]{H:Buch}. Set $A:=\Gal(L|K(t))$, so $G$
  is the fixed group of $\bK(t)$. Let $J$ be the normalizer of $H$ in
  $A$. As $G$ is transitive on the $A$--conjugates of $H$ we obtain
  $A=GJ$. Also, $G\cap J=H$ because $H$ is self--normalizing in $G$.
  Let $E$ be the fixed field of $J$ in $L$. Then $E$ is linearly
  disjoint from $\bK(t)$ over $K(t)$, and $\bK E$ is the fixed field
  of $H$. By assumption, $E$ has genus $0$, and the hypothesis about
  the nontrivial intersection of some class $\cC_{\frakp_i}$ with $H$
  with $\frakp_i$ being $K$--rational implies $E$ has a $K$--rational
  place, thus is a rational field $K(x)$, e.g.\ by
  \cite[I.6.3]{Stich}. Write $t=f(x)$ with $f\in K(X)$, then $f$ is
  the desired function.\end{proof}

\begin{Definition} Let $\cC$ be a conjugacy class of $G$, $n$ the
  order of $G$, and $\zeta_n$ a primitive root of unity. For $\gamma$
  an automorphism of $\bQ$ let $m$ be with $\gamma(\zeta_n)=\zeta^m$.
  Then the class $\cC$ is said to be $K$--rational if $\cC^m=\cC$ for
  each $\gamma\in\Gal(\bQ|K)$. This property is well known to be
  equivalent to $\chi(\cC)\in K$ for all irreducible characters $\chi$
  of $G$. In particular, if the character values of the classes in
  $\cH$ as above are in $K$, and $\cB\subset K$, then $(\cB,\cH)$ is
  $K$--rational.\end{Definition}

\subsection{Topological interpretation}\label{genusg}
Though we do not really need it here, it might help understanding some
arguments in this paper if one also has the geometric interpretation
of a branch cycle description in mind.
 
Suppose that $K$ is a subfield of the complex numbers $\CCC$. As
$\CCC(t)\cap L=\hat K(t)$ (see \cite[Corollary 2, V, \S4]{Chevalley}),
we get $\Gal(\CCC L|\CCC(t))\cong\Gal(L|\hat K(t))$ by restriction.
For any holomorphic covering $\alpha:A\to B$ of Riemann surfaces,
denote by $\alpha^\star:\cM(B)\hookrightarrow\cM(A)$ the natural
inclusion of the fields of meromorphic functions on $B$ and $A$.
Associated to $\CCC E|\CCC(t)$ is a branched, holomorphic covering
$\pi:S\to\PP(\CCC)$ of degree $n=[E:K(t)]$ with $S$ a connected
Riemann surface and $\PP(\CCC)$ the Riemann sphere, such that the
extension $\cM(S)|\pi^\star(\cM(\PP(\CCC)))$ can be identified with
$\CCC E|\CCC(t)$. Observe that $\cM(\PP(\CCC))\cong\CCC(t)$.

Let $\cB=\{\frakp_1,\frakp_2,\ldots,\frakp_r\}$ be the set of branch
points of $\pi$.  Fix $\frakp_0\in\PP(\CCC)\setminus\cB$, and denote
by $\cG$ the fundamental group
$\pi_1(\PP(\CCC)\setminus\cB,\frakp_0)$. Then $\cG$ acts transitively
on the points of the fiber $\pi^{-1}(\frakp_0)$ by lifting of paths.
Fix a numbering $1,2,\ldots,n$ of this fiber. Thus we get a
homomorphism $\cG\to S_n$ of $\cG$ into the symmetric group $S_n$. By
standard arguments, the image of $\cG$ can be identified with the
geometric monodromy group $G$ defined above, thus we write $G$ for
this group too.

This identification has a combinatorial consequence. Choose a standard
homotopy basis of $\PP(\CCC)\setminus\cB$ as follows. Let $\gamma_i$
be represented by paths which wind once around $\frakp_i$ clockwise,
and around no other branch point, such that
$\gamma_1\gamma_2\cdots\gamma_r=1$. Then
$\gamma_1,\gamma_2,\ldots,\gamma_{r-1}$ freely generate $\cG$.

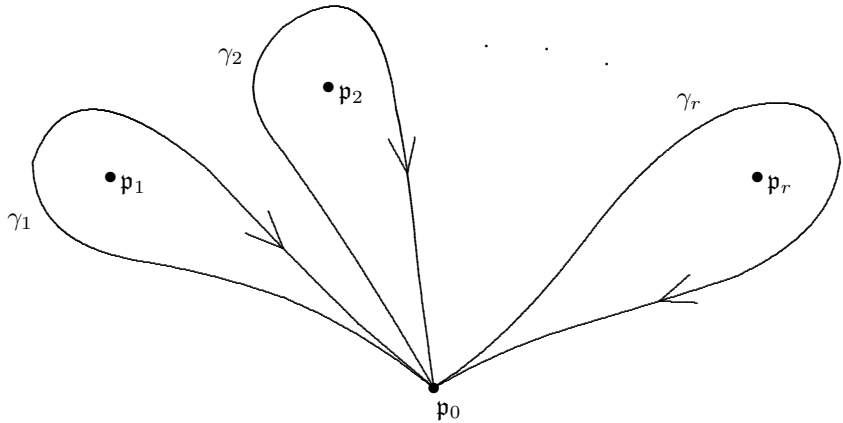
\setlength{\unitlength}{1mm}
\begin{picture}(130,70)

\curve(65,10,55,17,45,22)
\curve(45,22,35,25,25,27)
\curve(25,27,15,31.5,11.7,40)
\curve(11.7,40,20,47,35,39)
\curve(35,39,45,28.5,55,18.5)
\curve(55,18.5,65,10)

\curve(43,33.5,45,28.5)
\curve(45,28.5,40,30.5)

\curve(65,10,55,26.5,45,41.3)
\curve(45,41.3,41,50,45,58)
\curve(45,58,55,59.7,60,47)
\curve(60,47,61.5,38.5,63,25)
\curve(63,25,65,10)

\curve(59,43,61.5,38.5)
\curve(61.5,38.5,62.5,43.3)

\curve(65,10,75,18,85,29.5)
\curve(85,29.5,95,40.7,105,47)
\curve(105,47,115,47,119,40)
\curve(119,40,115,31.5,105,24.7)
\curve(105,24.7,95,21.5,85,18.5)
\curve(85,18.5,75,15,65,10)

\curve(99,25,95,21.5)
\curve(95,21.5,100,21.3)

\put(65,10){\circle*{1.5}}
\put(67,7){\makebox(0,0){$\frakp_0$}}

\put(22,38){\circle*{1.5}}
\put(25,37){\makebox(0,0){$\frakp_1$}}
\put(51,50){\circle*{1.5}}
\put(54,49){\makebox(0,0){$\frakp_2$}}
\put(108,38){\circle*{1.5}}
\put(111,37){\makebox(0,0){$\frakp_r$}}

\put(72,55.5){\circle*{.5}}
\put(80,55){\circle*{.5}}
\put(88,53){\circle*{.5}}

\put(10,32){\makebox(0,0){$\gamma_1$}}
\put(38,54){\makebox(0,0){$\gamma_2$}}
\put(99,48){\makebox(0,0){$\gamma_r$}}
\end{picture}

Let $\s i$ be the image of $\gamma_i$ in $G$. It can be shown that
this construction is compatible with the algebraic definition of the
conjugacy classes $\cC_i$ (see \cite[Section I.5.4]{H:Buch}), thus
proving Proposition \ref{RET1}.

\subsection{Group theoretic translation of arithmetic exceptionality}
Here we provide the setup for the equivalent translation of arithmetic
exceptionality in terms of monodromy groups. The material in this
section is well known (mainly in the context of polynomials, with
immediate generalization to rational functions). Let $K$ be a number
field, and $f(X)\in K(X)$ of positive degree.

With $t$ a transcendental, let $L$ be a splitting field of $f(X)-t$
over $K(t)$, and $\hat K$ be the algebraic closure of $K$ in $L$.  Set
$A:=\Gal(L|K(t))$ --- the arithmetic monodromy group of $f$ --- and
$G:=\Gal(L|\hat K(t))$ --- the geometric monodromy group. Recall that
$A/G\cong\Gal(\hat K|K)$. We regard $A$ and its subgroups as
permutation groups on the set $\Omega$ of the roots of $f(X)-t$. If
$B$ is a group between $A$ and $G$, then we say that the pair $(B,G)$
is \emph{exceptional} if the diagonal of $\Omega\times\Omega$ is the
only common orbit of $B$ and $G$.

\begin{Theorem}[Fried, \cite{Fried:CM}]\label{togroups} In the setup
from above the following holds.\begin{itemize}
\item[(a)] $f$ is arithmetically exceptional if and only if there is a
  group $B$ with $G\le B\le A$ with $B/G$ cyclic such that $(B,G)$ is
  exceptional.
\item[(b)] There is a constant $C\in\NNN$ depending on $f$ with the
  following property. Let $\frakp$ be a place of $K$ with
  $\abs{K_\frakp}>C$. Then $f_\frakp$ is defined, and $f_\frakp$ is
  bijective on $\PP(K_\frakp)$ if and only if $(B,G)$ is exceptional,
  where $B/G$ is the decomposition group of a place of $\hat K$ lying
  above $\frakp$.\end{itemize}\end{Theorem}

\begin{proof} There is a constant $C'$ such that for
  $\abs{K_\frakp}>C'$ the following holds. The place $\frakp$ is
  unramified in $\hat K$, $f_\frakp$ is defined, and $f_\frakp(X)-t$
  is separable. Let $A_\frakp$ and $G_\frakp$ be the arithmetic and
  geometric monodromy group of $f_\frakp$ over $K_\frakp$
  respectively.
  
  Let $B$ be the subgroup of $A$ such that $B/G$ induces the
  decomposition group of a place of $\hat K$ lying above $\frakp$.
  Then, if $C'$ was big enough, we have the following result by Fried
  \cite[Lemma 1]{Fried:HIT}:\begin{align*}A_\frakp &\cong B\\ G_\frakp
    &\cong G.\end{align*}
  
  Furthermore, it is known (see \cite{FGS}) that there is a bound
  $C''$ such that if $\abs{K_\frakp}>C''$, then $f_\frakp$ is
  bijective on $\PP(K_\frakp)$ if and only if $(A_\frakp,G_\frakp)$ is
  exceptional.  This proves (b) and the `only if' part of (a).
  
  To prove the `if' part of (a), assume that there is $B$ as in (a).
  Let $bG$ generate $A/G$. By Chebotar\"ev's density theorem there are
  infinitely many primes $\frakp$ of $K$ with a prime $\frakP$ of
  $\hat K$ above, such that $bG$ induces the Frobenius automorphism on
  $\hat K_\frakP|K_\frakp$.  In particular, $B$ then is the
  decomposition group of $\frakP$, and the assertion follows from
  (b).\end{proof}

\subsection{Remark about exceptional functions over finite fields}

Let $\FFF_q$ be a finite field with $q$ elements, where $q$ is a power
of the prime $p$. Let $f(X)\in\FFF_q(X)$ be a rational function which
is not the $p$--th power of another rational function. We say that $f$
is \emph{exceptional} if $f$ is bijective on $\PP(\FFF_{q^m})$ for
infinitely many $m\in\NNN$. Let $A$ and $G$ be the arithmetic and
geometric monodromy group of $f$. Then $A/G$ is cyclic, and $f$ is
exceptional if and only if $(A,G)$ is exceptional with respect to the
action on the roots of $f(X)-t$, this follows from proof of Theorem
\ref{togroups}.

In \cite{FGS} the possible monodromy groups of exceptional polynomials
have been classified, and the question about existence of actual
examples has been answered positively for all non--affine groups in
\cite{CohenMatthews}, \cite{LZ}, and for some affine groups in
\cite{GM:ExPol}.  

Indeed, aside from the affine examples, all indecomposable exceptional
polynomials are now known (see \cite{GZ} and \cite{GRZ}). This can be
extended to arbitrary fields using \cite{GSRam}.

We hope to eventually achieve a similar result for exceptional
rational functions over finite fields.  Our result Theorem
\ref{classificationexcept} shows that the key ingredient is to prove
that the generic example of exceptionality cannot be of genus $0$
(i.e.\ the case where $G$ is a Chevalley group and the action is on
the cosets of a subfield group).

If all the ramification of
$\overline{\FFF_q}(X)|\overline{\FFF_q}(f(X))$ is tame (that is all
the inertia groups have order not divisible by $p$), then we get the
exact analogue of Proposition \ref{RET1} by Grothendieck's existence
theorem, see \cite[XIII, Corollary 2.12]{SGA1}. As a result, we obtain
the same group--theoretic candidates. Of course, the question about
actual existence, which used various arithmetic arguments and results
applicable for the rationals or number fields, is quite different and
will not be pursued in this paper.

\section{Group theoretic exceptionality}\label{groupex}
\newcommand{\GF}{\rm GF}
\newcommand{\ep}{\mathop{\epsilon}}

\subsection{Notation and definitions}\label{nota}
Here we gather definitions and notation from finite group theory which
are used throughout this paper. Notation used only locally is defined
where it first appears.
\begin{description}
\item[General notation:] For $a,b$ elements of a group $H$ set
$a^b:=b^{-1}ab$. Furthermore, if $A$ and $B$ are subsets of $H$, then
$A^b$, $a^B$ and $A^B$ have their obvious meaning. If $G$ is a
subgroup of $H$, then for a subset $S$ of $H$ let $\Cen_G(S)$ denote
the centralizer of $S$ in $G$ and $\Nor_G(S)$ denote the normalizer
$\{g\in G|\;S^g=S\}$ of $S$ in $G$.

If $A,B,\dots$ is a collection of subsets or elements of $H$, then we
denote by $\gen{A,B,\dots}$ the group generated by these sets and
elements.

Let $\pi$ be a set of primes. A $\pi$--group is a group with order
only divisible by primes in $\pi$. Denote $\PPP\setminus\pi$ by
$\pi'$. The subgroup $\O{\pi}(H)$ of $H$ is the maximal normal
$\pi$--subgroup. If $\pi$ consists of a single prime $p$, then we
write $\Op$ rather than $\mathbf{O}_{\{p\}}$.

\item[Permutation groups:] Let $G$ be a permutation group on a finite
  set $\Omega$. We use the exponential notation $\omega^g$ to denote
  the image of $\omega\in\Omega$ under $g\in G$. The stabilizer of
  $\omega$ in $G$ is denoted by $G_\omega$.  If $G$ is transitive on
  $\Omega$, then the action of $G$ on $\Omega$ is isomorphic to the
  coset action of $G$ on $G/G_{\omega}$, the set of right cosets of
  $G_{\omega}$ in $G$.
  
  Let $\chi(g)$ be the number of fixed points of $g$ on $\Omega$, and
  $\ind(g)$ the index, which is $\abs{\Omega}$ minus the number of
  orbits of $\gen{\omega}$, see Section \ref{genusg}.

\item[Specific groups:] We denote by $C_n$ and $D_n$ the cyclic and
dihedral group of order $n$. The alternating and symmetric group on
$n$ letters is denoted by $A_n$ and $S_n$ respectively (also if we
view it as an abstract rather than permutation group). We will
generally follow the notation used in \cite{ATLAS} for the almost
simple groups.

\item[Exceptionality:] Let $G$ be a normal subgroup of $A$, with both
$A$ and $G$ being transitive permutation groups on the finite set
$\Omega$. Then $(A,G,\Omega)$ is called \emph{exceptional} if the
diagonal of $\Omega\times\Omega$ is the only common orbit of $A$ and
$G$ on $\Omega\times\Omega$. This definition is equivalent to the
following: Let $\omega\in\Omega$, then every $A_\omega$--orbit except
$\{\omega\}$ breaks up into strictly smaller $G_\omega$--orbits.

Motivated by Theorem \ref{togroups}, we call the triple $(A,G,\Omega)$
\emph{arithmetically exceptional}, if there is a subgroup $B$ of $A$
which contains $G$, such that $(B,G,\Omega)$ is exceptional, and $B/G$
is cyclic.
\end{description}

\subsection{Primitive groups}\label{PrimGroups}

Recall that if $J$ is a finite group, then $F(J)$ is the Fitting
subgroup of $J$, the maximal normal nilpotent subgroup of $J$.  We
denote the maximal normal $p$-subgroup of $J$ by $\Op(J)$.  A subgroup
$H$ is subnormal (written $H\sn J$) if there is a composition series
for $J$ with $H$ as a term.  A group is called quasisimple if it is
perfect (i.e.\ its own commutator subgroup) and simple modulo its
center.  A component of $J$ is a quasisimple subnormal subgroup.  Any
two distinct components commute and any component commutes with
$F(J)$.  The central product of all the components is denoted by
$E(J)$.  The generalized Fitting subgroup is $\gFit(J):=E(J)F(J)$.  An
important property is that $\gFit(J)$ contains its own centralizer
(and so $J/Z(\gFit(J))$ embeds in $\Aut(\gFit(J))$).

We say $J$ is almost simple if $\gFit(J)$ is a nonabelian simple
subgroup.  Note that this is equivalent to $L\subseteq
J\subseteq\Aut(L)$ for $L$ nonabelian simple (and so
$L=\gFit(J)$). See \cite{Asch:Book}, \cite{Huppert1}, or
\cite{Suzuki1} for more details on these concepts.

We use the Aschbacher--O'Nan--Scott classification of primitive groups
(cf.\ \cite{AschScott}, \cite{LPS}). Let $A$ be a group acting
primitively and faithfully on $\Omega$. Let $M$ be a point stabilizer
and set $E=\gFit(A)$. We have to consider the following cases (there
are other ways of dividing up the possibilities):\begin{description}
\item[(Aff)] $E$ is an elementary abelian $p$-group, $A=ME$
(semidirect), $M$ acts irreducibly on $E$ and $E$ acts regularly on
$A/M$.  \item[(Reg)] $A$ contains a regular normal nonabelian
subgroup.  \item[(AS)] $A$ is almost simple.  \item[(Prod)]
$E=L_1\times\dots\times L_t$ with $t>1$ and $L_i\cong L$ a nonabelian
simple group. Moreover, the $L_i$ are conjugate in $A$ and $E\cap
M=J_1\times\dots\times J_t$ where $J_i=M\cap L_i$.  \item[(Diag)]
$E=L_1\times\dots\times L_t$ with $t>1$ and $L_i\cong L$ a nonabelian
simple group. Moreover, the $L_i$ are conjugate in $A$ and $E\cap M$
is a product of diagonal subgroups of $E$.  \end{description}

\subsection{General results on exceptionality}\label{genex}

We begin with some general properties of exceptionality which do not
require a finer analysis of various classes of groups.

The following result is straightforward and well known.

\begin{Lemma}\label{chig} Let $G$ act transitively on a finite set
$\Omega$. If $H$ is a point stabilizer and $g\in G$, then
$$\chi(g)=\sum_{i=1}^m[\Cen_G(g_i):\Cen_H(g_i)],$$ where
$g_1,\dots,g_m$ is a set of representatives for the $H$-conjugacy
classes of $g^G\cap H$.\end{Lemma}

\begin{proof}$$\chi(g)=\abs{\Cen_G(g)}\abs{g^G\cap H}/\abs{H}=
\abs{\Cen_G(g)}\sum\abs{g_i^H}/\abs{H}=\sum[\Cen_G(g_i):\Cen_H(g_i)].$$
\end{proof}

The next two results are essentially in \cite[13.1]{FGS}. See also
\cite[Section 3]{GW}.

\begin{Lemma}\label{chicount} Let $A$ be a finite group with $G$ a
normal subgroup with $A/G$ cyclic generated by $xG$. Suppose $A$ acts
on a finite set $\Omega$. Then $\sum_{g\in G}\chi(xg)=r\abs{G}$, where
$r$ is the number of common orbits that $A$ and $G$ have on
$\Omega$.\end{Lemma}

In particular, if $A$ and $G$ are both transitive on $\Omega$,
then the sum above is $\abs{G}$, i.e.\ the average number of fixed
points of an element is $1$. An immediate consequence is:

\begin{Lemma}\label{chiex} Let $G$ be a normal subgroup of
the finite group $A$ with $A/G$ cyclic generated by $xG$. Let $\Omega$
be a transitive $A$-set. Assume that $G$ is also transitive on
$\Omega$. The following are equivalent:\begin{itemize}
\item[(i)] $(A,G,\Omega)$ is exceptional;
\item[(ii)] $\chi(xg)\le1$ for all $g\in G$;
\item[(iii)] $\chi(xg)=1$ for all $g\in G$;
\item[(iv)] $\chi(xg)\ge1$ for all $g\in G$.
\end{itemize}\end{Lemma}

Let $M$ be a subgroup of $A$. Note that $G$ is transitive on $A/M$ if
and only if $A=MG$. If in addition $A/G$ is cyclic then this is
equivalent to the existence of $x\in M$ with $A/G$ generated by
$xG$.

\begin{Lemma}\label{Ki} Assume that $(A,G,A/M)$ is exceptional and
that $xG$ is a generator of $A/G$ with $x\in M$.\begin{itemize}
\item[(a)] $\Cen_A(x)$ and $\Nor_A(\gen{x})$ have a unique fixed
point on $A/M$ (and therefore are contained in $M$);
\item[(b)] If $g\in A$ and $gxg^{-1}\in M$, then $g\in M$;
\item[(c)] Let $M_0=\Nor_A(\gen{x})$ and define $M_{i+1}$ inductively
by $M_{i+1}=\{g\in A:gxg^{-1}\in M_i\}$. Then $M_i\subseteq M$ for
each $i$.\end{itemize}\end{Lemma}

\begin{proof} By Lemma \ref{chiex} and exceptionality, $x$ has a
unique fixed point, whence (a) holds. Let $g\in A$ with $gxg^{-1}
\in M$. So $g$ fixes the unique fixed point of $x$, thus $g\in M$.
Claim (c) follows from (b).\end{proof}

The following is well known and easy.

\begin{Lemma}\label{excomp} Let $G$ be a normal subgroup of $A$ with
$A/G$ cyclic. Suppose $M<U<A$ and $A=GM$. Then $(A,G,A/M)$ is
exceptional if and only if $(A,G,A/U)$ and $(U,G\cap U,U/M)$ are
exceptional.\end{Lemma}

\begin{proof} Choose a generator $xG$ for $A/G$ with $x\in M$.
Suppose $(A,G,A/U)$ is not exceptional. Then some
element $xg$ has no fixed points on $A/U$ and so has none
on $A/M$, whence $(A,G,A/M)$ is not exceptional.

So we may assume that $(A,G,A/U)$ is exceptional. Then
each $xg$ has a unique fixed point on $A/U$. So $xg$ is conjugate to
some $y\in U$. The number of fixed points of $y$ on $A/M$ is
precisely the number of fixed points on $U/M$; the claim
follows.\end{proof}

The next results provide methods for producing exceptional triples.

Suppose $A/G$ is cyclic of order $e$. Note that we may always choose a
generator $xG$ for $A/G$ such that the order of $x$ divides a power of
$e$ (by replacing $x$ by a power of itself) and necessarily that $e$
divides the order of $x$.

Before stating the next theorem, we make the following easy observation
-- if $G$ is a finite group and $x$ is an automorphism of $G$, then
$\Cen_G(xg)$ is a $\pi'$-group for  all $g \in G$ if and only if
$x$ preserves no $G$-conjugacy class of $\pi$-elements in $G$.

\begin{Theorem}\label{expi} Let $G$ be a normal subgroup of the finite
group $A$ with $A/G$ generated by $xG$. Let $\pi$ be a set of primes,
and write $x=yz=zy$ where $y$ is a $\pi$--element and $z$ is a
$\pi'$--element. Assume that $\Cen_G(xg)$ is a $\pi'$--group for every
$g\in\Cen_G(y)$. Then $(A,G,A/\Cen_A(y))$ is exceptional.\end{Theorem}

\begin{proof} Set $M=\Cen_A(y)$. Note that $G$ is transitive on $A/M$,
since $x\in M$. It suffices to prove that $xg$ has at most $1$ fixed
point on $A/M$ for each $g\in G$. We may assume that $xg$ has at least
one fixed point and so we may take $g\in M$ (note that $x\in M$).
Write $xg=y'z'=z'y'$ where $y'$ is a $\pi$--element and $z'$ is a
$\pi'$--element in $M$. Note that $y\in\Cen_A(xg)$ (since $g\in M$)
and $yG=y'G$ (as the decomposition into commuting $\pi$-- and
$\pi'$--elements is unique). But the $\pi$--element $y{y'}^{-1}$
commutes with $xg$, so is also a $\pi'$--element by assumption, thus
$y=y'$. We contend that any $A$--conjugate of $xg$ which is contained
in $M$ is an $M$--conjugate. So assume that $(xg)^a\in M$, where we
may assume $a\in G$ as $A=\Cen_A(xg)G$. So $(xg)^a=x\tilde g$ with
$\tilde g=a^{-x}ga\in\Cen_G(y)$. As above, it follows that the
$\pi$--part of $x\tilde g$ is $y$. Thus $a\in\Cen_A(y)=M$, proving our
claim.\end{proof}

\begin{Corollary}\label{excor} Let $G$ be a normal subgroup of the
finite group $A$ with $A/G$ generated by $xG$. Let $H$ be a subgroup
of $G$ normalized by $x$. Let $M=\gen{H,x}$. If the order of $x$ is
relatively prime to $\abs{\Cen_G(x)}$, then $(A,G,A/M)$ is exceptional
if and only if $\Cen_G(x)\subseteq M$.\end{Corollary}

\begin{proof} The condition $\Cen_G(x)\subseteq M$ is necessary by
Lemma \ref{Ki}(a).

We show that even $(A,G,A/\Cen_A(x))$ is exceptional, the claim then
follows from Lemma \ref{excomp}. Since the order of $x$ is prime to
$\abs{\Cen_G(x)}$, it follows that $\Cen_G(xg)\subseteq\Cen_G(x)$
for every $g\in\Cen_G(x)$:
if $c^{-1}xgc = xg$ then $c^{-1}xc = x$ by the uniqueness of the
decomposition of $xg$ into commuting $\pi$-- and $\pi'$--elements.
Thus, we may apply the previous result
with $y=x$.\end{proof}

The next two lemmas show that exceptionality descends in certain cases
and will allow us to use induction in analyzing this property.

\begin{Lemma}\label{desc1} Let $G$ be a normal subgroup of $A$ with
$A/G$ generated by $xG$. Let $M$ be a subgroup of $A$ with $(A,G,A/M)$
exceptional. If $A_1$ is a subgroup of $A$ such that $x\in A_1\cap M$,
then $(A_1,G\cap A_1,A_1/(M\cap A_1))$ is exceptional.\end{Lemma}

\begin{proof} Note that $A=GA_1$ as $x\in A_1$. We get
$A_1=G\gen{x}\cap A_1=(G\cap A_1)\gen{x}\subseteq(G\cap A_1)(A_1\cap
M)$, so $G\cap A_1$ is transitive on $A_1/(M\cap A_1)$. Since $xg\in
xG$ for every $g\in G\cap A_1$, $xg$ has at most one fixed point on
$A_1/(M\cap A_1)$. By Lemma \ref{chiex} this implies that $(A_1,G\cap
A_1,A_1/(M\cap A_1))$ is exceptional.\end{proof}

\begin{Lemma}\label{desc2} Let $G$ be a normal subgroup of $A$ with
$A/G$ generated by $xG$ with $\Omega$ a transitive $A$-set. Let $L$ be
a nontrivial normal subgroup of $A$ contained in $G$ such that $L$ is
transitive on $\Omega$ (e.g.\ this is so if $A$ is primitive on
$\Omega$).  Then $(A,G,\Omega)$ exceptional implies that
$(B,L,\Omega)$ is exceptional where $B=\gen{L,x}$.\end{Lemma}

\begin{proof} Every element in the coset $xL$ has a unique fixed point
on $\Omega$.\end{proof}

We make two more useful observations.

\begin{Lemma}\label{ppower}  Let $G$ be a normal subgroup of $A$ with
$A/G$ generated by $xG$. Assume that $(A,G,A/M)$ is exceptional and
that $x\in M$.  Let $b$ be a power of the prime $r$.
Then $[\Cen_A(y):\Cen_M(y)] \equiv 1 \pmod r$ for
$y= x^b$.
\end{Lemma}

\begin{proof} Let $\omega \in A/M$ be the fixed point of $x$.
Consider the orbit of $\Cen_A(y)$ containing $\omega$.
Of course, $x$ acts on this orbit and $y$ acts trivially.
So $x$ acts on this orbit as an $r$-element.
  The size of
this orbit is $[\Cen_A(y):\Cen_M(y)]$.  Since $x$ has
a unique fixed point $\omega$  on this orbit, the assertion follows.
\end{proof}

\begin{Lemma}\label{cong1} Assume that $(A,G,A/M)$ is exceptional and
$A/G$ is a power of the prime $p$. Then $[A:M]\equiv1\pmod{p}$.
\end{Lemma}

\begin{proof} Let $H=M\cap G$. Then $[M:H]=[A:G]$, a power of
$p$. Also, $M/H$ fixes no $H$--orbits other than $\{M\}$, by
exceptionality. The assertion follows.
\end{proof}

\subsection{Examples of exceptionality} We use Corollary \ref{excor}
to provide some examples of exceptional triples.

\begin{Example}\label{L^t} Let $G=L^t$ with $t>1$ and $L$ a finite
group. Let $x$ be the automorphism of $G$ of order $t$ which
cyclically permutes the coordinates of $G$.
Let $D=\Cen_G(x)$ (the
diagonal of $G$). Let $A$ be the semidirect product of $G$ and
$\gen{x}$. Set $M=\gen{D,x}$.  Then $(A,G,A/M)$ is exceptional if and
only if $(t,\abs{L})=1$.\end{Example}

\begin{proof} If $t$ is relatively prime to $\abs{L}$, this follows by
Corollary \ref{excor}. If $t$ is not relatively prime to $\abs{L}$, we
can choose $1\ne h\in L$ with $h^t=1$. It follows that $xg$ and $x$
are conjugate in $A$ but not in $M$ for $g=(h,\dots,h)$. Lemma
\ref{chig} implies that $x$ has more than $1$ fixed point.\end{proof}

In fact, one can modify the previous example.

\begin{Example}\label{L^t+}  Let $G=L^t$ with $t$
prime  and $L$ a finite group.  Let $u$ be the automorphism of
$G$ of order $t$ which cyclically permutes the coordinates of
$G$. Suppose that $v$ is an automorphism of $L$ of order
prime to $t$ such that $v$
preserves no conjugacy class of elements of order $t$ in $L$.
Let $x=(v,\ldots,v)u \in \Aut(G)$. Let $D=\Cen_G(u)$ (the
diagonal of $G$). Let $A=\gen{G,x}$.
 Set $M=\gen{D,x}$.  Then $(A,G,A/M)$ is exceptional.
\end{Example}

\begin{proof} Note that $x$ preserves no conjugacy class
of $t$-elements in $G$ (if it did, there would be $t$
classes of $L$ permuted by $v$, whence each fixed by $v$,
contradicting the hypotheses).  Now apply Theorem \ref{expi}
with $u=y$ and $\pi = \{t\}$.
\end{proof}

\begin{Example}\label{Chev} Let $x$ be the field automorphism of the
simple Chevalley group $G:=L(q^a)$ determined by the $q$-th power map
on the field of $q^a$ elements. If $a>1$ and $(a,\abs{\Cen_G(x)})=1$,
then $(A,G,A/M)$ is exceptional for $A$ the semidirect product of $G$
and $\gen{x}$ and $M$ any subgroup of $A$ containing $\Cen_A(x)$ which
is normalized by $x$.\end{Example}

The typical case for $M\cap G$ above is the normalizer
of $L(q^b)$ in $G$ for some $b|a$ (and if we consider
primitive actions $a/b$ prime). For example, if $G=L_2(p^a)$
with $a$ odd and $x$ is the Frobenius automorphism, then
we get examples as long as $(a,p(p^2-1))=1$ and $a>1$.

The next example is an elaboration of the previous one
and gives almost all examples of exceptionality for almost simple groups.

\begin{Theorem}\label{Chev2} Let $G=X(q)$ be a finite
simple Chevalley group defined over the field $F_q$
with $q=p^a$.  Let $x$ be a field or graph-field
automorphism of $G$.
     Let $r$ be a prime different from $p$ such that
    $\langle x \rangle$ contains a field automorphism of order $r$ and
    assume that one of the following conditions hold:
     \begin{itemize}
     \item
     there are no $x$-invariant $G$-conjugacy classes of elements
     of order $r$ in $G$; or
     \item $\Cen_X(x)$ contains no elements of order $r$, where $X$
     is the corresponding algebraic group.
     \end{itemize}
     Let $A=\gen{G,x}$ and $M=\Cen_A(y)$, where $y$ is an element of
     order $r$ in $\langle x \rangle$ (so $M$ is a subfield group).
     Then $(A,G,A/M)$ is exceptional.
\end{Theorem}

\begin{proof}  Note that the existence of an $x$-invariant class
  is equivalent to the fact that $r$ divides $|\Cen_G(xg)|$ for some
  $g \in G$.  Since $\Cen_X(xg)$ is conjugate to $\Cen_X(x)$ (by
  Lang's theorem), the second condition above implies the first.
  
  So $r$ does not divide the order of $\Cen_G(xg)$ for any $g \in G$.
  Now apply Theorem \ref{expi}.
\end{proof}

Since one knows the outer automorphism groups of Chevalley groups
quite explicitly (see \cite{GLS3}), it is easy to determine when there
exists an invariant class of $r$-elements.

\begin{Lemma}\label{tab}  Let $G$ be a finite simple Chevalley group of characteristic
$p$.  Let $r$ be a prime.  Assume that there does not exist a conjugacy
 class of $G$ consisting of
elements of order $r$ that is $\Aut(G)$-invariant.
\begin{enumerate}
\item [(i)] There exists an
 automorphism
of $G$ that fixes no conjugacy class of $r$-elements;
\item [(ii)]   $r \ne 2$, and $r=3$  if and only if
$G=L_2(3^a), Sp_4(2^a)$  or $Sz(2^a)$;
\item [(iii)]   $r=p$ if and only
if $G=L_2(q)$ with $q$ odd;
\item[(iv)]   if $G$ is untwisted and has no
graph automorphisms, then such a class exists if and only if
$r$ divides the order of the corresponding group over the prime
field;
\item [(v)] for the other groups, the set of $r$ for which there
  exists an $\Aut(G)$-invariant class of $r$-elements are precisely
  those primes $r$ which divide $i(G)$, with $i(G)$ given in Table
  \ref{Table:1}, page \pageref{Table:1}.
\end{enumerate}
\end{Lemma}

\begin{proof}  The fact that $r \ne 2$ follows by \cite{FGS}.
If $r=p$, this follows easily by induction:
If $L$ is rank $1$ or $2$, the result is clear.
In general, choose a Levi subgroup whose class is $\Aut(G)$-invariant and
apply induction.

So assume that $r \ne p$.  Thus, $r$-elements are semisimple and so
all $r$-classes are invariant
under inner-diagonal automorphisms (because any such element is
contained in a maximal torus
of ${\rm Inndiag}(G)$).

Thus, we need only work modulo inner-diagonal automorphisms.
The complete set
of primes $r$ for which there is an $\Aut(G)$-invariant class of $r$-elements
comes from inspection and Lang's Theorem.  In the table, we give a subgroup of
$G$ that is centralized by a subgroup of $\Aut(G)$ which induce all
automorphisms
and so certainly all such $r$ have the property.
Conversely, by considering field
and graph-field automorphisms (and extending to the algebraic group and
using Lang's theorem),
we see that if there exists an invariant class of $r$-elements,
then $r$ does divide
the order of the group given.

Finally, we prove (i).  Assume that there are no $r$-invariant
classes under $\Aut(G)$. If the automorphism group modulo the
inner-diagonal automorphisms is cyclic, the result is clear.  In
the remaining cases which are the untwisted groups of type $A_n,
n > 1$, $B_2$ or $F_4$ (in characteristic $2$), $D_n$ and $E_6$,
one checks from the table that for any given $r$, one of a field
or graph-field automorphism fixes no $r$-class.
\end{proof}

\begin{table}\label{Table:1}
\caption{$r$-invariant classes for $r \ne p$}
\renewcommand{\arraystretch}{1.4}
\setlength\tabcolsep{15pt}
\begin{tabular}{@{}llp{1.8cm}l}
\hline\noalign{\smallskip}
$G$  & & &  \quad  $i(G)$   \\
\noalign{\smallskip} \hline \noalign{\smallskip}
$L_n(p^a), U_n(p^a), n$ odd     & & &                 $|O_n(p)|$  \\
$L_n(p^a), U_n(p^a), n >2$ even   & & &                  $|\PSp_n(p)|$\\
$Sp_4(2^a), Sz(2^a)$               & & &                $|Sz(2)|$\\
$G_2(3^a), {^2}G_2(3^a)$            & & &             $|{^2}G_2(3)|$\\
$F_4(2^a), {^2}F_4(2^a)$             & & &               $|{^2}F_4(2)|$\\
$O_8^+(p^a)$, ${^3}D_4(p^a)$          & & &             $|G_2(p)|$\\
$O_8^-(p^a)$                           & & &            $|O_7(p)|$\\
$O_{2m}^{\epsilon}(p^a), m >  4$         & & &        $|O_{2m-1}(p)|$\\
$E_6(p^a), {^2}E_6(p^a)$                 & & &           $|F_4(p)|$\\
\noalign{\smallskip} \hline \noalign{\smallskip}
\end{tabular}\\
  \label{table}
\end{table}

There are some other examples for rank $1$ groups.  In
particular, we mention some interesting examples for
$L_2$ (some of these come up in the classification of
exceptional polynomials \cite{FGS}).

The next examples give rise to interesting families of
exceptional polynomials.

\begin{Example}\label{l2example1}
 Let $L=L_2(2^a)$ with $x$ a field
automorphism of odd order $a> 1$.  Set $A=\langle L,x \rangle$ and
$M$ the normalizer of a nonsplit torus of $L$. Then $(A,L,A/M)$
is exceptional.
\end{Example}

\begin{proof}  We can identify $A/M$ with set of subgroups of order
$3$ in $L$.  It suffices to show that every element $y$ in the
coset $xL$ has at most one fixed point.  Since $A/L$ has odd
order, we may take $y$ to be of odd order and so we need show that
$y$ centralizes at most one subgroup of order $3$ in $L$.  By
Lang's Theorem $y$ and $x$ and conjugate over the algebraic group
and have centralizer $\PGL_2(2)=S_3$, whence the result.
\end{proof}

\begin{Example}\label{l2example1a}
  Let $L=L_2(3^a)$ with $a > 1$ odd.  Let $x$ be the field-diagonal
  automorphism of order $2a$.  Set $A=\langle L,x \rangle$ and $M$ the
  normalizer of a nonsplit torus of $L$. Then $(A,L,A/M)$ is
  exceptional.
\end{Example}

\begin{proof} We can identify $A/M$ with the set of involutions of $L$.
  Note that for any $g \in L$, $C:=\Cen_L(xg)$ contains no unipotent
  elements. By Lang's theorem, the centralizer of $xg$ in the
  algebraic group is $\PGL_2(3)$.  If $xg$ centralizes no involution,
  then $C=1$ (since it contains no unipotent elements).  Then $x$ and
  $xg$ are conjugate in $A$, whence $\Cen_L(x)=1$. On the other hand,
  $x$ acts as an involution on $\PSL_2(3)$ and so has a nontrivial
  centralizer.  See also Proposition \ref{fpf}.
\end{proof}

\begin{Example}\label{l2example2}
  Let $L=L_2(p^a)$ with $p > 2$ and $a$ even.  Let $x$ be the product
  of the diagonal automorphism and any field automorphism of even
  order.  Set $A=\langle L,x \rangle$ and $M$ the normalizer of a
  split torus of $L$.  Then $(A,L,A/M)$ is exceptional.
\end{Example}

\begin{proof}  We can identify $A/M$ with the set of
involutions of $L$.  Note that for any $g \in L$,
$C:=\Cen_L(xg)$ contains
no unipotent elements. By Lang's theorem, the centralizer
of $xg$ in the algebraic group is $\PGL_2(p^b)$ where
$a/b$ is even. In particular, $C$ is isomorphic to a subgroup
of $\PGL_2(p^b)$.

Suppose that $C$ has odd order.  Then $C$ is cyclic and contained in a
split torus of $L$ (since its order divides $p^{2b}-1$).  Thus there
exists a unique central involution in $\Cen_L(C)$ (since $a$ is even).
It follows that $xg$ centralizes this involution, whence $C$ cannot
have odd order.  Thus, $xg$ has at least one fixed point, whence
exceptionality holds.
\end{proof}

\begin{Example}\label{VH} Let $G=VH$ where $V$ is a finite
dimensional vector space over the field $\FFF_q$ of $q$ elements and
$H$ is a group of $\FFF_q$-automorphisms of $V$. Let $1\ne r$ be a
divisor of $q-1$ with $H$ containing no element of order
$r$. Choose an element $x\in GL(V)$
of order $r$ acting as a scalar on $V$. Set $M=H\times\gen{x}$ and
$A=VM$. Then $(A,G,A/M)$ is exceptional.\end{Example}

The next example is an exceptional triple where $A/G$ is not cyclic.
There are minor variations on the theme by taking products of
simple groups.

\begin{Example}\label{A5} Let $G$ be a nonabelian simple group. Let
$A=G\times H$ with $H\cong G$. Let $M$ be the diagonal subgroup of
$A$. Then $(A,G\times\bn1,A/M)$ is exceptional. Moreover, if $G=A_5$,
then $(G,A/M)=(G,G)$ may correspond to a covering of genus $0$ (with
inertia groups of order $2$, $3$ and $5$).\end{Example}

\subsection{Nonabelian regular normal subgroups}

The following result is an immediate consequence of
\cite[Lemma~12.1]{FGS}. It was observed with another proof in a
later paper of Rowley \cite{Rowley}.

\begin{Proposition}\label{fpf} Let $G$ be a finite group and $x\in\Aut(G)$
with $\Cen_G(x)=1$. Then $G$ is solvable.\end{Proposition}

\begin{proof} Take a counterexample with $\abs{G}$ minimal. Let $A$ be
the semidirect product $G\rtimes\gen{x}$. The hypothesis is equivalent
to $Gx$ is a single conjugacy class in $A$. This hypothesis persists
on $x$--invariant sections of $G$. So we may assume that $G$ is
characteristically simple. By \cite[12.1]{FGS}, there exists an
involution $t$ in $G$ with $t^G=t^A$. Thus, $A=G\Cen_A(t)$ and so
$gx\in\Cen_A(t)$ for some $g\in G$, so $\Cen_G(gx)$ is
nontrivial. Since $x$ is conjugate to $gx$, also $\Cen_G(x)$ is
nontrivial, a contradiction.\end{proof}

\begin{Lemma}\label{reg} Let $G$ be a nontrivial normal subgroup of
the finite group $A$ acting primitively and faithfully on
$\Omega$. Suppose $(B,G)$ is exceptional, with $G\le B\le A$, and
$B/G$ cyclic. Then $A$ does not contain a regular normal nonabelian
subgroup.\end{Lemma}

\begin{proof} As $A$ is primitive, then every nontrivial normal
subgroup acts transitively. Suppose $N$ is a regular normal subgroup
of $A$. As $G$ is normal in $A$, a minimal normal subgroup of $A$
contained in $G$ is either $N$ or centralizes $N$, so is regular in
either case. Thus we may assume that $G$ contains $N$. Let $M$ be the
stabilizer in $A$ of a point in $\Omega$. Then $N$ is a set of coset
representatives of $M$ in $A$. As $Mnm=Mn^m$ for $n\in N,m\in M$, we
may identify the action of $M$ on $\Omega$ with the conjugation action
of $M$ on $N$. In particular, if $x\in M$, then we may identify the
fixed points of $x$ in $\Omega$ with $\Cen_N(x)$. Since $N$ is a
direct product of nonabelian simple subgroups, it follows by the
previous result that $\Cen_N(x)$ is nontrivial for any automorphism
$x$ of $N$. In particular, no element of $A$ has a unique fixed point
(either it fixes no points or fixes at least 2 points), contrary to
exceptionality.\end{proof}

\subsection{Product action}\label{ProdAction}

\begin{Lemma}\label{L:prod} Suppose $E:=\gFit(A)=L_1\times\dots\times
  L_t$ with $L_i$ a non-abelian group, $t>1$, and with the set of
  $L_i$ being the complete set of $A-$conjugates of $L_1$. Let $M$ be
  a maximal subgroup of $A$ such that $M\cap E = R_1\times\dots\times
  R_t$ where $R_i=M\cap L_i\ne1$.  Let $G$ be a normal subgroup of $A$
  with $A=GM$.

Suppose $(A,G,A/M)$ is (arithmetically) exceptional. Then
$(A_1,L_1,A_1/M_1)$ is (arithmetically) exceptional for some group
$A_1$ with $\gFit(A_1)=L_1$, $M_1$ maximal in $A_1$ and $R_1=M_1\cap
L_1$.
\end{Lemma}

\begin{proof} Set $X=A/M$. Since $E$ is transitive, we may
assume that $E=G$. There is no harm in enlarging $A$ to the maximal
subgroup of $\Aut(E)$ which still acts on $A/M$ (namely, we may assume
that $A=EN$ where $N$ is the normalizer of $M$ in $\Aut(E)$).  Then
$A=A_1\wr S_t$ where $\gFit(A_1)=L_1$ and $M=M_1\wr S_t$ where $M_1$
is a maximal subgroup of $A_1$. Let $A_i$ be a conjugate of $A_1$ via
$g_i \in M$ which maps $L_1$ to $L_i$.  Since $g_i \in M$, $R_i =
R_1^{g_i}$.  Let $M_i = M \cap A_i$ (which is the $g_i$-conjugate of
$M_1)$.  We may identify $A/M$ with $A_1/M_1
\times\dots\times A_t/M_t$. Set $X_i=A_i/M_i$.

First assume that the hypothesis of exceptionality holds and suppose
$(A_1,L_1, A_1/M_1)$ is not exceptional. Then $M_1$ and $R_1$ have a
nontrivial common orbit $Y_1$ on $X_1$. Let $Y_i$ the corresponding
orbit in $X_i$. Note that $S_t$ permutes these orbits and so
$Y=Y_1\times\dots\times Y_t$ is an $M$-orbit.  By construction, it is
also an $M \cap E$-orbit.  Thus, it is a nontrivial common $M, M \cap
G$-orbit, contradicting the exceptionality of $(A,G,A/M)$.  This
proves the claim.

Now assume the hypothesis of arithmetic exceptionality holds and
keep notation as above. Thus there is a group $B$ with $G\le B\le
A$, such that $(B,G,A/M)$ is exceptional and $B/G$ is cyclic. Let
$xE$ be a generator for $B/E$. Suppose that $y=x^d$ is the
smallest power of $x$ which normalizes $L_1$ (note that this
depends only on the coset $xE$ and not on $x$). We can then view
$y$ as acting on $A_1/M_1$. We claim that every element of $yL_1$
has a unique fixed point on $A_1/M_1$.  It suffices to assume
that $x$ is transitive on the set of $L_i$ (if not, consider a
single $x$-orbit and note that the number of fixed points of $x$
on $Z_1\times Z_2$ is just the product of the number of fixed
points -- so for every such decomposition preserved by $x$, it
will have a unique fixed point on both $Z_1$ and $Z_2$).

Then, by conjugating we may assume that $x = s(y, 1, \ldots, 1)$
where $s \in S_t$.  The number of fixed points of $x$ is
precisely the same as the number of fixed points of $y$ on
$X_1$.  Since $xg$ has a unique fixed point for each $g \in E$,
it follows that if $g \in L_1$, then $x=s(yg, 1, \ldots, 1)$ has
a unique fixed point, whence $yg$ does as well.  This proves the
claim.\end{proof}

\subsection{Diagonal action}\label{DiagAction}

As we have seen above, the product action case can be reduced to a
smaller case, so in considering the diagonal case, the essential
situation is where a point stabilizer is a full diagonal subgroup.

\begin{Theorem} Suppose $A$ acts primitively and
faithfully on $\Omega$ and is of diagonal type but does not
preserve a product structure on $\Omega$.  Suppose that $G$ is
normal in $A$ and $A/G$ is cyclic. If $(A,G,\Omega)$ is
exceptional, then the following hold:
\begin{itemize}
\item  $F^*(A)$ is the direct product of $t$ copies of a nonabelian simple
group $L$;
\item  $G$ normalizes each component of $A$;
\item $t$ is a prime and there exists an automorphism
$v$ of $L$  such that $ v$ preserves no
 conjugacy class of elements of order $t$
 and $\Cen_L(vg)=\Cen_L((vg)^t)$ for all $g \in L$ (in particular,
this is the case if $t$ does not divide $|L|$ taking $v=1$ or if
$v$ is a $t'$-element preserving no $t$-class).
\end{itemize}
If all the conditions above hold, then there exists
a corresponding exceptional triple $(A,G,\Omega)$.
\end{Theorem}

\begin{proof} Assume that $(A,G,\Omega)$ is exceptional. Let $M$ be
the normalizer in $A$ of the standard diagonal.  Since $A$ does
not preserve a product structure on $\Omega$, it follows that a
point stabilizer of $E:=F^*(A)$ is a full diagonal subgroup which
we may identify with the standard diagonal subgroup. Moreover,
$A$ acts primitively on the set of components. Suppose that $G$
does not act trivially on the set of components. Then $G$ is
transitive on the set of components and there exists $x \in A$
generating $A/G$ with $x$ normalizing a component $L$. Since $L$
acts semiregularly on $\Omega$ and $\Cen_L(x) \ne 1$, there is
a regular $L$-orbit containing the fixed point of $x$.  The set
of $x$-fixed points on this orbit can be identified with  $\Cen_L(x)$,
a contradiction.

Thus, $A/G$ is a primitive
cyclic group on the set of $t$ components, whence $t$ is prime.
Let $u$ be the automorphism of $E$ cyclically permuting the
components of $L$. So we may assume that $x=u(v,\ldots,v)$ where
$v \in \Aut(L)$ (since $x$ normalizes the diagonal subgroup of
$E$).

Since $x$ normalizes no other diagonal subgroup, it follows that
$\Cen_A(x) \le M$.  Let $g_1 \in \Cen_L(v^t)$.  Then define $g_i,
i=2,\ldots,t$ by
$$
g_{i+1}vg_i^{-1} = v,
$$
for $i=1,\ldots,t_1$.  It follows that $(g_1,\ldots,g_t)$
commutes with $x$.  Thus, $g_i=g_1$ for all $i$, whence $g_1$
commutes with $v$.  Thus, $\Cen_L(v)=\Cen_L(v^t)$ as claimed.

Applying the same argument to $vg$ yields the fact that
$\Cen_L(vg)=\Cen_L((vg)^t)$ for all $g \in G$.

Conversely, we now show that the triple $(A,G,A/M)$ is exceptional.
We will show that every element $xg$ for $g \in G$
has at least one fixed point, whence the result holds.

So consider $y:=x(g_1,\ldots,g_t)=u(vg_1,\ldots,vg_t)$. We need
to show that this is conjugate to $z:=u(vh,\ldots,vh)$ for some
$h \in L$.  Note if we take $t$-th powers we see that $vg_t\cdots
vg_1$ must be conjugate to  $(vh)^t$.

Let $J=\langle L,v \rangle$. We first show that there exists $h
\in H$ such that $vg_t\cdots vg_1=(vh)^t$. Consider the mapping
$va \mapsto (va)^t$ for $a \in L$. If $(va)^t=v^t$, then $va \in
\Cen_J(v^t)=\Cen_J(v)$. Thus, $\Cen_J(v) = V \times V'$ where $V$
is the Sylow $t$-subgroup of $\langle v \rangle$ and $V'$ is a
$t'$-group (this follows from the fact that $v$ fixes no
$t$-class of $L$ and that $\Cen_J(v)=\Cen_J(v^t)$). Clearly, the
$t$-th power map is injective on $vV'$.  Thus, $(va)^t=v^t$
implies that $a=1$. Since $vg$ has the same property as $v$ for
any $g \in L$, the same argument shows that the $t$-th power map
is injective on $vL$. Thus, the image of the $t$ power on $vL$ is
precisely $v^tL$. So choose $h$ satisfying the equation above.

We now show that $y$ and $z$ are conjugate
(for this choice of $h$). Choose $s_1,\ldots, s_t \in L$ satisfying
$$
s_2vg_1s_1^{-1}=s_3vg_2s_2^{-1} = \ldots s_tvg_{t-1}s_{t-1}^{-1}=vh,
$$
with $s_1 = 1$.

Note that each $s_i$ is uniquely determined.

Multiplying and telescoping we see that
$$
(vh)^{t-1} = s_tvg_{t-1}\cdots vg_1s_1^{-1}.
$$

We now compute that
$$
s_1vg_ts_t^{-1} = s_1(vh)^t(vg_{t-1}\cdots vg_1)^{-1}s_t^{-1}
  = (vh)^t (vh)^{1-t}=vh.
$$

A straightforward computation shows that
conjugating $u(vg_1,\ldots,vg_t)$ by
$(s_1,\ldots, s_t)$ gives
$u(vh,\ldots,vh)$.

Thus every element in the coset $xG$ has at least one
fixed point as claimed.
\end{proof}

One can show that any automorphism $v$ of a simple group satisfying
the conditions above must have order prime to $t$ (and so the example
is already given in Example \ref{L^t+}). The only proof that we can
see is by inspection using the classification.

By \cite{FGS}, it follows that $t > 2$ and moreover $t$ may be taken
to be $3$ if and only if $L=Sz(q)$, $L_2(3^a)$ or $Sp_4(q)$ with $q$
even.  One can use this result to give an alternative proof of the
nonexistence of exceptional triples of genus zero of diagonal type in
the next section.

\subsection{Almost simple groups}\label{ASAction}

In this section, we consider the exceptionality of triples $(A,G,A/M)$
in the case where $\gFit(A)$ is a non-abelian simple group $L$. We
essentially classify all such arithmetically exceptional examples.

We first assume that $L=L_2(p^a)$ with $p$ prime and $p^a > 3$.  The
examples include the important examples discovered in \cite{FGS} with
$p=2$ or $3$, $a > 1$ odd, $A=\Aut(L)$ and $M$ the normalizer of a
nonsplit torus.  These lead to interesting families of exceptional
polynomials (see \cite{PM:Ex28}, \cite{CohenMatthews}, \cite{LZ},
\cite{GZ}, \cite{GRZ}).

\begin{Theorem}\label{l2aodd} Let $G=L_2(2^a)$ with $a > 1$ odd. Let
  $H$ be a proper subgroup of $G$. The following are
  equivalent:\begin{itemize}
  \item[(i)] There exist subgroups $M$ and $A$ of $\Aut(G)$ with $G<A$
    such that $A=GM$ and $G\cap M=H$ with $(A,G,A/M)$ exceptional;
  \item[(ii)] $L_2(2)\subseteq H$ and $(3,[G:H])=1$;
  \item[(iii)] $H$ contains a dihedral subgroup of $G$ of order
    $2\cdot 3^b$ where $3^b$ is the order of a Sylow $3$-subgroup of
    $G$.\end{itemize}\end{Theorem}

\begin{proof} Assume (i). Let $x\in\Aut(G)=\PgL(2,2^a)$ denote the
  Frobenius automorphism generating the full group of field
  automorphisms. Let $y$ be a power of $x$ with $yG$ a generator for
  $A/G$. It follows by Lemma \ref{Ki} that $\Cen_A(x)\subseteq M$.
  Suppose that $3$ divides the index $[G:H]$. Let $T$ be a Sylow
  3-subgroup of $H$. Let $S$ be the (unique) Sylow 3-subgroup of $G$
  properly containing $T$. Then $S$ is cyclic. By the Frattini
  argument, we may choose $yg$ normalizing $T$, with $g\in H$. We may
  replace $yg$ by $(yg)^{2}$ and assume that $yg$ has odd order (since
  $a$ is odd, $(yg)^2G$ still generates $A/G$). Since the normalizers
  of $S$ and $T$ coincide, it follows that $yg$ normalizes $S$. Since
  $yg$ has odd order, this implies that $yg$ centralizes $S_1/T$,
  where $[S_1:T]=3$. By Lemma \ref{Ki}, this implies that
  $S_1\subseteq H$, a contradiction. Thus (i) implies (ii).

  Assume (ii). Then $H$ contains a Sylow 3-subgroup of $G$ and an
  involution inverting an element of order $3$. This implies that $H$
  contains the dihedral group claimed in (iii).

  Assume (iii). Let $D$ denote this dihedral group. Note that any
  subgroup $H$ of $G$ containing $D$ is normalized by a conjugate of
  the field automorphism generating the full group of field
  automorphisms.  Thus, we may assume that $A=\Aut(G)$.

  Replacing $D$ by a conjugate, we may assume that $x$ normalizes $D$.
  Let $J=\gen{D,x}$. It suffices to show that $A/J$ is exceptional.
  Let $H$ be the unique normalizer of a nonsplit torus in $G$,
  containing $D$. Let $M=\gen{H,x}$. It is straightforward to check
  that $(M,H,M/J)$ is exceptional (e.g., pass to ${\bar
    M}:=M/{\core_M(J)}$ and observe that this group is a Frobenius
  group), whence by Lemma \ref{excomp} it suffices to show that
  $(A,G,A/M)$ is exceptional.

  We can identify $A/M$ with the set of subgroups of order $3$ in $G$
  (since $M$ is the normalizer of a subgroup of order $3$ and any two
  such are conjugate in $G$). It suffices to show that $xg$ for $g\in
  M$ has at most one fixed point on $A/J$. We may square this element
  to assume that $xg$ has odd order. The fixed points of $xg$ are the
  subgroups of order $3$ centralized by $xg$. By Lang's Theorem, $xg$
  is conjugate (in the algebraic group) to $x$, whence its centralizer
  (in the algebraic group) is $L_2(2)$. So it centralizes at most one
  subgroup of order 3 as desired.\end{proof}

\begin{Theorem} Let $G=L_2(2^a)$ with $a$ even. Let $H$ be a proper
  subgroup of $G$. There exist subgroups $M<A\le\Aut(G)=\PgL_2(2^a)$
  with $A=GM$ and $M\cap G=H$ with $(A,G,A/M)$ exceptional if and only
  if $H=L_2(2^{b})$ where $(6,a/b)=1$.\end{Theorem}

\begin{proof} Let $x$ denote a generator for the full group of field
  automorphisms. Let $y$ be the smallest power of $x$ in $A$. Assume
  $(A,G,A/M)$ is exceptional. We may assume that some power of $y\in
  M$. Then $\Cen_G(y)\subseteq H$. In particular, $H$ is an
  irreducible subgroup of $G$. Since $x$ has even order, it follows
  that $xg$ and $x$ are conjugate in $A$ for $g$ an involution in
  $\Cen_G(x)$. Indeed, we can find an involution $z \in \Cen_G(g)$
  with $zxz=xg$. It follows by Lemma \ref{Ki} that $z \in M$. In
  particular, the order of $H$ is divisible by $4$. The only such
  subgroups of $G$ are $L_2(2^b)$ with $b$ dividing $a$. In
  particular, the $G$-conjugacy class of such a subgroup is
  $A$-invariant, whence we may assume that $A=\Aut(G)$ and
  $M=\Nor_A(H)$.

  We claim that $(6,a/b)=1$. Suppose not, and let $p \in \{2,3\}$ with
  $p$ dividing $a/b$. Put $U=\Nor_A(L_2(2^{a/p})$. Then $(A,G,A/U)$ is
  exceptional by Lemma \ref{excomp}. Write $a=p^f d$, where $d$ is
  prime to $p$. Set $y=x^{p^f}$. Then $\Cen_G(y) \cong
  L_2(2^{p^{f}})$. Moreover, if $y'\in M$ is conjugate to $y$ in $A$,
  it follows (by Lang's theorem) that $\Cen_{M}(y') \cap G$ is
  isomorphic to a subgroup of $\Cen_G(y)$. In particular,
  $[\Cen_A(y'):\Cen_M(y')]$ is divisible by $p$ for any such $y'$. It
  follows that the number of fixed points of $y$ is divisible by $p$
  and also for $x$, a contradiction.

  All that remains to prove is that $(A,G,A/M)$ is exceptional for $H
  \cong L_2(2^b)$ whenever $(6,a/b)=1$, where
  we take $A=\Aut(G)$ and
  $M=\Nor_A(H)$. This follows by Theorem \ref{Chev2}.
\end{proof}

\begin{Theorem}\label{L2odd} Let $G=L_2(p^a)$ with $p$ odd
  and $p^a > 3$. Let $X = \Aut(G)$. Let $H$ be a proper subgroup of
  $G$. There exists a subgroup $A$ of $\Aut(G)$ and a subgroup $M$ of
  $A$ such that $A=MG$, $M \cap G=H$, $A/G$ is cyclic and the triple
  $(A,G,A/M)$ is exceptional if and only if $A \le \Nor_X(H)G$ and one
  of the following holds:
\begin{itemize}
  \item[(a)] $H=L_2(p^b)$ with
    $(a/b,p^2-1)=1$;
  \item[(b)] $a$ is even and $H$ is contained
in the normalizer $N$ of a split torus of
    $G$ and $[N:H]$ is not divisible by any prime
    dividing $|L_2(p)|$;
  \item[(c)] $p=3$, $a$ is odd or $a/2$ is odd,
    and $[G:H]$ is odd.
\end{itemize}\end{Theorem}

\begin{proof} Let $H$, $A$ and $M$ be given satisfying the
  hypotheses.

We will handle the case $p=3$ and $a$ is odd at the end of
the proof.  So we exclude this case for the moment.

  Let $f$ be the largest power of $2$ dividing $a$. Let
  $U = L_2(p^f) \le G$. Note that all subgroups of $G$ isomorphic to
  $U$ are conjugate in $G$. Thus, we may choose $x \in A$ which
  normalizes $U$ and generates $A/G$. Since $x$ has a fixed point on
  $A/M$, we may replace $M$ by a conjugate and assume also that $x \in
  M$. It follows by Lemma \ref{desc1} that $(A_1,U,A_1/M_1)$ is
  exceptional where $A_1=\gen{U,x}$ and $M_1=M \cap A_1$.

  Note that $M_1$ may contain a central subgroup of $A_1$ (i.e.\ the
  $\Cen_A(U)$) and so we may pass to a homomorphic image which we
  still call $A_1$. Now $A_1/U $ is a 2-group, whence it follows that
  $[A_1:M_1]$ is odd.

  We claim that one of the following holds:
\begin{itemize}
\item[(i)] $A_1=M_1$; or
\item[(ii)] $f > 1$ and $M_1 \cap U$ contains both a Sylow
  $2$-subgroup of $U$ and a cyclic subgroup of order $(p^2-1)/2$; or
\item[(iii)] $p^f=9$ and $M_1$ contains a Sylow $2$-subgroup of $A_1$.
\end{itemize}

Suppose $A_1 \ne M_1$. We may replace $x$ by a power and assume that
$x$ is a 2-element (since $A_1/G_1$ is a $2$-group). Let $J_1$ be a
maximal subgroup of $A_1$ containing $M_1$. This implies (using the
fact that $M_1$ has odd index in $A_1$) that $H_1:=J_1 \cap U$ is the
normalizer of a torus, or $H_1$ is $A_4$ or $S_4$. If $f > 1$, then a
Sylow $2$-subgroup of $U$ has order at least $8$ and if $f > 2$, then
a Sylow $2$-subgroup has order at least $16$. Since $H_1$ has odd
index in $U$, this implies that if $H_1=A_4$ or $S_4$, then either
$f=1,$ or $f=2$ and $H_1=S_4$.

On the other hand, $M_1 \cap U$ contains elements of order $(p \pm
1)/2$. Moreover, if $a$ is even, then $x$ preserves the conjugacy
class of elements of order $p \pm 1$, and so $M_1 \cap U$ contains
such elements. If $p > 5$, then $A_4$ and $S_4$ cannot contain such
elements. The same holds if $p=5$ and $f > 1$. Thus, the only
possibility for $H_1=A_4$ is with $f=1$ and $p=5$ (note that the case
$p=3$ with $f=1$ is excluded). It is straightforward to see that the
triple is not exceptional.

Similarly, if $H_1=S_4$, the only possibility  is $p^f=9$.
This case can occur.  We deal with this case below.

So for now we assume that $p^f \ne 3$ and $p^f \ne 9$.

Thus (i) or (ii) holds. Note that it is straightforward to check that
if (ii) holds and not (i), then $x$ must induce a field-diagonal
element on $U$ (otherwise, $x$ induces a diagonal automorphism and all
semisimple classes are fixed by $x$ or $x$ induces a field
automorphism and all unipotent classes are fixed -- in either case, we
see that $A_1=M_1$).

First assume that (i) holds. We continue to exclude the case that
$p=3$ and $a$ is odd.

Thus we may assume that $L_2(p^f) \subseteq H$. It follows that $H
\cong L_2(p^b)$ with $a/b$ odd. Let $r$ be an odd prime dividing
$a/b$. We claim also that $r$ does not divide $p^2-1$. Suppose it
does. It suffices to assume that $H=L_2(p^{a/r})$. Let $y$ denote the
Frobenius automorphism on $G$ which we assume normalizes $H$.  Then we
may take $x=y^i$ or $x=y^iz$ where $z$ induces the diagonal
automorphism on $G$ (and normalizes $H$). Since $r|p^2-1$ then $r$
divides the order of a torus, split or nonsplit, and we may take $z$
to centralize this torus. Then we may choose $g$ in this torus such
that $g\notin H$ but $g^r \in H$. Then $gxg^{-1} \in M$, whence $g \in
M$ by Lemma \ref{Ki}, a contradiction.

We prove next that if $H=L_2(p^b)$ with $(a/b,p^2-1)=1$, there is an
exceptional triple. By induction (and the transitivity property for
exceptionality), we may assume that $a/b=r^c$ for $r$ a prime and
$r^{c+1}$ not dividing $a$. We distinguish two cases, $r=p$ and $r\ne
p$.

First consider the case that $r \ne p$. Let $x$ be the element in
$\Aut(G)$ generating the full group of field automorphisms (and
normalizing $H$). Write $x=yz=zy$ where $y$ is an $r$-element and $z$
is an $r'$-element. We may take $H=\Cen_G(y)$. We will show that $xg$
has at most $1$ fixed point for every $g \in G$. Assume that $xg$ has
a fixed point. So we may assume that $g \in H$. Write $xg=y'z'=z'y'$
where $y'$ is an $r$-element and $z'$ is a $r'$-element. By Lang's
Theorem, $\Cen_G(xg)$ is an $r'$-group.  Thus, $y'$ is the unique
Sylow $r$-subgroup in $\Cen_A(xg)$. Also, $y \in \Cen_A(xg)$ (since $g
\in H$) and $y$ and $y'$ have the same order $r^c$. Since $yG=y'G$,
this implies $y=y'$.

Now suppose that $xg' \in M$ is conjugate to $xg$ in $A$. Arguing as
above, we see that the $r$-part of $xg'$ is $y$. Then the conjugacy in
$A$ implies conjugacy in $\Cen_A(y)=M$. Moreover, $\Cen_A(xg)
\subseteq \Cen_A(y)=M$, whence Lemma \ref{chig} implies that $xg$ has
one fixed point.

Now consider the case $r=p$. Let $\sigma$ generate the full group of
field automorphisms. Let $\tau$ denote a diagonal automorphism which
centralizes a split torus normalized by $\sigma$. Let $x=\sigma\tau$.
We argue as above. So let $x=yz=zy$ where $y$ is a $p$-element and $z$
is a $p'$-element. We may take $H=\Cen_G(y)$.  Since $\tau$ switches
the two conjugacy classes of elements of order $p$ in $G$ and $\sigma$
fixes the classes, $xg$ switches the two classes for any $g \in G$. It
follows that $\Cen_G(xg)$ is a $p'$-group. Now we argue as above to
conclude that $xg$ has at most $1$ fixed point for all $g \in G$.

Now consider (ii). The only possible subgroups $H$
either are contained in the normalizer $N$ of a  split torus
or are a subfield
group.  The latter we handled above.  It is clear from (ii)
that $H$ must contain (up to conjugacy) the normalizer of
the  split torus of $L_2(p^2)$.
Let $x$ be the
product of a diagonal automorphism and a field automorphism of even
order (recall $a$ is even in this case). In this case the normalizer
of the split torus  is precisely
the centralizer of an involution in $G$. Thus, we can identify the set
we are acting on with the set of involutions in $G$. Now one computes
 that for any such $x$ normalizing $N$, $\Cen_G(x)$ is cyclic
of order $q_0 \pm 1$ for a suitable power $q_0$ of $p$
(see the proof of Example \ref{l2example2}). In particular,
$x$ centralizes a unique involution. Thus, each element in the coset
$xG$ fixes at most $1$ point and the cover is exceptional.
It is straightforward to see that $(N,H,N/H)$ is exceptional
precisely given the conditions in (b).

Next  consider the case $p=3$ and $a > 1$ is odd. Let $x \in
\Aut(G)$ be the commuting product of a diagonal automorphism of order
$2$ and the field automorphism of order $a$. So let $y=x^i \in A$ with
$yG$ generating $A/G$. Now $\Cen_G(x)$ is a generated by an involution
$t$, and $F:=\O2(\Cen_G(x^2))$ is an elementary abelian group of order
$4$ which is $y$-invariant. Since $y$ centralizes $F/\langle t
\rangle$, it follows by Lemma \ref{Ki} that $F \subseteq H$. Denote
the elements of $F$ by $\{1,t,u,tu\}$.

Note that the overgroups of $F$ in $G$ are the subgroups of
$\Cen_G(t)$ (a dihedral group of order $3^a+1$) and the subgroups
$L_2(3^b)$ with $b|a$. These latter subgroups are just the
centralizers of $x^{2a/b}$. Thus, it suffices to show that $(A,G,A/F)$
is exceptional for $A=\Aut(G)$. Suppose $xg$ has at least one fixed
point. Thus, we may assume that $g \in F$. Then $x$ and $xt$ are
conjugate in $M=\gen{F,x}$ and $xu$ and $xtu$ are conjugate in $M$.
Since $xu$ has order $4a$ and $x$ has order $2a$, $x$ and $xu$ are not
conjugate in $A$. Note that $\Cen_G(x)=\gen{t}$. Also $\Cen_G(xu)
\subset \Cen_G(x^2) =L_2(3)$. On the other hand, $xu$ cannot
centralize an element of order $3$ (since it switches the two
conjugacy classes of elements of order $3$ in $G$). Thus
$\Cen_G(xu)=\Cen_G(x) \subset F$. It follows by Lemma \ref{chiex} that
$xg$ has a unique fixed point for every $g \in F$, whence the same is
true for every $g\in G$, as desired.

Finally, consider the case that $p^f=9$.  As we noted above,
it follows that $H \cap G$ contains a Sylow $2$-subgroup of
$L_2(9)$ (which is a  Sylow $2$-subgroup of $G$).
  Let $x$ be the product of a field automorphism
of order $a/2$ and a commuting diagonal automorphism of order $2$.
Let $A=G\langle x \rangle$. As above, it follows that $(A,G,A/M)$
is exceptional with $M=L_2(9)\langle x \rangle$ (because
$(a/2,p^2-1)=(a/2,2)=1$). Also, one checks easily that
$(M,L_2(9),M/H)$ with $H$ a Sylow $2$-subgroup of $M$ is
exceptional (this example is given in terms of  $A_6$ for the
main theorem). Thus, $(A,G,A/H)$ is exceptional and the same is
true for any overgroup of $H$.
\end{proof}

The previous results (and their proofs -- which indicate the possibilities
for $A$) give the classification of all quadruples
$(A,B,G,A/M)$ where $G=L_2(q)$ or $PGL_2(q)$ with $B/G$ cyclic,
$A$ primitive on $A/M$ and $(B,G,A/M)$ exceptional.

\begin{Theorem}\label{L2ex}
  Let $L=L_2(q)$ with $q > 3$. There exists a quadruple $(A,B,G,A/M)$
  with $\gFit(A)=L$, $A$ primitive on $A/M$, $G=L_2(q)$ or
  $\PGL_2(q)$, $G < B \le A$ with $B/G$ cyclic and $(B,G,A/M)$
  exceptional, if and only if one of the following holds:
\begin{itemize}
\item[(a)] $q=2^a$ and\begin{itemize}
  \item[(i)] $M\cap L=L_2(2^{b})$ with $3 < a/b=r$ a prime and $B/G$
    generated by a field automorphism $x$ such that $\Cen_L(x)$ is not
    divisible by $r$ (in particular, $r$ is odd and does not divide
    $4^{a/e}-1$ where $e$ is the order of $x$);
  \item[(ii)] $a>1$ is odd and $M\cap G$ is dihedral of order
    $2(q+1)$;\end{itemize}
\item[(b)] $q=p^a$ with $p$ odd and\begin{itemize}
  \item[(i)] $M \cap L= L_2(p^b)$ with $r:=a/b$ prime and $r$ not a
    divisor of $p(p^2-1)$, and $B$ any subgroup containing an
    automorphism $x$ which is either a field automorphism of order $e$
    or the product of a field automorphism of order $e$ and a diagonal
    automorphism, such that $r$ does not divide $p^{2a/e}-1$;
  \item[(ii)] $M \cap L=L_2(p^{a/p})$ and $B$ a subgroup of $\Aut(G)$
    such that $p|[B:G]$ and $B$ contains a diagonal automorphism (so
    in fact $2p|[A:G]$);
  \item[(iii)] $a$ is even, $M \cap L$ is the dihedral group of order
    $p^a-1$ and $B/G$ is generated by a field-diagonal automorphism;
  \item[(iv)] $p=3$, $a$ odd, $M \cap L$ is dihedral of order $p^a+1$
    and $B=\Aut(G)$.
\end{itemize}\end{itemize}\end{Theorem}

We consider the remaining rank one groups in the next result.

\begin{Theorem}\label{rank 1}  Let $L$ be a finite simple group of
Lie type and rank $1$ defined over $\FFF_q$ and not of type $L_2$.
Assume that
  $(A,B,G,M)$ is such that $L \n A \le \Aut(L)$, $L \le G < B \le A$
  with $B/G$ cyclic, $M$ maximal in $A$, $M$ does not contain $L$
  and $(B,G,B/B\cap M)$ exceptional.  Then one of the following holds
  (and all such give rise to examples):
  \begin{itemize}
  \item[(i)]  $L=Sz(2^a)$ with $a > 1$ odd and $M \cap L = Sz(2^{a/b})$
  with $b$ a prime not $5$ and $a \ne b$  or  $M \cap L$ the normalizer of a Sylow
  $5$-subgroup (which is the normalizer of a torus);
  \item[(ii)] $L=Re(3^a)$ with $a > 1$ odd and $M \cap L=Re(3^{a/b})$
  with $b$ a divisor of $a$ other than $3$ or $7$;
  \item[(iii)]   $L=U_3(p^a)$ and $L \cap M=U_3(p^{a/b})$ with $b$ a prime
     not dividing $p(p^2-1)$.
  \item[(iv)] $L=U_3(2^a)$ with $a > 1$ odd and $M \cap L$ is the subgroup preserving
  a subspace decomposition into the direct sum of $3$ orthogonal nonsingular $1$-spaces.
  \end{itemize}
\end{Theorem}

\begin{proof}  The rank one groups other than $L_2$ are $Sz$, $Re$ and $U_3$.
There is no harm in assuming that $G=L$.

If $L=Sz(2^a)$ with $a > 1$, then the outer automorphism group
consists of field automorphisms and so $Sz(2) \le M$.  The
maximal overgroups of $Sz(2)$ are subfield groups and the
normalizer of a  maximal torus containing a Sylow $5$-subgroup.
We may assume that $A/L$ is generated by a field automorphism $x$.

If $M \cap L$ is a maximal subfield group, it is of the form
$Sz(2^{a/b})$ with $b$ an odd prime. If $b \ne 5$, $(A,L,A/M)$ is
exceptional by Theorem \ref{Chev2}.  If $b=5$ and $z$ is the $5'$-part
of $x$, then $\Cen_L(z)$ is divisible by $b$ and so by Lemma
\ref{ppower}, $M$ contains a Sylow $5$-subgroup of $\Cen_L(z)$, which
is not the case.

Note that $M \cap L$ is not maximal if $a=b$ is a prime other than
$3$.  If $a=b=3$, then $Sz(2)$ is the normalizer of the Sylow
$5$-subgroup and we are in the remaining case.

Suppose that $M \cap L$ is the normalizer of the maximal torus containing
a Sylow $5$-subgroup. Thus, $M \cap L = T.4$ where $T$
is cyclic. We may assume that $x$ is a generator
for the full group of field automorphisms (for this will leave
invariant the normalizer of a maximal torus).
We need to show that all elements in the coset
$x(M \cap L)$ have a unique fixed point.  By raising $x$ to a power,
it suffices to assume that $x$ has odd order and so we need only
consider elements in the coset $xT$. Since $T=T_1 \times T_2$,
where $T_1$ is a Sylow $5$-subgroup of $T$ and $[x,T_2]=T_2$,
it suffices to consider elements $xt$ with $t \in T_1$.

Any such element centralizes the subgroup of order $5$ in $T$.  If
some power of $xt$ is a nontrivial element of $T$, then its
centralizer is contained in $T$.  Any two such $L$ conjugate elements
therefore would be conjugate via the normalizer of $T$ and so have a
unique fixed point.  If $\langle xt \rangle$ does not intersect $T$,
then by \cite[7.2]{GorLy}, $xt$ is conjugate to $x$.  Since $xt$ and
$x$ have centralizers isomorphic to $Sz(2)$ which have a common
(normal) subgroup of order $5$, it follows that they are conjugate in
$M$.  Thus, $x$ also has a unique fixed point.

Next consider $L=Re(3^a)$ with $a$ odd.
We may take $a > 1$ (since $Re(3)' \cong L_2(8)$).
Then $A/L$ is generated by a field automorphism.  Thus, $M \cap L$ contains $Re(3)$.
The only overgroups of $Re(3)$ are subfield groups.
So $M$ contains $Re(3^{a/b})$ with $b$ an odd prime.  Let $x$ be the generator
of the group of field automorphisms.  Write $x=yz=zy$ where $y$ is a $b$-element and
$z$ is a $b'$-element.  First we show that $b \ne 3$ or $7$.  If so, then
$\Cen_L(z)$ is divisible by $b$ and so by Lemma \ref{ppower},
$M$ contains a Sylow $b$-subgroup
of $\Cen_L(z)$, which is not the case. If $b \ne 3$ or $7$, then
we do have exceptionality as in Theorem  \ref{Chev2}.

Next consider $L=U_3(p^a)$ with $p $ odd and $p^a \ne 3$.
Let $x$ be a generator for $A/G$.
If $x$ does not involve an
outer diagonal automorphism, then $x$ normalizes $O_3(p^a)$ (and acts as a field
automorphism) and so by the proof
for $L_2$,  $M$ contains
$O_3(p^{a/b})$ with $b$ prime to $p(p^2-1)$.  In particular,
$M \cap L$ acts irreducibly and primitively.  Also, $M$ contains transvections
as the class of transvections is stable under the full automorphism group.
It follows that $M$ contains $U_3(p^{a/b})$ with $b$ prime and $(b,p^2-1)$.
 Arguing as in the previous paragraph shows that $b \ne p$.
The subfield examples with $(b,p(p^2-1)=1$ do give examples by
Theorem \ref{Chev2}.

Suppose that $x$ does involve an outer diagonal automorphism.  Note that
$p \ne 3$ in this case.  Then $A/L$ is
generated by an element of order dividing
$3a$ and acts as a field automorphism on $SU_2(p^a)$ (the stabilizer of
a nonsingular vector).  It follows that $M$ contains $SU_2(p^{a/b})$
with $(b,p(p^2-1))=1$ (by the proof of the $L_2$ case) and is
a subfield group allowed in the conclusion  or that
$M \cap L$ acts reducibly and so $M$ is the stabilizer of a nonsingular
$1$-space.  A straightforward argument shows that this latter case is impossible
(see the argument in the next result for the general unitary case).

Next consider $U_3(2^a)$ with $a > 1$.  Let $x$ be a generator for
$B/L$.  If $a$ is even, we use the result for $U_2=L_2$ and see
that $L \cap M$ contains $L_2(2^{a/b})$ with $(b,6)=1$.  It follows
that either $M \cap L$ is the stabilizer of a nonsingular $1$-space,
a contradiction as above, or $M$ is a subfield group satisfying the conclusion
(and those subfield groups do give rise to examples by Theorem \ref{Chev2}.

So consider $a > 1$ odd.

By Theorem \ref{l2aodd}, it follows that $L \cap M$ contains
the dihedral group of degree $3^c$ where $3^c$ is the order of the
Sylow $3$-subgroup of $L_2(2^a)=U_2(2^a)$.  Since $M$ is maximal, this
implies that one of the following holds:
\begin{itemize}
\item[(i)]   $M$ is the stabilizer of a nonsingular $1$-space;
\item[(ii)]  $M$ contains $U_3(2^{a/b})$, with $b$ an odd prime; or
\item[(iii)] $M \cap L$ acts imprimitively on the natural module.
\end{itemize}

We have seen before that (i) is impossible.  If (ii) holds, we only
need to show that $b \ne 3$.  If so, write $x=yz=zy$ where $y$ is a
$3$-element and $z$ is a $3'$-element.  Then $[\Cen_B(z):\Cen_M(z)]$ is
divisible by $3$, contradicting Lemma \ref{ppower}.

So suppose that (iii) holds.  If $a$ is prime to $3$, then
let $x$ generate $A/L$ with $x$ of order $3a$ (the product
of a commuting diagonal automorphism of order $3$ and a field
automorphism of order $a$).  We first see that
$(A,L,A/H)$ with $H$ the normalizer of $U_3(2)$ is exceptional
by Theorem \ref{expi}.
On the other hand, $(H, U_3(2), H/D)$ with $D \cap L$
a maximal imprimitive subgroup of $L$ (of order $18$) is also
exceptional.  Thus, $(A,L,A/E)$ is exceptional for any overgroup
$E$ of $D$ including the stabilizer of a system of three mutually
orthogonal nonsingular $1$-spaces, leading to the example as given
in the conclusion.

Now suppose that $3|a$.  Let $f$ be the maximal $3$ power dividing $a$.
Let $x$ be the product of the field automorphism of order $a$ and
the diagonal automorphism of order $3$ and let $x$ generate
$A/L$.  So $A/L$ has order $a$.

Note that $(A,L,A/\Nor_A(U_3(2^f))$ is exceptional.
By the transitivity of exceptionality, it therefore suffices to assume that
$a=f$.  Note that any element in $xL$ has the centralizer of order
$2\cdot 3^b$ for some $b$.  Any element of order a power of $3$ has a fixed
point as $M$ contains a Sylow $3$-subgroup of $A$.  So we need consider
elements $xg$ which have even order and in particular are contained
in the centralizer of a transvection of $L$.  The $3$-part of such an element
is contained in the normalizer of a the torus of a Borel subgroup.
One sees easily that such an element permutes an orthonormal base (indeed,
fixing one $1$-space and permuting the other $2$).  Thus, every element
in $xL$ has at least $1$ fixed point, whence the result.

Finally we show that if $b$ is prime and $(b,p(p^2-1)=1$, we do get an example.
Let $x$ be of order $2a$ (essentially generating the field automorphism
on the field on size $p^{2a}$).  Observe that there are no $x$-stable conjugacy
classes of $L$ of elements order $b$ in $L$. Thus $\Cen_L(xg)$ has order prime
to $b$ for every $g \in L$, whence the result follows by Theorem \ref{expi}.
\end{proof}

Before doing the general case, we will need the following result.

\begin{Lemma}\label{Sylowstructure}  Let $G$ be a simple Chevalley
  group over the field of $q=p^a$ elements.  Let $r$ be a prime
  dividing the order of $G$ such that $r$ does not divide the order of
  the Weyl group of $G$.  Let $R$ be a Sylow $r$-subgroup of $G$.
  Then $R$ is a homocyclic abelian group and if $q=q_0^r$ and
  $M=\Cen_G(\sigma)$ with $\sigma$ a field automorphism of order $r$,
  then $M \cap R$ is contained in the Frattini subgroup of $R$.
\end{Lemma}

\begin{proof}  There is no harm in replacing $G$ by the corresponding
  universal group.  It is well known (see \cite[10.1]{GorLy}) that $R$
  is a homocyclic abelian group (i.e.\ a direct product of cyclic
  groups of the same order) with the exponent being equal to the
  $r$-part of $\Phi_m(q)$, for a suitable $m$, where $\Phi_m$ is the
  $m$th cyclotomic polynomial.  See the references for the value of
  $m$ (the important point for us is that while $m$ depends on $q$, it
  is the same for $q_0$ and $q_0^r$).  If $R=1$, there is nothing to
  prove.  Otherwise, the $r$-part of $\Phi_m(q)$ is strictly greater
  than the $r$-part of $\Phi_m(q_0)$ which is the exponent of $R \cap
  M$ (note that $M$ may pick up some extra diagonal automorphisms, but
  our assumption is that there are no outer diagonal automorphisms of
  order $r$).  Thus, we see that the exponent of $R \cap M$ is
  strictly less than that of $R$ and since $R$ is homocyclic, this
  implies that $R \cap M$ is contained in the Frattini subgroup of
  $R$.
\end{proof}

The previous result will be used in conjunction with the following
fact -- if an $r'$-group $S$ acts on an abelian $r$-group $R$, then $R
= [S,R] \times \Cen_R(S)$ and so if $\Cen_R(S)$ is contained in the
Frattini subgroup of $R$, then $\Cen_R(S)=1$.

\begin{Corollary}\label{subfield1}  Let $L=L(q_0^r)$ be a finite simple
Chevalley group in characteristic $p$ and $r \ne p$ a prime that does
not divides  the
order of the Weyl group of $L$.
 Let $L \le A \le \Aut(L)$ with $A/L$ generated by $x$.
Let $M=\Nor_A(L(q_0))$.
If $(A,G,A/M)$ is exceptional, then there is no conjugacy class of
$r$-elements invariant
under $x$.
\end{Corollary}

\begin{proof}  If such a class exists, then by replacing $x$ by $xg$ for
some $g \in L$
we may assume that $\Cen_L(x)$ contains an element $z$ of order $r$.

Then $xg'$ normalizes a Sylow $r$-subgroup $R$ of $\Cen_L(z)$ for some
$g' \in \Cen_L(z)$
and so we may assume that $x$ does (by  exceptionality,
this element has a fixed point and so we may still assume that $x \in M$).
Since the Sylow $r$-subgroup of $L$ is abelian,
$R$ is a full Sylow $r$-subgroup of $L$.
Since $\Cen_L(x) \le M$, $z \in M$ and
so by the previous result, $z$ is in the Frattini subgroup of $R$.

Write $x=uv=vu$ with $u$ an $r$-element and $v$ an $r'$-element.  It
follows that $\Cen_R(v)$ is not contained in the Frattini subgroup of
$R$ and so is not contained in $M$.  By Lemma \ref{ppower}, it follows
that the action is not exceptional.
\end{proof}


\begin{Theorem}\label{Lieplus} Let $L$ be a finite simple group of
  Lie type, defined over $\FFF_q$ and of rank at least $2$ other than
  $Sp_4(2)'=A_6$.  Assume that $(A,B,G,M)$ is such that $L \n A \le
  \Aut(L)$, $L \le G < B \le A$ with $B/G$ cyclic, $M$ a maximal
  subgroup of $A$ not containing $L$, and $(B,G,B/B\cap M)$
  exceptional.  Then $M \cap L$ is a subfield group, the centralizer
  in $L$ of a field automorphism of odd prime order $r$.  Moreover,
  \begin{itemize}
  \item [(a)] $r \ne p$;
  \item [(b)] if $r=3$, then $L$ is of type $\Sp_4$ with $q$ even; and
  \item [(c)] there are no $\Aut(L)$-stable $L$-conjugacy classes of
   $r$-elements.
  \end{itemize}

\end{Theorem}

\begin{proof}

  The proof proceeds by induction in the following manner. We may
  assume that $\gFit(A)=G$ is simple (replace $G$ by $\gFit(A)$ and
  $B$ by $\gen{G,x}$ where $B/G$ is generated by $xG$). We will choose
  a subgroup $J$ of $G$ such that the $G$-conjugacy class of $J$ is
  stable under $\Aut(G)$. Thus, $xg$ normalizes $J$ for some $g \in
  G$.  Replace $x$ by $xg$ to assume that $x$ normalizes $J$. Since
  $xg$ has a fixed point, it is in some conjugate of $M$ (which we
  shall call $M$). Set $B_1=\gen{J,x}$ and $M_1=M \cap B_1$. Since $x
  \in M$ and $x$ normalizes $J$, $B_1=LM_1$. Then $(B_1, J, B_1/M_1)$
  is exceptional (see \ref{desc1}).
  
  For example, if $G=L_n(q)$, $n \ge 3$ and let $x=y^iz^jd$ where $y$
  is the standard field automorphism, $z$ is the standard graph
  automorphism and $d$ induces a diagonal automorphism and centralizes
  $J=\SL_{n-1}(q)$. Then letting notation be as above, $(B_1, L,
  B_1/M_1)$ is exceptional and we may apply induction to show that
  often $M_1 \ge J$. Moreover, in this situation $x$ induces a field
  or graph-field automorphism on $J$ and so we may avoid diagonal
  automorphisms. This will help us to identify $M$ and show that
  exceptionality does not hold.
  
  We shall also use the following fact. If $g \in G$ is contained in
  an $\Aut(G)$-stable conjugacy class, then $xG \cap \Cen_A(g)$ is
  nonempty. Thus, $g$ has a fixed point, as each element in $xG$ has a
  unique fixed point. For example, if $\gFit(A)=L_n(q)$ with $n > 2$,
  then the class of transvections is stable under $\Aut(L_n(q))$.
  Thus, $M$ contains transvections and we may apply results of
  \cite{Kantor:RootElements} and \cite{McL:Some} to help us identify
  the possibilities for $M$. More generally, the class of long root
  elements is usually stable under $\Aut(G)$.  Thus in general $M$
  will contain long root elements and we can use results about
  subgroups containing long root elements.

  Another tool to use to identify a subgroup as a subfield group
  is the result of \cite{BGL} that (with very few exceptions, all
  of rank $1$ in characteristic $2$ or $3$) any
  subgroup of $X(q)$ containing $X(p)$ is a subfield group.

  We also note that conclusion (c) of the theorem implies (a) and (b)
  (see \ref{tab}).

  We now deal with the simple groups of Lie type family by family.

  First consider the case that $G=L_3(q)$.  As noted,
  $M$ contains transvections.  By \cite{FGS}, $M$ is not parabolic.
  So by
  \cite{Kantor:RootElements}, it follows that $M \cap G$ is either a subfield
  group or is the full stabilizer of a non-incident
  unordered point, hyperplane pair (and the
  graph automorphism is present in $A$) or is the unitary subgroup.
  In the last case, $q=p^a$ with $a$ even and
  we apply induction to $L_3(p^{a/2})$ to conclude
  that $M$ contains $L_3(p)$ and so is a subfield group.

  We give a general argument below ruling out the second case.

  So we may assume that $M$ is a subfield group and by maximality
  must be of the form $L_3(q_0)$ with $q=q_0^r$ for some
  prime $r$.  If $x$ does not involve a graph automorphism
  or $p=2$, then
  we may assume that $x$ normalizes $SL_2(q)$ and acts as a
  field automorphism on $SL_2(q)$.  The proof of the $L_2(q)$
  case shows that $r$ is not a divisor of $p(p^2-1)$. On the
  other hand, if stable classes of $r$-elements exist, then
  $r$ is a divisor of $p(p^2-1)$, a contradiction.

  If $x$ involves the graph automorphism (and so not a diagonal automorphism)
  and $p$ is odd, we apply the same argument to $O_3(q)$.

Next suppose that $G=L_n(q)$ with $n > 3$. Write
$q=p^a$. Let
$J=SL_{n-1}(q)$ stabilizing a 1-space and a complementary
hyperplane. Then choose $x$ generating $A/G$ and
normalizing $J$.   By induction on $n$,  and Lemma \ref{desc1},
$SL_{n-1}(p^b) \le M$ for some $b$ and $a/b$ is not divisible
by any prime $r$.

By \cite{Kantor:RootElements}, it follows that the only
possibilities for $M \cap G$ are the stabilizers
of a point or hyperplane,
$J$ or $L_n(p^{a/r})$ with $r$ a prime.
If $M \cap G$ is the stabilizer of a point or hyperplane,
$G$ is $2$-transitive and so exceptionality fails.
If $M \cap G=J$, assume that $M$ is the
stabilizer of the pair $\{U_0,W_0\}$ in $\Omega = A/M$ with $U_0$ a
point and $W_0$ a hyperplane not containing it. Then $M$ has an orbit
on $\Omega$ of size $q^{n-2}(q^{n-1}-1)/(q-1)$ consisting of pairs
$\{U,W\}$ with $U_0 \subset W$ and $U \subset W_0$. This orbit is
visibly transitive also under the action of $G \cap M$, contrary to
exceptionality.

So $M$ is a subfield group.  We need only show that there are no
$\Aut(L)$-stable $L$-classes of elements of order $r$.  Assume
the contrary. By induction,
there are no such classes in $J$, whence the only possibility is
that $r$ is a primitive prime divisor of $p^n-1$ and in particular
$r > n$ and so $r$ satisfies the hypotheses of the previous lemma.
Thus, (c) holds.

Next consider $G=\PSp_{2m}(p^a)$ with $m \ge 2$. First assume that $m=2=p$.
Write $a=ef$ with $e$ odd and $f$ a power of $2$.

If
$e=1$, then $A/G$ has   order a power $2$ and so $M$ contains
a Sylow $2$-subgroup of $A$.  Since parabolic actions are not
exceptional, this implies that  the graph automorphism
must be present and $a=1$ (when $a=1$, $G$ is not the full
symplectic group).  In fact, we do get an example for
$a=1$ (the $A_6$ example).

Thus,   we see that $M$ contains $\PSp_{4}(2^f)$ and
so $M$ is a subfield group (even if $f=1$ as the example is only
when $G=A_6$ not $S_6$).    The argument as for the linear groups
shows that (c) holds.

    Next consider $m > 2=p$.  In this case, $A/G$ is generated by a field
    automorphism and so $Sp_{2m}(2) \le M$.  It follows that $M$ is a subfield
    group, $Sp_{2m}(2^{a/r})$ with $r$ prime.
    By considering the action on $Sp_{2m-2}(2^a)$, we see that
    there no invariant $r$-classes in $G$ unless possibly when $r$ is a divisor
    of $p^{2m}-1$ but not a divisor of $p^{2k}-1$ for $k < m$.  In particular,
    $r > m$ and so does not divide the order of the Weyl group and so the
    previous result applies.

    Now consider $p$ odd.  Essentially the same argument shows that if $M$
    is a subfield group, (c) holds.  So it suffices to prove that $M$ is a
    subfield group.  By induction, it suffices to prove that there are no
    examples when $a$ is a power of $2$ (since the class of $\PSp_{2m}(p^f)$
    is invariant under the full automorphism group and any overgroup of this
    is a subfield group).  So assume $a=f$ is a power
    of $2$.  Then $A/G$ is a $2$-group and so $M$ has odd index in $A$.

    Suppose $m=2$ and $a=1$.  The centralizer of an outer involution is
    $L_2(p^2).2$,  which is maximal in $L$ of even index.  If $a > 1$, the
    only maximal subgroup of odd index is the centralizer of an involution
    (viewing this group as the orthogonal group, it is the stabilizer of a nonsingular
    $1$-space of $+$ type).  There is a unique largest subdegree for $G$
    of size $(q^2-1)(q+1)$
    which is necessarily invariant under the full point stabilizer in $A$.

    Now suppose that $m > 2$.  By induction, $Sp_{2m-2}(p^a) < M$.  In particular,
    $M$ contains transvections.  There are no such proper subgroups
    of odd index  by \cite{Kantor:RootElements}.

Next take $G=U_n(q)$, $q=p^a$ with $n \ge 4$.  We need to show that
$M$ is a subfield group.

Let
$J=SU_{n-1}(q)$ be the derived subgroup of the stabilizer of a
nonsingular $1$-space. We may choose $x$ to be a product of a diagonal
automorphism (possibly trivial) and a field automorphism (of order dividing $2a$)
such that the diagonal automorphism centralizes $J$.
If $J \le M$,  we see (cf.\ \cite{Kantor:RootElements}) that $M$ is the
full stabilizer in $A$ of a non-singular $1$-space $U$. Then $M$ has
an orbit on the set of nonsingular $1$-spaces, consisting of those
$1$-spaces perpendicular to $U$, and $M \cap G$ is visibly transitive
on this, contrary to exceptionality.

So assume that $J$ is not contained in $M$.  We see by induction that
unless $p=2$, $a$ is odd, and  $n=4$, then  $M \cap J$ contains a subfield
group of $J$.  It follows (cf.\ \cite{Kantor:RootElements}) that
$M$ is a subfield group.  Moreover, there are no invariant $r$-classes
contained in $J$ and so arguing as in the linear case shows that (c) holds.

All that remains to consider is the case $n=4$ and $p=2$.  In this case,
there are no diagonal automorphisms.  If $a=1$, the result follows
(since the action would need to be parabolic), and so we see that
$U_4(2) \le M$, whence $M$ is a subfield group.

If $G$ is $\Omega_n(p^a)$ with
$n$ odd, then we may take $n \ge 7$ and $q$ is odd.
Considering subgroups $J_{\ep}$ of type $O_{n-1}^{\ep}(p^a)$ stabilizing
suitable hyperplanes (with $\ep = \pm$), we see that $\Omega_n(p) \le M$
and so $M$ is a subfield group.  The argument as above shows that (c) holds.

Next consider $G=P\Omega^{\ep}_n(p^a)$ with $n$ even, $n \ge 8$. Assume
first that $x$ does not involve a triality automorphism. We show
that $M$ contains $P\Omega^{\ep}_n(p)$ and thus is a subfield group.
By induction, we may reduce to the
case that $a$ is a power of $2$ and then $M$ contains a Sylow $2$-subgroup of
$G$.  Since the action is not parabolic, this forces $p$ odd.
Consider $J_{\delta}$, a
subgroup of type $O_{n-2}^{\delta}$ stabilizing a suitable subspace of
codimension $2$ (with $\delta = \pm$). By induction, $M$ contains
$J(p)$ (the subfield subgroup of $J$). It follows that either $M$
is a subfield group or the stabilizer of a nonsingular $1$-space.
Write $\Omega$ for the set of such $1$-spaces.
Then there is an $M$-orbit on $\Omega$ consisting of all those
elements perpendicular to the one fixed by $M$. However, $G \cap M$ is
transitive on this set, contrary to exceptionality.

Finally consider
$G=\gFit(A)=P\Omega^+_8(p^a)$ with $x$ involving the triality
automorphism.  Since the $3$-cycles are self-centralizing in $S_3$ and
$S_4$, it follows that $x$ must be just a triality or a triality-field
automorphism.    Thus, $G_2(p) \le M$.  If $M$ is a subfield group,
it follows that $r$ does not divide $|G_2(p)|$ and so (c) holds.
It therefore suffices to consider the case that $a=1$ (for the groups over the
prime field form an invariant class under $A$ and if we show
that $P\Omega^+_8(p) \le M$, then $M$ is a subfield group).
So assume that $a=1$.  Thus, $M \ge G_2(p)$ and
$[G:M]$ has order prime to $3$.  In particular, $p \ne 3$.  If
$p=2$, the only possibility for $M$ is the normalizer of a Sylow $3$-subgroup.
This does not contain $G_2(2)$.
If $p > 3$, then we still observe that no proper
subgroup of $G$ containing $G_2(p)$ has index prime to $3$.

This completes the proof for classical groups.

If $G=G_2(p^a)$, assume first that $p \ne 3$
or that $p=3$ and the graph automorphism is not involved.
Then $A/G$ is generated by a field automorphism and so
$G_2(p) \le M$.  By induction, we may assume that $a=r$
and thus $M$ contains a Sylow $r$-subgroup of $G$, whence
$r$ does not divide $|G_2(p)|$ and so there are no invariant
$r$-classes in $G$.

So suppose that $p=3$ and the graph automorphism is involved.
It follows that ${^2}G_2(3) \le M$.  If $a=1$, it also follows
that $M$ has odd index in $A$.  There are no such subgroups
and so $M \ge G_2(3)$ and hence $M$ is a subfield group
$G_2(p^{a/r})$. As above, we may assume that $a=r$.  We need
only show that $r$ is not $2$, $3$ or $7$.  If $r=7$,
we apply Lemma \ref{subfield1}.  If $r=2$, $M$ must contain
a Sylow $2$-subgroup which it does not.  If $r=3$, note
that $|\Cen_G(x^3):\Cen_M(x^3)|$ is divisible by $3$
and so exceptionality cannot hold by Lemma \ref{ppower}.

Suppose that $G={^2}F_4(2^a)'$, with $a$ odd.  If $a=1$,
then $M$ contains a Sylow $2$-subgroup of $G$.  Considering
centralizers, we see that $M \cap G$ contains elements of order
$3$ and $5$.  There is no such proper subgroup of $G$.
Thus, $M$ contains ${^2}F_4(2)'$ and so is a subfield group
${^2}F_4(2^{a/r})$ with $r$ prime.  We claim that
 $r$ does not divide $|{^2}F_4(2)|$, whence (c) holds.

By induction, we may assume that $a=r$.  It follows
that $M$ contains a Sylow $r$-subgroup of $G$ (since
the outer automorphism group has order $r$).
This is not the case (for $r > 3$, apply Lemma \ref{subfield1};
for $r=2$, this is clear and for $r=3$, inspect).

If $G={^3}D_4(p^a)$, then $A/G$ is generated by a field automorphism $x$
and so $G_2(p) \le M$. If $a=1$, then $|A:G|$ has order prime to
$3$ (since $A/G$ has order $3$), which it does not.
Thus, $M$ contains ${^3}D_4(p)$ and is a subfield group
${^3}D_4(p^{a/r})$.  We claim that $r$ does not
 $|G_2(p)|$.  By induction, we may assume that $a=r$.
If $ p \ne r > 3$, apply  Lemma \ref{subfield1}. If $r=3$, then
$M$ contains a Sylow $3$-subgroup of $G$ (since $A/G$
is a $3$-group),  which is visibly not the case.
If $r=2$ or $p$, Lemma \ref{ppower} implies that
$\Cen_M(x^r)$ contains a Sylow $r$-subgroup of
$\Cen_A(x)$, which is not the case.

If $G=F_4(p^a),p \ne 2, E_7(2^a)$ or $E_8(p^a)$, then $G(p) \le
M$, whence $M$ is a subfield group. The standard argument shows
that (c) holds.

If $G=F_4(2^a)$, then we consider $J=F_4(2)$.  It follows that
$J \le M$ (since the outer automorphism group of $J$ has order $2$
and parabolic actions are not exceptional). It follows that
$F_4(2) \le M$, whence $M$ is a subfield group and as usual (c) holds.

If $G=E_6^{\ep}(p^a)$, we will first show that $E_6^{\ep}(p) \le M$, whence
$M$ is a subfield group, $E_6^{\ep}(p^{a/r})$ with $r$ prime.
If $x$ does not involve a diagonal automorphism,
then  $F_4(p)$ is contained
in the centralizer of some element in the coset $xG$.
By induction, we can reduce to the case that $a$ is a power of $2$
(for then there is a unique class of such subgroups).  It follows
that $M$ contains a Sylow $2$-subgroup of $G$. There is no
subgroup of odd index containing $F_4(p)$.

If $x$ does involve a diagonal automorphism, then the
graph automorphism is not involved (since it inverts
the diagonal automorphism).  By choosing different
elements in the coset $xG$, we see that
$M$ contains both $A_5^{\ep}(p)$ and $D_5^{\ep}(p)$,
whence $M$ is a subfield group.

We claim that $r$ does not divide $|F_4(p)|$.
If $p \ne r > 5$, then Lemma \ref{subfield1} applies.
If $p=r$ or $r \le 5$, then the class of the Levi subgroup
$A_5^{\ep}(p^a)$ is invariant and so by induction, there
are no invariant $r$-classes in the Levi subgroup, whence
in fact $r > 5$ and $r \ne p$.

Finally, consider $G=E_7(p^a), p$ odd.  Let $f$ be the largest
power of $2$ dividing $a$.  Suppose that $a=f$.  Then $A/G$ is a
$2$-group and so $M$ has odd index in $A$.  We can choose a Levi
subgroups of $G$ of types $E_6(p^a)$ and  $A_6(p^a)$ whose class
is invariant and so by induction, $M$ contains (up to conjugacy)
these subgroups.  There is no such proper subgroup.  Since the
class of $E_7(p^f)$ is invariant under $A$, by induction, $M$
contains $E_7(p^f)$ and so is the normalizer of $E_7(p^{a/r})$ for
a prime $r$.  If $r$ does not divide the order of the Weyl group,
the results above show that $r$ does not divide the order of
$|E_7(p)|$.  Otherwise $r \le 7$ and we restrict to the Levi
subgroup of type $A_6(p^a)$ to obtain a contradiction.

This completes the proof of the Theorem.
\end{proof}

We finally look at the case that $\gFit(A)$ is sporadic or
alternating.  If $n=6$, then $A_6 = L_2(9)$ and the previous results
apply.

\begin{Theorem}\label{Sporadic}\label{spor} Assume $\gFit(A)$ is $A_n$
  with $6\ne n\ge5$ or is a sporadic simple group. If $M$ is a proper
  subgroup of $A$, then $(A,G,A/M)$ is not exceptional.\end{Theorem}

\begin{proof} Assume that $(A,G,A/M)$ is an exceptional triple. If
  $\gFit(A)=A_n$ with $6 \ne n \ge 5$ then $A=S_n$ and $G=A_n$.  Then
  $xG$ contains transpositions. This implies that the centralizer of a
  transposition is contained in $M$. Since this is a maximal subgroup,
  the only possibility is that $M$ is this centralizer whence has
  orbits of size $2$ and $n-2$ (in the natural degree $n$ permutation
  representation). On the other hand, $xG$ also contains either an
  $n$-cycle or an $(n-1)$-cycle. This is a contradiction.

  If $\gFit(A)=L$ is sporadic, then since $\Out(L)$ has order at most
  $2$, it follows that $A=\Aut(L) \ne G=L$. This implies that $M$ must
  contain a Sylow 2-subgroup of $L$. Moreover, $M$ also contains the
  centralizer of an outer involution.
  There are $12$ of the sporadic simple groups which have outer
  automorphisms. In all but two, namely $M_{12}$ and $HN$,
  there is a centralizer of an outer involution which is a maximal
  subgroup of $G$ and is not of odd index. In the remaining two cases,
  it is straightforward to check that there is no maximal subgroup
  of odd index containing the centralizer of an outer involution.
  \end{proof}

\section{Genus $0$ condition}

In this section, we will use our earlier classification
of exceptional triples to determine the group theoretic
possibilities for arithmetic exceptionality.
We first discuss some general properties of genus $0$ systems.

\subsection{Genus $0$ systems in finite permutation groups}

In the following we assume that $G$ is a transitive permutation group
on a finite set $X$ of size $n$, and let $\cT=(x_1,x_2,\dots,x_r)$ be
a generating system of $G$ with the product relation $x_1x_2\dots
x_r=1$. Recall that for $x\in G$, the index $\ind(x)$ is defined to be
$n$ minus the number of cycles of $x$, and in accordance with the
results in Section \ref{AriPrep}, we say that $\cT$ is a \emph{genus
  $g$ system} if
\[
\sum_{i=1}^r\ind(x_i)=2(n-1+g).
\]
We define the \emph{type $\abs{\cT}$ of} $\cT$ as
$(\abs{x_1},\dots,\abs{x_r})$.

While determining possible genus $0$ systems in specific groups using
(in-)equa\-lities for the index of the elements, one can often rule out
candidates which do fulfill the above index formula by showing that
they cannot generate the given group. For instance, we have the well
known consequence from the theory of polyhedral groups.

\begin{Proposition}\label{poly} Suppose that the type of $\cT$ is
  $(2,2,2,2)$, $(2,3,6)$, $(3,3,3)$, $(2,4,4)$, $(2,3,3)$, $(2,3,4)$,
  or $(2,2,k)$ ($k\ge2$ arbitrary), then $G$ is solvable. If
  $\abs{\cT}=(2,3,5)$, then $G\cong A_5$.\end{Proposition}

The formulation of the main result Theorem \ref{main} makes a
distinction for the genus of a splitting field of $f(X)-t$ over the
field $K(t)$. In terms of the Galois group, this genus is given for
the above group $G$ with respect to the regular action. Note that in
this action, the index of $x\in G$ is given by $\abs{G}(1-1/\abs{x})$,
so this regular action genus $g_{\text{regular}}$ is given by
\[
\sum_{i=1}^r(1-\frac{1}{\abs{x_i}})=
\frac{2}{\abs{G}}(\abs{G}-1+g_{\text{regular}}).
\]

\begin{Definition} We say that $\cT$ is sub--Euclidean or Euclidean, if
  $g_{\text{regular}}$ is $0$ or $1$, respectively.\end{Definition}

A classical result of Zariski (see \cite[II.4]{Magnus:Tess}) basically
classifies the sub--Euclidean and Euclidean systems:

\begin{Proposition}\label{sph_euc}
\begin{itemize}
\item[(a)] If $\cT$ is sub--Euclidean, then one of the following
  holds.\begin{itemize}
\item[(i)] $\abs{\cT}=(n,n)$ for $n\in\NNN$, $G$ is cyclic.
\item[(ii)] $\abs{\cT}=(2,2,k)$ for $k\ge2$, $G$ is dihedral of order
  $2k$.
\item[(iii)] $\abs{\cT}=(2,3,3)$, $G\cong A_4$.
\item[(iv)] $\abs{\cT}=(2,3,4)$, $G\cong S_4$.
\item[(v)] $\abs{\cT}=(2,3,5)$, $G\cong A_5$.
\end{itemize}
\item[(b)] If $\cT$ is Euclidean, then $G''=1$, and $\abs{\cT}$ is
  $(2,2,2,2)$, $(2,3,6)$, $(3,3,3)$, or $(2,4,4)$.
\end{itemize}
\end{Proposition}

As a consequence of part (b) of the previous proposition, we obtain

\begin{Lemma}\label{sph_euc_cor}
  Let $A$ be primitive group acting on a set of size $n$, and $G$ be a
  non--trivial normal subgroup which admits a Euclidean genus $0$
  system $\cT$. Then $n=p^e$ for an odd prime $p$ and $e=1$ or $2$,
  and $G=N\rtimes H$, where $N=C_p^e$ acts regularly, and $H$ is
  cyclic of order the least common multiple of $\abs{\cT}$.
  Furthermore, $n\equiv1\pmod{\abs{H}}$.
\end{Lemma}

\begin{proof}
  By primitivity of $A$, we know that $G$ is transitive. Set $N=G'$.
  By the previous proposition, $N$ is abelian. Clearly
  $N\trianglelefteq A$, and also $N\ne1$, for otherwise $G$ were
  abelian, hence regular, contrary to $\cT$ inducing a genus $1$
  system in the regular action. It follows that $N$ is transitive too,
  and because $N$ is abelian, $N$ is regular and a minimal normal
  subgroup of $A$ (the latter by primitivity of $A$), so $N=C_p^e$ for
  some prime $p$ and $e\ge1$. The results about polyhedral groups in
  \cite[II.4]{Magnus:Tess} show that $N=G'$ is generated by two
  elements. Hence $e\le2$. Let $M$ be a stabilizer in $A$ of a point,
  and $H=G\cap M$. Then $H\cong G/G'$ is an abelian normal subgroup of
  $M$. If $H$ acts irreducibly on $N$, then $H$ is cyclic and
  semi--regular on $N\setminus\{1\}$. If $H$ is reducible on $N$, and
  $N_0<N$ is $H$--invariant, then $N_0^m$ is $H$--invariant for each
  $m\in M$. But primitivity of $A$ implies that $M$ is irreducible on
  $N$, and from that we get easily that $H$ acts a group of scalars on
  $N$, if we identify $N$ with $\FFF_p^e$ and $M$ as a subgroup of
  $\GL_e(p)$. So $H$ is cyclic and semi--regular in either case, and
  the claim follows.
\end{proof}

The determination of genus $0$ systems relies on the determination of
the exact values or lower bounds of the index. If $x\in G$, then
clearly
\[
n-\ind(x)=\frac{1}{\abs{x}}\sum_{i=1}^{\abs{x}}\chi(x^i),
\]
where $\chi(y)$ denotes the number of fixed points of $y$. So if we
know, for instance, the permutation character and character table of
$G$, then we can compute the index function. In small groups, which we
do not want to treat by hand, we make frequent use of the computer
algebra system GAP \cite{GAP}.

The next is a reduction result for genus $0$ systems which preserve a
product structure. For a set $X$ denote by $S(X)$ the symmetric group
on $X$.

\begin{Lemma}\label{P} Let $X=Y\times Z$ be a finite set with $\abs{Y}=y$
  and $\abs{Z}=z$ with $y,z>1$. Let $G$ be a transitive subgroup of
  $S(Y)\times S(Z)$. Let $x_i=(y_i,z_i)$ be a genus $0$ system for $G$
  (where $y_i\in S(Y)$ and $z_i\in S(Z)$), that is
  $\sum\ind(x_i)=2(yz-1)$. Then
\begin{itemize}
\item[(a)] $\sum\ind(y_i,Y)=2(y-1)$;
\item[(b)] Either some $x_i$ must have no fixed points on $X$, or
$G''=1$ and if $d_i$ is the order of $x_i$, then $y_i$ and $z_i$ each
have order $d_i$ with $\sum(1-1/d_i)=2$;
\item[(c)] $\sum \chi(y_i,Y)\ind(z_i,Z)<2z$ (where $\chi(y_i,Y)$
denotes the number of fixed points of $y_i$ on $Y$).
\end{itemize}
\end{Lemma}

\begin{proof} It is straightforward to compute
  (cf.~\cite[3.2(b)]{GurNeu}) that
  \begin{equation}\label{indpro}\ind(x_i)\ge
    z\cdot\ind(y_i,Y)+\chi(y_i,Y)\ind(z_i,Z).\end{equation}

  In particular, $\sum\ind(x_i)\ge z(\sum\ind(y_i,Y))$. If (a) does
  not hold, then $\sum\ind(y_i,Y)\ge2y$, whence $\sum\ind(x_i)\ge2yz$,
  contrary to the hypothesis.

  Then \eqref{indpro} yields that $2(yz-1)\ge 2z(y-1)+\sum
  \chi(y_i,Y)\ind(z_i,Z)$, and (c) holds.

  Suppose that each $\chi(x_i)>0$. Then $\chi(y_i,Y)>0$ for each $i$
  (since the number of fixed points for $x_i$ is just the product of
  the number of fixed points for $y_i$ and for $z_i$).

  If $\chi(y_i,Y)>1$ for some $i$, then \eqref{indpro} gives a
  contradiction. Similarly, we see that $\chi(z_i,Z)=1$ for each $i$.

  Moreover, by \cite[3.2(a)]{GurNeu}, it follows that each $x_i$ has a
  unique fixed point and acts semi-regularly on all other points
  (i.e.\ all other orbits have size $d_i$, the order of $x_i$). In
  particular, $y_i$ and $z_i$ each have order $d_i$ as well. Thus,
  $2yz-2=\sum (d_i-1)(yz-1)/d_i$, whence $\sum(d_i-1)/d_i=2$, so
  $G''=1$ by Proposition \ref{sph_euc}(b), and $(x_1,\dots,x_r)$ is
  Euclidean.\end{proof}

\begin{Proposition}[\protect{\cite[Reduction Proposition,
    p.~4279]{Magaard:Sporadic}}]\label{Kay} Let $G$ be a nonsolvable
  primitive group acting on $X:=Y^t$ with $t>1$ and preserving the
  product structure on $X$. Assume that $|Y| \ge 100$ and if $g \in G$
  stabilizes $Y$ and acts nontrivially on $Y$, then $|\chi(g,Y)|/|Y|
  \le 1/7$.  Then there is no genus $0$ system of generators for the
  action of $G$ on $X$.\end{Proposition}

\begin{proof} Suppose $\cT$ is a genus $0$ system for $G$ for the
  action on $X$. First, we apply \cite{GurNeu} to observe that if $T$
  is the normal subgroup of $G$ which acts coordinatewise on $X$, then
  $G/T$ is either cyclic of order $a$ or dihedral of order $2a$ with
  $a \le 4$.  Now the argument in \cite{Magaard:Sporadic} shows that
  $\cT$ consists of $3$ elements. Moreover, if each element of $\cT$
  acts semiregularly as a permutation of the coordinates, then the
  argument in \cite{Magaard:Sporadic} shows that every nontrivial
  power of each element in $\cT$ fixes at most $|X|/49$ points. This
  implies that $\cT$ must contain elements of order $2$, $3$ and $7$,
  whence $G$ is perfect. This contradicts that $G/T$ is a nontrivial
  solvable group.

  The remaining cases are when $G/T$ is dihedral of order $6$ or $8$
  and $t=3$ or $4$ respectively. The latter case is eliminated by
  using \cite[Lemma K]{Magaard:Sporadic} to conclude that $\sum
  (d_i-1)/d_i < 2 + 149/2400$, where the $d_i$ are the orders of the
  elements of $\cT$.  Then this implies that $\cT$ has type $(2,3,b)$ with
  $7 \le b \le 9$ or $(2,4,5)$.  However, if $G/T$ is dihedral of
  order $8$, every element in $\cT$ has even order and we obtain a
  contradiction.

  Similarly in the case $\bar\cT$ (the image of $\cT$ in $G/T$) is of
  type $(2,2,3)$, we see that the only possibilities for the type of
  $\cT$ are $(2,3,8)$ or $(2,3,10)$.  In the latter case, the third
  element must be of the form an involution on $2$ components and an
  element of order $5$ on the third --- it is easily seen that these
  elements do not generate (because $3$ conjugates of the square of
  the third element generate and these $3$ conjugates generate a
  solvable group --- namely, a dihedral group on each copy of $Y$).

  The remaining case is when $\cT$ has type $(2,3,8)$.  Again we see
  that $3$ conjugates of the third element squared generate a
  nonsolvable group.  This forces at least $2$ of the coordinates of
  the square of this element to have order $4$. Then every nontrivial
  power of this element fixes at most $|X|/49$ points, whence the
  genus is positive.
\end{proof}

An elementary consequence of \cite[Lemma 2.3, Prop.~2.4]{GT} allows us
to exclude the existence of genus $0$ systems for various actions.

\begin{Lemma}
  \label{n/85}
  Let $G$ be a non--solvable permutation group of degree $n$, which is
  not isomorphic to $A_5$. If $G$ admits a genus $0$ system, then
  there is $1\ne g\in G$ such that $g$ fixes more than $n/85$ points.
\end{Lemma}

\subsection{Diagonal action}

We need the following extension of Aschbacher's result
\cite{Asch:GT}.

\begin{Theorem}\label{diag} Suppose $A$ acts primitively and
faithfully on $\Omega$ and is of diagonal type. Let $G$ be a
nontrivial normal subgroup of $A$. Suppose that $G$ admits a genus $0$
system. Then $\gFit(A)=A_5\times A_5$. In particular, $(A,G,\Omega)$ is
not exceptional.\end{Theorem}

\begin{proof} Let $E=\gFit(A)$. By the Aschbacher--O'Nan--Scott theorem,
$E$ is the unique minimal normal subgroup of $A$, whence is contained
in $G$. Since $\gFit(G)$ is $A$-invariant, it follows that
$\gFit(G)=\gFit(A)$.

Write $E=L_1\times\dots\times L_t$, where the $L_i$ are the components
of $A$ with $L_i\cong L$ a nonabelian simple group. Let $M$ be a point
stabilizer in $A$. Then $H=G\cap M$ is a point stabilizer in $G$.
There is an $A$-invariant partition $\Delta$ of $I=\{1,\ldots,t\}$ so
that $H\cap E=\prod_{\delta\in\Delta}L(\delta)$, where $L(\delta)\cong
L$ is a full diagonal subgroup of $\prod_{i\in\delta}L_i$. Of course,
each $\delta$ has the same cardinality and has size greater than $1$
(but the partition $I$ is allowed). Since $A$ is primitive, it follows
that $A$ is transitive on the set of the $L_i$ and moreover,
$A_{\delta}$, the stabilizer of the block $\delta$, acts primitively
on $\delta$.

Note that $M=\Nor_A(E\cap H)$. If $\gamma$ is a subset of $I$, let
$E(\gamma)$ denote the direct product of the $L_i$ as $i$ ranges over
$\gamma$.

We consider several cases.

Suppose first that $G_{\delta}$ is not transitive on $\delta$. Since
$A_{\delta}$ is primitive on $\delta$, this implies that $G_{\delta}$
is trivial on $\delta$. Since $G$ is normal in $A$, it follows that
this holds for all $\delta$. Let $F(\delta)$ be a subgroup of
$E(\delta)$ containing $L(\delta)$ which is maximal subject to being
$H_{\delta}$--invariant. Note also that $F(\delta)$ must miss at least
$2$ of the components in $E(\delta)$ (otherwise,
$F(\delta)=E(\delta)$).

Moreover, we may do this in such a way so that if $g\delta=\delta'$,
then $F(\delta')=hF(\delta)h^{-1}$ where $h\in H$ maps $\delta$ to
$\delta'$ (since $G=EH$ and $E$ acts trivially on $I$, there is no
loss in taking $g\in H$ -- moreover, $h$ is uniquely determined up to
an element in $H_{\delta}$ which normalizes $F(\delta)$ so that this
conjugate is uniquely determined). Let $F$ be the direct product of
the various $F(\delta)$. This is $H$-invariant and so we can consider
the subgroup $HF$. Let $S$ be a maximal subgroup of $G$ containing
$HF$. Then $S\cap E(\delta)=F(\delta)$ or $E(\delta)$ for each
$\delta$. Since $G=HE$, $S$ cannot contain $E$.  So, $E(\delta)$ acts
nontrivially on the cosets of $F$ for some $\delta$. In particular, we
see that there are $2$ components acting nontrivially on the cosets of
$F$ which are not $G$-conjugate. Thus, if $C=\core_G (S)$, we see that
$\gFit(G/C)$ is not a minimal normal subgroup, but that $G/C$ acting
on the cosets of $F$ has genus $0$ (by L\"uroth's theorem). This
contradicts \cite{Shih}.

Suppose next that $G$ is transitive on $I$. Let $F$ be a maximal
subgroup containing $H$. Then $F$ cannot contain any components of $G$
(since $G=HE$ and therefore $H$ is transitive on the components --
thus, if $F$ contains a component, it contains $E$ and so $HE=G$).
Since $E$ is the unique minimal normal subgroup of $G$ (under the
assumption that $G$ is transitive on $I$), $G$ acts faithfully on the
cosets of $F$. So $G$ is primitive and faithful on the cosets of $F$.
Now apply \cite{Asch:GT} to conclude that $\gFit(G)=A_5\times A_5$ as
desired.

So we may assume that $G$ is not transitive on $I$ but that
$G_{\delta}$ is transitive on $\delta$. If $I'$ is a $G$-orbit of $I$,
it follows that $I'$ is a union of certain subsets of $\Delta$.
Arguing as in the previous paragraph, we see that $L=A_5$ and that
$I'$ has cardinality $2$. Thus, $\Delta$ is precisely the set of
$G$-orbits on $I$ each of cardinality $2$. Denote these orbits by
$\delta_i$ for $i=1,\dots,s$ (and so $t=2s$).

It follows that $G\le G_1\times\dots G_s$ where $t=2s$ and
$G_i=G/\Cen_G(E(\delta_i))$. Then $\gFit(G_i)=A_5\times A_5\le G$ and we
may identify $\Omega$ with $\Omega_1\times\dots\times\Omega_s$.

Write $x_i=(y_{i1},\dots,y_{is})$. We may now apply Lemma \ref{P} and
use the specific information given in \cite[Section~19]{Asch:GT} to
conclude that if $s>1$, then the genus is larger than zero. So $s=1$
and $\gFit(G)=A_5\times A_5$ as desired.

Finally, since $G$ is transitive on the components, it follows that
$A/G$ embeds in $\Out(A_5)$. In particular, $A/G$ has order at most
2. So $(A,G,\Omega)$ exceptional implies that $|\Omega|$ would have
odd order. However, $|\Omega|=60$. This completes the
proof.\end{proof}

\subsection{Product action}\label{prod2}

In this section we will use some results of the almost simple case
given in the next section.

\begin{Theorem} Suppose $\gFit(A) = L_1 \times \ldots \times L_t$
  with $L_i \cong L$ a nonabelian simple group with $t > 1$.  Let $M$
  be a maximal subgroup of $A$ such that $M$ does not contain
  $\gFit(A)$ and $M \cap \gFit(A)= R_1 \times \ldots \times R_t$ where
  $R_i = M \cap L_i \ne 1$.  Let $G$ be a normal subgroup of $A$ with
  $A=GM$.  If $(G,A/M)$ has genus $0$, then $(A,G,A/M)$ is not
  arithmetically exceptional.
\end{Theorem}

\begin{proof} Set $X=A/M$. Identify $X=Y\times\dots\times Y$ ($t$ copies
  of $Y$ where $Y\cong L_1/(L_1\cap M)$). So $n=\ell^t$ where
  $\ell=|Y|$. It follows by \cite{GT} that some nontrivial element in
  $N_A(L_1)/C_A(L_1)$ fixes more than $\ell/85$ points on $Y$.
  
  Assume the result is false. By Lemma \ref{L:prod} $(A_1,L,Y)$ is
  arithmetically exceptional (in the notation of the lemma). By
  Theorems \ref{Lie85} and \ref{spor}, it follows that $L\cong
  L_2(q)$. We now apply Theorem \ref{L2} to conclude that $L_1$ and
  $L_1 \cap M$ must be as given in Theorem \ref{L2}(a).

  Let $s$ be the size of a nontrivial orbit of $G$ on the components
  (note that they are all of equal size since $G$ is normal in $A$).
  Exclude the cases $p^a=8$ or $9$ for the moment. Then $s$ cannot be
  greater than $1$ by Lemma \ref{P}, Proposition \ref{Kay}, and
  Theorem \ref{L2}(d).

  So $L=L_2(q)$ with $q=8$ or $9$ and $\ell=28$ or $45$ respectively.

  Let $T$ be the subgroup of $G$ which normalizes each component.  We
  first show that the minimal genus of $T$ on $X$ is large.

  We first consider the action of $T$ on $W:=Y\times Y$. Let $\cT$ be a
  generating set for $T$ with the product of the elements being $1$.

  First suppose that $\cT$ contains an element of the form $(1,g)$. Let
  $\cT'$ be the remaining elements of $\cT$.  By restricting our attention
  to the first component, we see that $\sum_{\cT'} \ind(h_1,h_2) \ge
  2\ell^2 - 2\ell$.  Thus, $\sum_\cT \ind(h_1, h_2) \ge 2\ell^2 - 2\ell
  + \ind_Y(g)\ell$.  Since $\ind_Y(g) \ge 12 $ (case $L=L_2(8)$) or
  $18$ (case $L=L_2(9)$), it follows that $g(T,W,\cT) - 1 \ge 5 \ell$ or
  $8 \ell$.

  Now assume that $\cT$ contains no such elements.  An inspection of the
  possible elements and using the fact that in each component, the
  orders of the elements in $\cT$ are not Euclidean or sub--Euclidean,
  c.f.~Proposition~\ref{sph_euc}, we see that $g(T,W,\cT) > 0$.

  So we may assume that $G\ne T$. Again, by Lemma \ref{P}, this
  implies that $g(G,W,\cT)=0$ for $W=Y^s$ for some generating set $\cT$ of
  $G$ (with the product of the elements being $1$).  We show that this
  cannot happen.  For ease of notation, we assume that $G$ is
  transitive on the components (this amounts to passing to a quotient
  group -- this quotient group need not be exceptional but we will no
  longer use that hypothesis).

  By \cite{GurNeu}, we know that the action of $G$ on the orbit of
  components is cyclic of order $a$ or dihedral of order $2a$ with $a
  \le 4$.  Moreover (by reordering if necessary), $\bar\cT$ has type
  $(a,a)$ or $(2,2,a)$ respectively.

  Let $\cT=(x_1, \ldots, x_r)$.  First assume that $G/T$ is cyclic.  We
  follow the notation in \cite{GurNeu}. There is a generating set $F$
  for $T$ consisting (up to conjugation) of $x_1^a, x_2^a$ and $a$
  conjugates of $x_i$ for each $i > 2$.

  First consider the case that $L=L_2(9)$ and $\ell=45$.

  Note that each power of $x_i$ for $i <3$ fixes at most $n/5^a$
  points.

  Suppose that $a=4$.  If $x_3$ has order at least 3, then $\ind(x_3)
  \ge (2/3)n$ and $\ind(x_i)\ge 3n/4-45/2-45^2/4$ for $i=1,2$. It
  follows easily that the genus is positive.  Similarly, if $r \ge 4$,
  the genus is positive.  So we may assume that $x_3$ has order 2.
  Then we may assume that $x_1$ has order greater than $4$ (since $G$
  is not Euclidean). Then $\ind(x_1) \ge (7/8)(624/625)n$ and again we
  see that the genus is positive.

  Now consider the case that $a=3$.  If $x_3$ has order at least $4$,
  then $\ind(x_i) \ge (2/3)(n-45)$ for $i \le 2$ and $\ind(x_3) \ge
  (32/45)n$, whence the genus is positive.  Similarly, the result
  follows if $r \ge 4$ (since then $\ind(x_3) + \ind(x_4) \ge
  (4/5)n$).  So assume that $r =3$.  If $x_3$ has order $3$, then we
  may assume that $x_1$ has order at least $6$. Then, $\ind(x_1) \ge
  (5/6)(124/125)n$, $\ind(x_2) \ge (2/3)(n-45)$ and $\ind(x_3) \ge
  (2/3)n$, whence the result follows.  If $x_3$ has order $2$, then
  either both $x_1$ and $x_2$ have order at least $6$ or we may assume
  that $x_1$ has order at least $9$. Since $\ind(x_3) \ge (2/5)n$, it
  is straightforward to compute as above that the genus is positive.

Next, consider the case $a=2$.  Then $\ind(x_i) \ge (n-45)/2$ for $i
\le 2$ and $\ind(x_i) \ge (2/5)n$ for all $i$.  It follows that $r \le
4$.  If $r=4$ and say $x_4$ is not an involution, then $\ind(x_4) \ge
(2/3)n$ and again we see that the genus is positive.  If $x_3$ and
$x_4$ are involutions, then we may assume that $x_1$ is not.  Then
$\ind(x_1) \ge 3n/4 - 171/4$ and again we see that the genus is
positive.  So $r=3$.  Since $G$ is not dihedral, we may assume that
$x_1$ is not an involution.  Straightforward computations show that
the genus is positive in all cases (first assume that $x_3$ has order
$>4$; if it has order $3$ or $4$, then either $x_1$ has order $> 6$ or
$4$ respectively or both $x_1$ and $x_2$ have order $4$; if $x_3$ has
order $2$, then we may assume that $x_1$ has order $ \ge 6$ and $x_2$
has order $\ge 4$).

Now assume that $G/T$ is dihedral of order $2a$ with $2 \le a \le 4$.
It is straightforward as above to see that $r=3$.  Since $G$ is not
dihedral, it follows that we may assume that one of $x_1$ and $x_2$ is
not an involution.  If $a=2$, then every nontrivial power of $x_i$
fixes at most $n/25$ points.  It follows easily from this observation
that the genus is positive.  The same is true for $a=3$ or $4$ unless
$t=a$.  A straightforward computation of the possibilities shows that
the genus is always positive (one needs to use the facts that $G$ is
neither Euclidean nor sub--Euclidean and that $T$ can be generated by
$x_1^2$, a conjugate of $x_2^2$ and $a$ conjugates of $x_3^2$).

A similar argument shows that the genus is positive when $L=L_2(8)$
(this case is much easier since each nontrivial element of $\Aut(L)$
fixes $0,1$ or $4$ points of $Y$).\end{proof}

\subsection{Almost simple groups}\label{as2}

In this section, we consider the exceptional examples that
are almost simple and determine which admit genus $0$ systems.

We will also obtain some additional information which we used in
Section \ref{prod2}.

The main result is:

\begin{Theorem}\label{AS} Let $A$ be a primitive permutation
  group of almost simple type and of degree $n$. Let $G$ be a
  nontrivial normal subgroup of $A$, where $G$ admits a genus $0$
  system. Suppose there is a subgroup $B$ of $A$ containing $A$ such
  that $B/G$ is cyclic, and $(B,G)$ is exceptional. Then one of the
  following holds:
\begin{itemize}
 \item[(a)] $n=28$, $G=\PSL_2(8)$ of
    types $(2,3,7)$, $(2,3,9)$, or $(2,2,2,3)$, and $B=A=\PgL_2(8)$.
  \item[(b)] $n=45$, $G=\PSL_2(9)$ of type $(2,4,5)$, $B=\M10$, and
    $A=B$ or $\PgL_2(9)$.
\end{itemize}
\end{Theorem}

\begin{Theorem}\label{L2} Let $A$ be a group with $\gFit(A)=L_2(p^a)$.
  Let $G$ be a normal subgroup of $A$ and $M$ a maximal subgroup of
  $A$.  Assume that $(A,G,A/M)$ is arithmetically exceptional. Set
  $n=[A:M]$.  Let $D$ be any subgroup of $A$ with $\gFit(A)=\gFit(D)$.
\begin{itemize}
\item[(a)] If there exists $1\ne x\in A$ such that $x$ fixes more than
  $n/85$ points, then one of the following holds:
\begin{itemize}
\item[(i)] $p^a=8$, $n=28$, $A=\Aut(G)$;
\item[(ii)] $p^a=9$, $n=45$, $G$ is simple, $A=\M10$ or $\Aut(G)$;
\item[(iii)] $p^a=32$, $n=496$, $G$ is simple and $A=\Aut(G)$;
\item[(iv)] $p^a=25$, $49$, $81$, or $121$, $n=p^a(p^a+1)/2$;
\item[(v)] $p^a=27$, $n=351$ or $819$, and $G$ is simple;
\end{itemize}
\item[(b)] If $(G,A/M))$ has genus $0$, then one of the following
  holds:
\begin{itemize}
\item[(i)] $p^a=8$, $[A:M]=28$, $A=\Aut(G)$, $G$ is simple and the
  ramification type is $(2,3,7)$, $(2,3,9)$, or $(2,2,2,3)$;
\item[(ii)] $p^a=9$, $[A:M]=45$, $G$ is simple, $A=\M10$ or $\Aut(G)$
  and the ramification type is $(2,4,5)$.
\end{itemize}
\item[(c)] If $(D,A/M)$ has genus $0$, then $p^a=8$ or $9$. If $D$ is
  simple, then the ramification type is as given in (b). Otherwise, if
  $p^a=8$, $D=\Aut(G)$, then the ramification type is $(2,3,9)$; and
  if $p^a=9$, the ramification type is $(2,3,8)$ if $D=\PGL_2(9)$, or
  $(2,4,6)$ if $D=\Aut(G)$;
\item[(d)] In case (a), subcases (iii)--(v), every nontrivial element
  in $\Aut(\gFit(A))$ fixes at most $n/7$ points.
\end{itemize}
\end{Theorem}

\begin{proof} We apply the  result of the exceptionality
section and consider the various
  cases. Let $L=L_2(p^a)$ and $H=M\cap L$. Note that if each
  nontrivial element of $A$ fixes at most $n/85$ points, then there is
  no genus $0$ system for this action by Lemma \ref{n/85}. Thus,
  assuming (a), one need only check the cases listed in (a) when
  proving (b).
  
  If $q>113$, then each non--trivial element in $\Aut(L)$ fixes at
  most $n/85$ points (see \cite{GT}), except for those cases listed in
  Table 1 in \cite[p.267]{LS:MinDeg}. The only such case compatible
  with the previous theorem is $p^a=121$ with $n=7381$.

First consider the case that $H=L_2(p^{a/b})$ with $a/b$ a prime. The
only possibilities are $p=2$ with $a=5$ or $p=3=a$. In the first case,
$M$ is not maximal in $A$. In the second case, we verify easily using
GAP that the non--trivial element with the most fixed points is an
involution with $13$ fixed points, and that no group $D$ has a genus
$0$ system.

Next suppose that $p=2$ and $H$ is dihedral of order $2(2^a+1)$.  Then
$n=2^a(2^a-1)/2$. Since $2^a<113$ and $a$ is odd, it follows that
$a=3$ or $5$ and $G=L$.

First suppose that $a=5$. The following table gives the fixed point
numbers and indices of the elements $x\in A$ depending on $\abs{x}$.
This table can be easily extracted from the information in
\cite{ATLAS}, or be computed using GAP.

\[
\begin{array}{c|c|c|c|c|c}
\abs{x} &   1 &    2 &    3 &    5 &   >5 \\
\hline
\chi(x) & 496 &   16 &    1 &    1 & \le1 \\
\hline
\ind(x) &   0 &  240 &  330 &  396 & \ge444
\end{array}
\]

>From this table we see that the only possibility of a genus $0$ system
for $D$ would consist of $3$ elements of order $3$. But such a system
generates a solvable group, see Proposition \ref{poly}. Also, we see
that each nontrivial $x\in A$ fixes fewer than $n/7$ points.

If $a=3$, elements of odd order fix either $0$ or $1$ of the $28$
points and involutions fix $4$ out of $28$ points. Moreover, elements
of order $7$ have no fixed points. It follows easily that the only
possible types of genus $0$ systems are $(2,3,7)$, $(2,3,9)$, or
$(2,2,2,3)$, as in the conclusion of the theorem. Note that the type
$(2,3,9)$ can occur with $D=L$ and $D=\Aut(L)$ as claimed.

Now consider case (b)(iii) of Theorem \ref{L2ex}. Thus, $H$ is
dihedral of order $p^a-1$, $p$ is odd and $a$ is even. In particular,
$n=p^a(p^a+1)/2$. From Table 1 in \cite[p.267]{LS:MinDeg} we obtain
$p^a=9$, $25$, $49$ or $81$, or $121$.

In all but the case $p^a=81$, since $B/G$ must be generated by a
field--diagonal automorphism only, for any $B$ with $B/G$ cyclic and
$(B,G,A/M)$ exceptional, it follows that $G=L$. If $p^a=81$, then
$G=L$ or $G/L$ has order $2$ generated by a field automorphism.

If $p^a=9$, then $n=45$. One computes easily using GAP that the only
genus $0$ systems are as stated in the theorem. The arithmetic
exceptionality forces $\M10\le A$.

If $p^a=25$, then $n=325$. One computes using GAP that there are no
genus $0$ systems for $D$.

If $p^a=49$, then $n=1225$. We use the character table of the $3$
different quadratic extensions of $L_2(49)$ in order to compute the
indices of the elements in $\PgL_2(49)$.

\[
\begin{array}{c|c|c|c|c|c|c|c|c|c}
\abs{x} & 1 & 2 & 2 & 3 & 4 & 5 & 6 & 6 & \ge7 \\
\hline
\chi(x) & 1225 & 49 & 25 & 1 & 1 & 1 & 1 & 1 & \le1 \\
\hline
\ind(x) & 0 & 588 & 600 & 816 & 912 & 980 & 1012 & 1016 & \ge1050
\end{array}
\]

First note that a possible genus $0$ system consists of only $3$
generators.  As it may contain at most one involution, we see that the
only sum of three indices giving $2(1225-1)=2448$ comes from three
elements of order $3$. But that leads to a solvable group, see
Proposition \ref{poly}.

If $p^a=81$, then $n=3321$. Again, using character tables of the
cyclic extensions of $L_2(81)$, we obtain similarly as above

\[
\begin{array}{c|c|c|c|c|c|c|c|c|c}
\abs{x} & 1 & 2 & 2 & 3 & 4 & 4 & 5 & 6 & \ge8 \\
\hline
\chi(x) & 3321 & 81 & 41 & 0 & 9 & 1 & 1 & 0 & \le1 \\
\hline
\ind(x) & 0 & 1620 & 1640 & 2214 & 2466 & 2480 & 2656 & 2754 &
\ge2900
\end{array}
\]

We see that no genus $0$ system exists for a group between $L=L_2(81)$
and $\Aut(L)$ in the given action.

If $p^a=121$, then $n=7381$. Using GAP (or computation by hand), we
obtain

\[
\begin{array}{c|c|c|c|c|c|c|c|c}
\abs{x} & 2 & 2 & 3 & 4 & 5 & 6 & 6 & >6 \\
\hline
\ind(x) & 3630 & 3660 & 4920 & 5520 & 5904 & 6130 & 6140 & \ge6450
\end{array}
\]

>From $3\cdot3630+4920>14760=2(n-1)$ we get that there are only $3$
generators in a genus $0$ system. At least two elements have order
$\ge3$, and from that we obtain a match in the index relation only if
all elements have order $3$, a case which does not occur.

Finally, consider case (b)(iv) of the previous theorem.  Again, from
\cite[Theorem 1]{LS:MinDeg}, $3^a < 113$, so $a=3$ and $n=351$.

We compute a table as above of the relevant data for $\Aut(L_2(27))$.

\[
\begin{array}{c|c|c|c|c|c|c|c|c|c|c}
\abs{x} & 1 & 2 & 2 & 3 & 3 & 4 & 6 & 6 & 7 & \ge9 \\
\hline
\chi(x) & 351 & 15 & 13 & 3 & 0 & 1 & 3 & 1 & 1 & \le1 \\
\hline
\ind(x) & 0 & 168 & 169 & 232 & 234 & 259 & 288 & 289 & 300 & \ge312
\end{array}
\]

There are two matches for the index relation. One belongs to three
elements of order $3$, which is not possible by Proposition
\ref{poly}. The other match is $2(351-1)=168+232+300$, suggesting a
$(2,3,7)$ system. However, such a system clearly cannot generate a
group with a non--trivial abelian quotient, so it would have to
generate $L_2(27)$. However, the element of order $3$ of index $232$
is the field automorphism of $L_2(27)$, a contradiction.

Note that in all cases, the proofs and tables show that the number of
fixed points is at most $n/7$.
\end{proof}

\begin{Remark*} In all examples described in the previous
  result, $L\cap M$ is maximal in $L$.\end{Remark*}

The next theorem  is an immediate consequence of our earlier results
from the exceptionality section.

\begin{Theorem}\label{Lie} Let $L$ be a finite simple group of
  Lie type, defined over $\FFF_q$ and not of type $L_2$.  Assume that
  $(A,B,G,M)$ is such that $L \n A \le \Aut(L)$, $L \le G < B \le A$
  with $B/G$ cyclic, $M$ a maximal subgroup of $A$ not containing $L$,
  and $(B,G,B/B\cap M)$ exceptional.  Then either $q > 113$ or
  $(L,B,G,M)$ is one of the following:\begin{itemize}
  \item[(a)] $(\Sz(8), \Aut(\Sz(8)), \Sz(8), \Sz(2) \times 3)$;
  \item[(b)] $(\Sz(32), \Aut(\Sz(32)), \Sz(32), 25:20)$;
  \item[(c)] $(U_3(8), B, G, \Nor_A(9 \times 3))$, with $U_3(8)$ of
    index $1$ or $3$ in $G$, and $|B/G| = 3$;
  \item[(d)] $(L^{\ep}_3(32), B, G, \Nor_A(L^{\ep}_3(2)))$, with
    $L^{\ep}_3(32)$ of index at most $3$ in $G$, and $5$ dividing
    $|B/G|$;
  \item[(e)] $(G_2(32), \Aut(G_2(32)), G_2(32), G_2(2)\times 5)$;
  \item[(f)] $(\Sp_4(q), \Aut(\Sp_4(q)), G, \Nor_A(\Sp_4(q^{1/3})))$,
    with $q=8$ or $q=64$ and $\Sp_4(q)$ in $G$ of index $1$ in the
    former and at most $2$ in the latter case;
  \item[(g)] $(\Omega_8^+(32), \Omega_8^+(32).15, \Omega_8^+(32),
    \Nor_A(\Omega _8^+(2)))$;
  \item[(h)] $({^3}D_4(32), B, G, \Nor_A({^3}D_4(2)))$, with
    ${^3}D_4(32)$ of index $1$ or $3$ in $G$, and $5$ dividing
    $|B/G|$;
  \item[(i)] $(U_3(32), B, G, M)$, with $M \cap U_3(32)$ the stabilizer
  of a system of imprimitivity.
\end{itemize}
\end{Theorem}

We now show that none of the remaining examples of the preceding
theorem has a genus $0$ system.

\begin{Theorem}\label{Lie85} Let $A$ be a group with $\gFit(A)$ a
  simple group $L$ of Lie type not isomorphic to $L_2(q)$.  Let $G$ be
  a normal subgroup of $A$ and $M$ a subgroup of index $n$ in $A$ such
  that $M=\Nor_A(M\cap\gFit(A))$. Assume that $(A,G,A/M)$ is
  arithmetically exceptional.  Then, no nontrivial element of $A$
  fixes more than $n/85$ points on $A/M$.
\end{Theorem}

\begin{proof} Write $\chi$ for the permutation character of the action of
  $G$ on $A/M$. By \cite{LS:MinDeg}, if $q > 113$ then the action is
  not of genus zero by Lemma \ref{n/85}, because the fixed point ratio
  $\chi(g)/\chi(1)$ is $\le 1/85$ for $1\ne g\in G$. Thus $q \le 113$,
  so by Theorem \ref{Lie}, we only have to deal with the groups listed
  in the conclusion of Theorem \ref{Lie}. We shall show that
  $\chi(g)/\chi(1)\le 1/85$ for $1 \ne g \in A$. We use Lemma
  \ref{chig}.

  (a) If $L=\Sz(8)$ then $\chi(1)=1456$; the only elements of prime
  order in $L$ fixing points of $A/M$ have order $2$ or $5$ (a unique
  class of each), and they fix $16 (=64/4)$ and $1$ points
  respectively.  An outer automorphism has a unique fixed point.

  (b) If $L=\Sz(32)$ with $G \cap M = 25:4$ then $\chi(1)=325376$;
  the only elements of prime order in $G$ fixing points of $A/M$ have
  order $2$ or $5$, and they fix $256 (=2^{10-2})$ and $1$ points
  respectively.  All outer automorphisms have a unique fixed point.

  (c) If $L = U_3(8)$, then the degree is $34048$.  One computes
  using GAP the corresponding permutation character for each cyclic
  extension of $U_3(8)$ to see that every character value is at most
  $256$.  In particular, no element fixes more than $n/85$ points and
  so there are no genus $0$ systems by Lemma \ref{n/85}.

  (d) If $L = L_3^{\ep}(32)$ with $M = \Nor_A(L_3^{\ep}(2))$ then
  $\chi(1) > 2^{32}$. Take $\ep = +$ first. The elements of prime
  order fixing a point of $A/M$ have order $2$, $3$, $5$ or $7$ (and
  there at most two classes in $A \cap M$ of any of these). The
  largest possible centralizer in $G$ is that of an (inner)
  involution, of order $2^{15 (+1)}.31 < 2^{21}$ --- so
  $\chi(g)/\chi(1) < 2^{-10}$ for any $g \in A, g \ne 1$. If $\ep =
  -$, the calculation is similar: the relevant prime orders are $2$,
  $3$ and $5$, the largest centralizer in $G$ has order
  $2^{15}\cdot11\cdot|A/L|$, and we deduce that $\chi(g)/\chi(1) <
  2^{-9}$ for any $g \in L, g \ne 1$.

  (e) If $G=G_2(32)$ with $G \cap M = G_2(2)$, then $\chi(1) >
  2^{56}$. The elements of prime order in $G_2(2) \times 5$ have order
  $2$, $3$ (two classes), $5$ or $7$ (one class). The largest
  centralizer in $G$ is that of a central involution, of order
  $2^{30}\cdot31\cdot33$, so $< 2^{40}$.  The assertion follows as
  before.

  (f) If $L = \Sp_4(q)$ with $L \cap M = \Sp_4(q^{1/3})$, with $q=8$
  or $q=64$, we can use the explicit description of conjugacy classes
  in \cite{Enomoto}. If $q=8$ then $\chi(1) > 2^{20}$. The elements of
  prime order in $L \cap M = \Sp_4(2)$ have orders $2$ (three
  classes), $3$ (two classes) and $5$ (one class). The involutions fix
  $2^8\cdot21$ points in two of the classes and $2^8$ in the third.
  The elements of order $3$ fix $252$ points, the elements of order
  $5$ fix $13$ points.  The assertion follows. If $q=64$, then
  $\chi(1) > 2^{40}$ and exactly similar consideration shows that
  $\chi(g) \le 2^{16}\cdot117$ for all $1\ne g \in L$.  We need only
  consider outer involutions and field automorphisms.  If $q=64$ and
  $x$ is an involutory field automorphisms, then $x$ fixes
  $|\Sp_4(8)|/|\Sp_4(2)| < 2^{21}$ points.  If $x$ is a field
  automorphism of order $3$, then $x$ has $1+ 40 + 40$ fixed points
  for $q=8$.  If $q=64$, the contradiction is easier to observe.

  (g) If $L=\Omega _8^+(32)$ with $L \cap M = \Omega _8^+(2)$, then
  $n > 2^{112}$. The largest centralizer in $G$ is the parabolic
  subgroup corresponding to the central node, so all centralizers have
  orders less than $2^{90}$ and hence there is no problem.

  (h) Similarly, if $L={^3}D_4(32)$ then again $n > 2^{111}$
  whereas all the centralizers of elements in $L$ have orders in $G$
  less than $2^{92}$. The outer elements in $A$   have
  centralizers of order less than $2^{72}$.

  (i) Finally consider $L=U_3(32)$ with $M \cap U_3(32)$ the stabilizer
  of an orthogonal decomposition of the natural module. If exceptionality
  holds, then the proof shows that $B/G$ has order $15$.  Thus,
  $G=L$.
  
  Suppose that $g \in M \cap G$ has prime order.  If $g$ has order
  $11$, then $|g^G \cap M| \le 6$ and the centralizer is not larger
  than the stabilizer of a nonsingular $1$-space, whence $|g^G \cap
  M|/|g^G| < 1/85$.  If $g$ has order $3$, then $g$ is a regular
  semisimple element and so its centralizer has order at most $33^2$
  and $|x^G \cap M| \le 244$, whence the fixed point ratio is less
  than $1/85$. The only other class of prime elements in $G \cap M$
  are transvections.  Since the Sylow $2$-subgroup of $M \cap G$ has
  order $2$, the number of fixed points of $g$ is
  $|\Cen_G(g):\Cen_M(g)|=32^3 < 85|A:M|$.  This completes the
  proof.\end{proof}

\subsection{Affine action}\label{AffAction}

Recall that a primitive permutation group $A$ of degree $n$ is of
affine type if $A$ has an elementary abelian normal subgroup
$N=C_p^e$, which then is automatically regular. Thus $A=NM$ is a
semidirect product where $M$ is a point stabilizer and $M$ acts
irreducibly and faithfully on $N$ by conjugation. If $G$ is a
nontrivial (and hence transitive) normal subgroup of $A$, then $N\le
G$, since $N$ is the unique minimal normal subgroup of $A$.

The aim of this section is the following

\begin{Theorem}\label{Affine} Let $A$ be a primitive permutation
  group of degree $n=p^e$ of affine type, with minimal normal subgroup
  $N\cong C_p^e$. Let $G$ be a nontrivial normal subgroup of $A$ which
  admits a genus $0$ system $\cT$. Let $g_{\text{regular}}$ be the genus of
  $(G,\cT)$ with respect to the regular action of $G$.  Suppose there
  is a subgroup $B$ of $A$ containing $G$ such that $B/G$ is cyclic,
  and $(B,G)$ is exceptional. Then the following holds:\begin{itemize}
\item[(a)] If $g_{\text{regular}}=0$, then either\begin{itemize}
  \item[(i)] $n=p\ge3$, $G$ is cyclic, $\abs{\cT}=(p,p)$, or
  \item[(ii)] $n=p\ge5$, $G$ is dihedral of order $2p$,
    $\abs{\cT}=(2,2,p)$, or
  \item[(iii)] $n=4$, $G=C_2\times C_2$, $\abs{\cT}=(2,2,2)$, $B=A_4$,
    $A=A_4$ or $S_4$.\end{itemize}
\item[(b)] If $g_{\text{regular}}=1$, then $n=p$ or $n=p^2$ for $p$ odd, and
  one of the following holds.\begin{itemize}
\item[(i)] $\abs{\cT}=(2,2,2,2)$, $G=N\rtimes C_2$, and $n\ge5$, or
\item[(ii)] $\abs{\cT}=(2,3,6)$, $G=N\rtimes C_6$, and
  $n\equiv1\pmod{6}$, or
\item[(iii)] $\abs{\cT}=(3,3,3)$, $G=N\rtimes C_3$, and
  $n\equiv1\pmod{6}$, or
\item[(iv)] $\abs{\cT}=(2,4,4)$, $G=N\rtimes C_4$, and
  $n\equiv1\pmod{4}$.
\end{itemize}
\item[(c)] If $g_{\text{regular}}>1$, then also $G$ is primitive, and one of
  the following holds.\begin{itemize}
  \item[(i)] $n=11^2$, $\abs{\cT}=(2,3,8)$, $G=N\rtimes\GL_2(3)$ and
    $B=A=N\rtimes(\GL_2(3)\times C_5)$, or
  \item[(ii)] $n=5^2$, $\abs{\cT}=(2,3,10)$, $G=N\rtimes S_3$ and
    $A=B=N\rtimes(S_3\times C_4)$, or
  \item[(iii)] $n=5^2$, $\abs{\cT}=(2,2,2,4)$, $G/N$ is a Sylow
    $2$-subgroup of the subgroup of index $2$ in $\GL_2(5)$, $B/N$ is
    the normalizer of $G/N$ with $B/G$ of order $3$, and either $A=B$
    or $\abs{A/B}=2$, or
  \item[(iv)] $n=5^2$, $\abs{\cT}=(2,2,2,3)$, $G=N\rtimes D_{12}$ and
    $A/G$ is cyclic of order $2$ or $4$, and $B\le A$ properly
    contains $G$, or
  \item[(v)] $n=3^2$, $\abs{\cT}=(2,4,6)$, $(2,2,2,4)$, $(2,2,2,6)$,
    or $(2,2,2,2,2)$, $G=N\rtimes D_{8}$ and $\abs{A/G}=2$, or
  \item[(vi)] $n=2^4$, $\abs{\cT}=(2,4,5)$ or $(2,2,2,4)$, $G=N\rtimes
    D_{10}$, $B=N\rtimes(D_{10}\times C_3)$, and $A=B$ or
    $A/G=S_3$.\end{itemize}
\end{itemize}\end{Theorem}

The proof will be considerably long and difficult. For better
readability we divide the proof into smaller pieces. We keep the
hypothesis, in particular let $\cT=(x_1,x_2,\ldots,x_r)$ be a genus
$0$ generating system of $G$ with $x_i\ne1$ for all $i$. Let $d_i$ be
the order of $x_i$.

\subsubsection{$G$ is abelian} Since $G$ is abelian, the permutation
representation is the regular representation, so $G=N$. Thus each
element $x_i$ has order $p$, and the genus $0$ relation gives
$rp^e(1-1/p)=2(p^e-1)$, so either $r=2$ and $e=1$, or $r=3$ and
$p=2=e$, and (a) follows. Note that if $r=2$, then $p\ne2$ because
then $A=G$.

\subsubsection{$G''=1$} Now assume that $G''=1$, but $G$ is not
abelian. Since $A$ is irreducible on $N$, $N$ is a semisimple
$G$-module and $A$ permutes transitively the homogeneous
$G$-components of $N$. Let $Q$ be any minimal $G$-invariant subgroup
of $N$. So $N=Q\oplus P$ for some $G$-invariant subgroup $P$ (by
complete reducibility). By passing to $G/P$ and using \cite[3.8]{GT},
it follows that $G/\Cen_G(Q)$ is cyclic of order $1$, $2$, $3$, $4$ or
$6$.  Moreover, $\Cen_Q(x_i)=1$ for each $x_i$ which is nontrivial in
$G/P$.  Since $Q$ was arbitrary, it follows that since $N=\Cen_G(N)$
(because $N=\gFit(A)$), $G/N$ is abelian and $\abs{G/N}$ is not
divisible by any prime larger than $3$.

If $N=Q$, then $G$ acts irreducibly on $N$. Thus, $G/N$ cyclic (as
$G/N$ is abelian) and each element outside of $N$ has a unique fixed
point. It follows that $\cT$ is a Euclidean or sub--Euclidean system and
$G/N$ is cyclic of order $2$, $3$, $4$, or $6$.  Indeed, the same
argument applies as long as every element outside $N$ has a unique
fixed point.

Now suppose $N\ne Q$. Then we may identify $\Omega$ with $P\times Q$
with $G$ acting on $P$ and $Q$ (by viewing $G$ inside $\AGL(Q)\times
\AGL(P)$). Thus, Lemma \ref{P}(b) applies and so either each $x_i$ has
exactly $1$ fixed point, and as above $G/N$ is cyclic of order $2$,
$3$, $4$, or $6$ (because $G/N$ is faithful on $Q$) or we may assume
that $x_1$ has no fixed points.

We may assume that $x_1$ has no fixed points on $Q$ (as a $G$-set).
Since the action on $Q$ is of genus $0$, the argument yields that the
genus zero system for $G$ acting on $Q$ is sub--Euclidean and so must be
dihedral. In particular, $n$ is odd and $\abs{Q}=p$. By Lemma \ref{P},
it suffices to show that there are no such genus $0$ systems with
$n=p^2$. Then each $x_i$ has index $(1/2)(n-1)$, $(1/2)(n-p)$, $n-p$
or $(n-(p+1)/2)$ (in the corresponding cases that $x$ inverts $N$, $x$
is a reflection on $N$, $x$ is a translation or a commuting product of
a reflection and a translation). The only possible genus $0$ system is
$(2,2,p)$ or $(2,2,2,p)$ with $p=3$. In the first case, we get a
dihedral group and so $n\ne p^2$. In the latter case, we repeat the
argument to see that $n$ cannot be bigger than $9$. Then $A/G$ has
order $2$ and is not exceptional.

So we have shown that $G/N$ is cyclic of order $2$, $3$, $4$, or $6$
and $\cT$ is Euclidean or sub--Euclidean. Now use Lemma \ref{sph_euc_cor}
to see that in the Euclidean case, we get the configurations listed in
Theorem \ref{Affine}(b). If $\cT$ is sub--Euclidean, then either $G$ is
dihedral of order $2p$ and $\abs{\cT}=(2,2,p)$ --- the case listed in
(a)(ii) --- or $G\cong A_4$ by Proposition \ref{sph_euc} (the other
possibilities do not fulfill $G''=1$). But the pair $(S_4,A_4)$ is not
exceptional for the degree $4$ action.

\subsubsection{$G''\ne1$, $p>2$} Here we use the results in
\cite{GT} and \cite{Neubauer:CA1} to determine all possibilities. Note
that those results assume that $G$ is primitive. Using
\cite{Neubauer:CA1} and some fairly straightforward computations, we
will see that in fact $G$ must be primitive.

Let $M$ be the stabilizer of a point and set $H=G\cap M$. Then $A=NM$
(a semidirect product) and $G=NH$ (semidirect). Then $G$ is primitive
if and only if $G$ acts irreducibly on $N$. Since $A$ is irreducible
on $N$, it follows that $N=N_1\oplus\dots\oplus N_s$ where $N_i$ is
$G$-irreducible. Moreover, $A$ permutes transitively the isomorphism
classes of the $N_i$ (otherwise, $A$ would leave invariant a proper
subspace). In particular, each $N_i$ has the same cardinality.
Moreover, $G/C_G(N_i)$ is nonabelian for each $i$ (for otherwise, $G'$
acts trivially on some $N_i$ whence on each $N_i$ by transitivity --
then $G''=1$ contrary to hypothesis).

As a $G$-set, we may identify $\Omega$ with $N$ and so with
$N_1\times\dots N_s$. Thus, we may apply Lemma \ref{P}. In particular,
the action on each $N_i$ is genus $0$. We may now apply \cite[Theorem
1.5]{Neubauer:CA1} and Lemma \ref{P} to conclude that either $s=1$ or
the genus is positive unless possibly $\abs{N_1}=9$. If $p>11$, it
follows from \cite{GT} that $G'$ centralizes $N_1$, a contradiction as
above.  So $p\le11$.

Let us first consider the case that $\abs{N_1}=9$. If $N=N_1$, then
either $G$ is 2-transitive (and so exceptionality cannot hold) or
$G/N=D_8$ and $A/G$ has order $2$.  Then we may assume that $x_1$ and
$x_2$ are in $G$ but are not in $G'$. The elements in $G$ outside $G'$
are either involutions with $3$ fixed points or elements of order $6$
with no fixed points. Since $x_1x_2$ has order $4$, we see that if
$r=3$, we may assume that $x_1$ has order $2$, $x_2$ has order $6$ and
$x_3$ has order $4$. If $r\ge4$, then we see that we get the examples
$(2,2,2,4)$, $(2,2,2,6)$, and $(2,2,2,2,2)$ in the Theorem.

Now assume that $N\ne N_1$. We claim that in this case the genus is
positive. Write $x_i=(y_i,z_i)$ where $y_i$ is the permutation induced
by $x_i$ on $N_1$ and $z_i$ is the permutation induced on
$N_2\oplus\dots\oplus N_s$. If say, $y_1=1$, then
$\chi(y_1)\ind(z_1)\ge9(n/27)=n/3$ and Lemma \ref{P} implies that the
genus is positive. If $\chi(x_i)>0$ for each $i$, Lemma \ref{P} applies
as well. So we may assume that $\chi(y_i)=0$ for some $i$. This implies
that either $y_i$ has order $3$ and $\ind(y_i)=6$ or $y_i$ has order
$6$ and $\ind(y_i)=7$. If $r\ge4$, we see that $\sum \chi(y_i)\ge6$.
Since $\ind(z_i)\ge n/27$, we may apply Lemma \ref{P} to obtain a
contradiction. So $r=3$. Let $d_i'$ be the order of $y_i$. Since $\sum
(d_i'-1)/d_i'>2$ (or $G''=1$) and $\sum \chi(y_i)<6$ (or Lemma \ref{P}
applies), the only possibilities for $(d_1',d_2',d_3')$ are
$(2,4,6)$ or $(2,8,3)$ where $y_1$ is a noncentral involution in
$G/N$, $\chi(y_2)=1$ and $\chi(y_3)=0$.

Moreover, in the first case, $y_3^3$ is also a noncentral involution
in $G/N$. It follows that $y_3^2$ is a translation on $N_1$. This
implies that $G/C_G(N_1)$ is a $(2,4,2)$ group, i.e.\ it must be $D_8$
and contains no elements of order 8. Since $G/C_G(N_i)$ are all
isomorphic it follows that on each component we have a $(2,4,6)$ group
or a $(2,8,3)$ group.

By Lemma \ref{P}, we see that $3\ind(z_1)+\ind(z_2)<2n/9$. It follows
that $z_1$ must be an involution (otherwise, $\ind(z_1)\ge2n/27$,
contrary to the above inequality or $z_1$ would induce a transvection
on some component with fixed points and interchanging components we
would not have a triple as above). Moreover, $z_1$ must induce a
noncentral involution on each irreducible factor (or in some factor,
we would not have a triple as above). It follows that $z_2$ must
induce an element of order $3$, $4$, $6$ or $8$ on each factor. The
inequality above shows that if $s\ge3$, then the genus will be positive.
So we may assume that $s=2$. One computes directly in these cases that
$\sum\ind(x_i)>160$ and so the genus is positive.

So assume that $\abs{N_1}\ne9$. As we noted above, this implies that
$N=N_1$, that is, $G$ is primitive. We can now apply
\cite[1.5]{Neubauer:CA1} and consider the various possibilities. The
possibilities for the $d_i$ are given there. It is just a matter of
deciding whether exceptionality can hold.

If $p=11$, then $G/N=\GL_2(3)$. Since $A$ normalizes $G$, the only
possibility is that $A/G$ has order $5$.

If $p=7$, then $G/N=C_3\rtimes D_8$. Since $A$ normalizes $G$, the only
possibility is that $A/G$ has order $3$. However such a pair $(A,G)$
is not exceptional.

If $p=5$ and $n=125$, then we see that a point stabilizer has a unique
orbit of maximal size, whence exceptionality cannot hold.

Next consider $n=25$. There are $6$ possibilities for $G$ given in
\cite[1.5]{Neubauer:CA1}. Cases 4b, 4c and 4e are all $2$-transitive and
so cannot be exceptional. In case 4d, $G/N$ is self-normalizing in
$\GL_2(5)$ and so $A=G$ is not exceptional. In the remaining cases,
$A$ is contained in the normalizer of $G$ and it is elementary to
verify the theorem.

Finally, assume that $p=3$. The case $n=9$ has already been handled.
If $n=27$, then either $G$ is $2$-transitive or $G/N=S_4$. The only
possibility for $A$ (since $G$ is normal in $A$) would be
$A/N=S_4\times C_2$. It is straightforward to verify that this is
not exceptional.

The remaining case is $n=81$. If $G/N=S_5$, then $A/N=S_5\times C_2$
and we see that exceptionality fails. If $G/N=S_5\times C_2$, then
$G/N$ is self-normalizing in $\GL_4(3)$ and so $A=G$ is not
exceptional. The only other possibility is that $G/N=C_2\wr S_4$ is
the monomial subgroup of $\GL_4(3)$. Then the orbits of a point
stabilizer have distinct sizes and so exceptionality cannot hold.

This completes this subcase.

\subsubsection{$G''\ne1$, $p=2$} In this section we will show that
the only possibility is $e=4$, leading to case (c)(vi). We will us the
following result from \cite{GurNeu:Affine}, which extends previous
work of Guralnick, Neubauer and Thompson on genus $0$ systems in
affine permutation groups.

\begin{Proposition}
  Suppose that $A$ is an affine primitive permutation group of degree
  $2^e$ and $G$ is a nontrivial normal subgroup which admits a genus
  $0$ system. If $G''\ne1$, then $e\le6$.
\end{Proposition}

So we have to go through the possibilities $2\le e\le6$.

If $e\le2$, then $G''=1$ unless $G=S_4$ which is $2$-transitive, hence
exceptionality cannot hold.

If $e=3$ or $5$, then $G$ is primitive and it is easy to see that this
forces $G$ $2$--transitive.

Let $H$ be a point stabilizer in $G$.

Now suppose that $e=4$. Thus, $H$ is a subgroup of $A_8=L_4(2)$.
If $G$ has order divisible by $7$, then either $G$ is $2$-transitive
or $H$ has orbits of size $1$, $7$ and $8$. In either case,
exceptionality cannot hold. If $G$ has order divisible by $5$, then
the orbits of $H$ must be of size $1$, $5$, $5$, $5$ (by
exceptionality). If $G$ is not solvable, then this forces $H$ to
contain $A_5$. There are two conjugacy classes of $A_5$ in
$A_8$. One $A_5$ is transitive on the nonzero vectors of $N$
(namely $\SL_2(4)$) and the other has orbits of size $5$ and $10$.
Thus, for any $G$ containing $A_5$, exceptionality cannot hold.

Now suppose that $G$ is solvable. This implies that $H$ has either a
normal $3$ or $5$-subgroup. In the latter case, either $G$ is
$2$-transitive or $H$ has order $10$ or $20$. Since the point
stabilizer in $B$ must permute the three orbits of size $5$, it
follows that we have $\abs{B/G}=3$. Also, $B$ is then $2$--transitive,
and we easily see $\AGL_1(16)\le B$. From that we derive case (c)(vi).

So finally, we may assume that $F(H)=\O3(H)$. If $\O3(H)$ has order
$9$, then $\O3(H)$ has orbits of size $1$, $3$, $3$ and $9$ and so
exceptionality cannot hold. If $F(H)$ has order $3$, then since
$C_H(F(H))\subseteq F(H)$, it follows that $H$ has order $3$ or $6$.
Since $G''\ne 1$, the latter must occur. So $H\cong S_3$ and $H$ is
not irreducible on $N$. Thus, $N$ is the direct sum of $2$ copies of
the irreducible $2$-dimensional $H$-module. Then the $H$--orbits on
$N$ have size $1$, $3$, $3$, $3$ and $6$ and so exceptionality cannot
hold.

Finally, we have to look at $n=2^6$. We show that there are no
examples. First, if $A$ is solvable, then we use the list of solvable
primitive groups of degree $2^6$ provided by GAP \cite{GAP} to make an
exhaustive search. (The program used for this exhaustive search was
also applied to the list of the not necessarily solvable primitive
groups of degree $p^e<50$, and independently confirmed the list in
Theorem \ref{Affine}(c) for this range.)

So from now on assume that $A$ is not solvable. The nonsolvable
irreducible subgroups of $\GL_6(2)=L_6(2)$ are given in \cite{HarYam}.
Since $M$ has index $2^6$ in $A$, we are only interested in those
possibilities with a quotient of order not a power of $2$, see Lemma
\ref{cong1}. This restricts attention to the $A$ with $SL_2(8)\le A
\le\gL_2(8)$, $SL_3(4)\le A\le\gL_3(4)$, and $C_3\times L_3(2)\le A\le
L_2(2)\times L_3(2)$. In the first two of these cases, the action of
the subgroup $SL$ is already transitive on vectors, so the action of
$G$ is 2-transitive. In the last, the group $L_2(2) \times L_3(2)$
acts on the vector space as on the tensor product $U \otimes W$, with
$U$, $W$ of dimension $2$, $3$, respectively.  There are two orbits on
the non-zero vectors: If $u_1$ and $u_2$ are independent in $U$, the
orbit $\Omega_i$ of $u_1\otimes w_1+u_2\otimes w_2$ is determined by
the dimension $i$ of the subspace generated by $w_1, w_2$ in $W$. It
is clear that $\Omega_2$ is an orbit also under $L_3(2)$, so the
action is not exceptional.

\section{Dickson polynomials and R\'edei functions}\label{g=0}

In this section we study the rational functions of ramification type
$(2,2,n)$ and $(n,n)$, appearing in Theorem \ref{Affine}(a).

\begin{Definition} Let $0\ne\alpha\in\bK$ with $\alpha^2\in K$. Set
$\lambda(X)=\frac{X-\alpha}{X+\alpha}$. For $n\in\NNN$ set
$R_n(\alpha,X)=\lambda^{-1}(\lambda(X)^n)$. Then one immediately
verifies $R_n(\alpha,X)\in K(X)$. The function $R_n(\alpha,X)$ is
called \emph{R\'edei function}, see \cite{Nobauer}.\end{Definition}

Recall that the Dickson polynomial of degree $n$ belonging to $a\in K$
is defined implicitly by $D_n(a,Z+a/Z)=Z^n+(a/Z)^n$. Also, we remind
that two rational functions $f$ and $g$ over a field $K$ are said to
be equivalent, if they differ only by composition with linear
fractional functions over $K$.

\begin{Theorem}\label{Redei} Let $K$ be a field of characteristic $0$, 
  and $f\in K(X)$ be a rational function of degree $n$. Let
  $\abs{\cT}$ be the ramification type of $f$.
  \begin{itemize}
  \item[(a)] If $\abs{\cT}=(n,n)$, then $f$ is equivalent to $X^n$ or
    to $R_n(\alpha,X)$ for some $\alpha\in\bK$, $\alpha^2\in K$.
  \item[(b)] If $\abs{\cT}=(2,2,n)$, then $f$ is equivalent to
    $D_n(a,X)$ for some non--zero $a\in K$.
  \end{itemize}
\end{Theorem}

\begin{proof} In the following the phrase `we may assume' means that
  we replace $f$ by a suitable equivalent rational function. Let $n$
  be the degree of $f$. Part (b) is well known, see
  e.~g.~\cite{Fried:Schur}, \cite{Turnwald:Schur}. Thus suppose that
  the ramification type is $(n,n)$. First suppose that one of the two
  ramified places is rational. So we may assume that this place is at
  infinity. Then $f^{-1}(\infty)$ must also be rational, again we may
  assume that $f^{-1}(\infty)=\infty$. But this means that $f$ is a
  polynomial.
  
  Now suppose that none of the two ramified places are rational. As
  the absolute Galois group $\Gal(\bK|K)$ permutes these places, they
  are algebraically conjugate. Thus we may assume that they are given
  by $t\mapsto\alpha$ and $t\mapsto-\alpha$ for some
  $\alpha\in\bK\setminus K$ with $\alpha^2\in K$. Then the two
  elements $f^{-1}(\pm\alpha)$ are in $K[\alpha]$ and are interchanged
  by the involution in $\Gal(K(\alpha)|K)$. Thus we may assume without
  loss that $f^{-1}(\alpha)=\alpha$ and $f^{-1}(-\alpha)=-\alpha$. The
  maps $X\mapsto\frac{aX+\alpha^2}{X+a}$ for $a\in K$ fix $-\alpha$
  and $\alpha$, but are transitive on $\PP(K)$. So we may assume
  $f(\infty)=\infty$. Thus there is a polynomial $R(X)\in\bK[X]$ of
  degree less than $n$ such that
\begin{align*}
f(X)-\alpha &= \frac{(X-\alpha)^n}{R(X)}\\
f(X)+\alpha &= \frac{(X+\alpha)^n}{R(X)}.
\end{align*}
Now eliminate $R(X)$ from these two equations to obtain the
claim.\end{proof}

\begin{Proposition}\label{Redex} Let $n$ be a prime, and $R_n(\alpha,X)$
  be the R\'edei function as above. Then
  $K(\alpha\frac{1+\zeta}{1-\zeta})$ is the algebraic closure of $K$
  in a splitting field of $f(X)-t$, where $\zeta$ is a primitive
  $n$--th root of unity. In particular, $R_n(\alpha,X)$ is
  arithmetically exceptional over a number field $K$ if and only if
  $\alpha\frac{1+\zeta}{1-\zeta}\not\in K$.
\end{Proposition}

\begin{proof} Let $x$ be a root of $f(X)-t$, and $\zeta$ be a
  primitive $n$--th root of unity. It is immediate to verify that the
  other roots of $f(X)-t$ are given by
  $x_i=\omega_i\frac{x+\alpha^2/\omega_i}{x+\omega_i}$, where
  $\omega_i=\alpha\frac{1+\zeta^i}{1-\zeta^i}$. This shows that the
  algebraic closure of $K$ in said splitting field is generated by the
  $\omega_i$. But if $\sigma\in\Gal(\bK|K)$, then either
  $\omega_1^\sigma=\omega_i$ or $\omega_1=\alpha^2\omega_i$ for some
  $i$, so the field generated by $\omega_i$, $i=1,2,\dots,n-1$, is
  independent of $i$.\end{proof}

Similarly as above, one proves

\begin{Proposition} Let $n$ be a prime, and $\zeta$ be a primitive
  $n$--th root of unity. The algebraic closure of $K$ in a splitting
  field of $X^n-t$ or $D_n(a,X)-t$ ($0\ne a\in K$) is $K(\zeta)$ or
  $K(\zeta+1/\zeta)$, respectively.\end{Proposition}

So when $K$ is fixed, the field of constants of $D_n(a,X)$ depends
only on the degree $n$, in contrast to the case for the R\'edei
functions. This has the following surprising consequence.

\begin{Corollary}\label{RedeiCor} There are $3$ arithmetically
exceptional R\'edei functions of degree $3$ over $\QQQ$, such that
their composition is not arithmetically exceptional.\end{Corollary}

\begin{proof} Let $\zeta$ be a primitive third root of unity. We have
  $\frac{1+\zeta}{1-\zeta}=\sqrt{-1/3}$. For $m\in\ZZZ$ set
  $f_m(X)=R_3(\sqrt{-3m},X)$. Then, by the previous result, the field
  of constants of $f_m(X)$ is $\QQQ(\sqrt{m})$. Now, by the
  multiplicativity of the Legendre symbol, every rational prime splits
  in at least one of the fields $\QQQ(\sqrt{-1})$, $\QQQ(\sqrt{-2})$,
  and $\QQQ(\sqrt{2})$. So the composition of $f_{-1}$, $f_{-2}$, and
  $f_{2}$ is not arithmetically exceptional by Theorem
  \ref{togroups}(b).\end{proof}

\begin{Remark*} A composition of arithmetically exceptional
  polynomials over $\QQQ$ is always arithmetically exceptional. This,
  however, is not true anymore over number fields. R.~Matthews
  \cite[6.5]{Matthews} showed that if $n$ is the product of three
  distinct odd primes, then there is a subfield $K$ of $\QQQ(\zeta)$
  ($\zeta$ a primitive $n$--th root of unity) such that $X^n$ is not
  arithmetically exceptional over $K$, whereas $X^d$ is for each
  proper divisor $d$ of $n$.

The permutation behavior of the R\'edei functions on finite fields
has been investigated in \cite{Nobauer}.\end{Remark*}

\section{Rational functions with Euclidean ramification
  type}\label{g=1}

Let $K$ be a field of characteristic $0$, and $f\in K(X)$ be an
indecomposable rational function. This section is devoted to the case
that the Galois closure of $K(X)|K(f(X))$ has genus $1$.

From Lemma \ref{sph_euc_cor} we know that such an $f$ has degree $p$
or $p^2$ for an odd prime $p$. Furthermore, the possible ramification
types are $(2,2,2,2)$, $(2,3,6)$, $(3,3,3)$, and $(2,4,4)$. Recall
that the associated geometric monodromy groups are isomorphic to
$N\rtimes C_2$, $N\rtimes C_6$, $N\rtimes C_3$, and $N\rtimes C_4$,
respectively, where $N$ is elementary abelian of order $\deg f$. We
will use these facts throughout this section.

\subsection{Elliptic Curves}
In the following we summarize some of the most basic notions and
results about elliptic curves. There is a vast literature on this
subject. Our main reference are the two excellent volumes
\cite{Silverman1} and \cite{Silverman2} by Silverman.

Let $K$ be a field of characteristic $0$. A smooth projective curve
$E$ of genus $1$ over $K$ together with a $K$-rational point $O_E$ is
called an \emph{elliptic curve}. The curve $E$ is a (commutative)
algebraic group with $O_E$ the neutral element.

There is an isomorphism, defined over $\QQQ$, to a plane curve in
Weierstra\ss\ normal form, which is the projectivization of
$Y^2=X^3+aX+b$, where $a,b\in K$, and the discriminant of the right
hand side does not vanish.

If we choose such a Weierstra\ss\ model, then a (not unique) structure
as an algebraic group can be defined by the requirement that the
points $-P$, $-Q$ and $P+Q$ lie on a line, where $-(x,y):=(x,-y)$.

If $E$ and $E'$ are elliptic curves, and $\Phi:E\to E'$ is a
non-constant rational map defined over $K$, then $\Phi$ is
automatically a surjective morphism, and the modified map
$P\mapsto\Phi(P)-\Phi(O_E)+O_{E'}$ is a homomorphism of algebraic
groups. Such maps are called \emph{isogenies}. If furthermore $E=E'$,
then we talk of \emph{endomorphisms}, and an \emph{automorphism} is a
bijective endomorphism.

The structure of the automorphism group $\Aut(E)$ of an elliptic curve
is very easy: Suppose that $E$ is given in Weierstra\ss\ form
$Y^2=X^3+aX+b$.  There is always the automorphism of order $2$,
$P\mapsto-P$, which, expressed in coordinates, is
$(x,y)\mapsto(x,-y)$. If $ab\ne0$ then there are no other
automorphisms. Next suppose that $b=0$. Setting $i:=\sqrt{-1}$, we
have the automorphism $(x,y)\mapsto(-x,iy)$ of order $4$. This
generates the whole automorphism group. Now suppose that $a=0$, and
let $\omega$ be a primitive $3$rd root. Then $(x,y)\mapsto(\omega
x,-y)$ has order $6$, and generates the automorphism group.

A very specific class of endomorphisms is the multiplication by $m$
map $P\mapsto mP$. This map is usually denoted by $[m]$, and its
kernel by $E[m]$. The kernel, as a group, is isomorphic to $C_m\times
C_m$. More generally, if $\Phi:E\to E'$ is an isogeny, then we denote
the kernel by $E[\Phi]$. The order of the kernel is the degree of the
map $\Phi$.

If $K=\QQQ$, then the only endomorphisms of $E$ which are defined over
$\QQQ$ are the multiplication maps $[m]$. In general, there are two
possibilities: The only endomorphism are the maps $[m]$, or $K$
contains an imaginary quadratic field $k$, and the endomorphism ring
is isomorphic to an order of $k$. In this case one says that $E$ has
\emph{complex multiplication}. On this plane curve level, it is not
possible to understand these endomorphisms. Instead, one has to choose
an analytic isomorphism between $E(\CCC)$ and $\CCC/\Lambda$, where
$\Lambda$ is a lattice. Then the endomorphisms of $E$ are induced by
multiplication with $0\ne\alpha\in\CCC$ with
$\alpha\Lambda\subseteq\Lambda$. Without loss assume that $\Lambda$ is
generated by $1$ and $\omega\in\CCC\setminus\RRR$. If there is an
$\alpha\in\CCC\setminus\ZZZ$ with $\alpha\Lambda\subseteq\Lambda$,
then it follows easily that $\omega$ is imaginary quadratic, and the
set of $\alpha$'s with $\alpha\Lambda\subseteq\Lambda$ is an order in
the ring of integers of $k=\QQQ(\omega)$. Later we will be concerned
only with the case that the endomorphism ring is the maximal order
$\ok$ of $k$. If this is the case, then $\Lambda$ can be chosen to be
a fractional ideal of $k$. If $\alpha\in\ok$ corresponds to an
endomorphism $\Phi$ of $E$, then the algebraic conjugate of $\alpha$
corresponds to the dual endomorphism. In particular, the degree of
$\Phi$ is the norm $N_{k|\QQQ}(\alpha)$.

We have described the structure of the automorphism group of $E$. In
our development, we will obtain isomorphisms (as algebraic curves) of
$E$ to itself (defined over $\bK$) which not necessarily fix the
origin. Denote this group of isomorphisms by $\Isom(E)$. If
$\bar\beta\in\Isom(E)$, then $\beta:P\mapsto\beta(P)-\beta(0)$ is in
$\Aut(E)$. Thus $\Isom(E)$, as a group acting on $E(\bK)$, carries the
natural structure as a semidirect product $E(\bK)\rtimes\Aut(E)$.

\begin{Lemma}\label{L:betaE} Let $\beta\in\Aut(E)$ have order
  $m>1$. Then $1+\beta+\beta^2+\dots+\beta^{m-1}=0$.
\end{Lemma}

\begin{proof} $\varphi:=1+\beta+\beta^2+\dots+\beta^{m-1}$ is an
  endomorphism of $E$, so either $\varphi=0$, or $\varphi$ is
  surjective on $E(\bK)$. The latter cannot happen, because each
  element in the image of $\varphi$ is fixed under $\beta$. Thus
  $\varphi=0$.
\end{proof}

From that we obtain

\begin{Corollary} Let $\bar\beta\in\Isom(E)$ with
  $\bar\beta(P)=\beta(P+w)$, where $1\ne\beta\in\Aut(E)$ and $w\in
  E(\bK)$. Then $\bar\beta$ and $\beta$ have the same order.
\end{Corollary}

\begin{Lemma}\label{L:EGal} Let $\psi:E\to C$ be a covering of finite
  degree to a curve $C$. Suppose that $\psi$ does not factor as $E\to
  E'\to C$ where $E\to E'$ is a covering of elliptic curves of degree
  $>1$. Then $\psi(P)=\psi(P+w)$ for all $P\in E(\bK)$ and a fixed
  $w\in E(\bK)$ implies $w=0$.
\end{Lemma}

\begin{proof} Suppose that $w\ne0$. As the fibers of $\psi$ are
  finite, we obtain that $w$ has finite order. Thus $E/\gen{w}$ is an
  elliptic curve $E'$, and $\psi$ obviously factors through $E\to
  E/\gen{w}$.
\end{proof}

\begin{Corollary}\label{C:EPP} Let $1\ne\beta\in\Aut(E)$ be defined
  over $K_\beta\supseteq K$, $w\in E(\bK)$, and set
  $\bar\beta(P):=\beta(P+w)$ for $P\in E$. Then the quotient curve
  $E/\gen{\bar\beta}$ has genus $0$. Suppose that $\bar\psi:E\to
  E/\gen{\bar\beta}$ is defined over $K_\beta$. Then $w\in
  E(K_\beta)$.
\end{Corollary}

\begin{proof} Clearly $\bar\beta$ has a fixed point on $E$. This point 
  is ramified in the covering $E\to E/\gen{\bar\beta}$. However, by
  the Riemann-Hurwitz genus formula, a covering of genus $1$ curves is
  unramified. Thus $E/\gen{\bar\beta}$ has genus $0$.
  
  Set $G_{K_\beta}:=\Gal(\bK|K_\beta)$. Apply $\gamma\in G_{K_\beta}$
  to $\bar\psi(P)=\bar\psi(\bar\beta(P))$ for $P$ a generic point on
  $E$ to obtain
\begin{align*}
\bar\psi(\beta(P)+w)               &= \bar\psi(\bar\beta(P))\\
                                   &= \bar\psi(\bar\beta^\gamma(P))\\
                                   &= \bar\psi(\beta(P+w^\gamma))\\
                                   &= \bar\psi(\beta(P)+w^\gamma),
\end{align*}
so $w^\gamma=w$ for all $\gamma$ by Lemma \ref{L:EGal}, hence $w\in
E(K_\beta)$. (Note that $E\to E/\gen{\bar\beta}$ is a cyclic cover,
and each proper intermediate curve corresponds to a non-trivial
subgroup of $\gen{\bar\beta}$, so has genus $0$ as well by the
argument from the beginning of the proof.)
\end{proof}

\subsection{Non--existence results}

In this section we associate rational functions with Euclidean
ramification type to isogenies or endomorphisms of elliptic curves.

\begin{Theorem}\label{T:2222main} Let $K$ be a field of characteristic 
    $0$, $f\in K(X)$ be an indecomposable rational function of degree
  $n$ with Euclidean ramification type $\abs{\cT}=(2,2,2,2)$. Then,
  upon linearly changing $f$ over $K$, the following holds:
  
  There exist elliptic curves $E$, $E'$ over $K$, a $K$-rational
  isogeny $\Phi:E\to E'$ of degree $n$, a point $w\in E(K)$, such that
  the following diagram commutes, where $\bar\psi:E\to\PP\isom
  E/\gen{\bar\beta}$ with $\bar\beta(P)=-(P+w)$ for $P\in E$, and
  likewise $\bar\psi:E'\to\PP\isom E'/\gen{\bar\beta'}$ with
  $\bar\beta'(P')=-(P'+\Phi(w))$ for $P'\in E'$:
\[
\begin{array}{ccc}
E & \stackrel{\Phi}\longrightarrow & E'\\
\Big\downarrow \bar\psi & & \Big\downarrow \bar\psi'\\
\PP & \stackrel{f}\longrightarrow & \PP
\end{array}
\]
\end{Theorem}

\begin{Remark*} In general we cannot assume $w=0$ by changing the
  neutral element of the group $E$, for in order to do so, there must
  be a $K$-rational point $y\in E$ with $2y=w$. Indeed, if
  $\psi:E\to\PP$ is the map corresponding to $w=0$, then
  $\bar\psi(P)=\lambda(\psi(P+y))$ for $P\in E$, $2y=w$, and a
  suitable linear fractional function $\lambda(X)\in\bK(X)$. Note that
  we cannot take $\lambda(X)=X$, because the map $P\mapsto\psi(P+y)$
  is not defined over $K$ unless $y$ is $K$-rational.
  
  Similar remarks apply in the following cases.
\end{Remark*}

\begin{Theorem}\label{T:CMmain} Let $K$ be a field of characteristic
  $0$, $f\in K(X)$ be an indecomposable rational function of degree
  $n$ with Euclidean ramification type $\abs{\cT}=(2,3,6)$, $(3,3,3)$,
  or $(2,4,4)$. Let $j\in\{3,4,6\}$ be the maximum of the entries of
  $\abs{\cT}$. Then, upon linearly changing $f$ over $K$, the
  following holds:
  
  There is an elliptic curve $E$ over $K$ with $\beta\in\Aut(E)$
  defined over $K_\beta\supseteq K$ of order $j$ and $w\in
  E(K_\beta)$, a $K$-rational covering $\bar\psi:E\to\PP\isom
  E/\gen{\bar\beta}$ where $\bar\beta(P)=\beta(P+w)$, and an
  endomorphism $\Phi$ of $E$ defined over $K$ of degree $n$ with
  $\Phi(w)=w$, such that the following diagram commutes:
\[
\begin{array}{ccc}
E & \stackrel{\Phi}\longrightarrow & E\\
\Big\downarrow \bar\psi & & \Big\downarrow \bar\psi\\
\PP & \stackrel{f}\longrightarrow & \PP
\end{array}
\]
More precisely, we may assume the following.
\begin{enumerate}
\item[(a)] If $\cT=(2,3,6)$, then $E$ is given by $Y^2=X^3-1$, $w=0_E$,
  and $\bar\psi((x,y))=x^3$.
\item[(b)] If $\cT=(3,3,3)$, then $E$ is given by
  $Y^2=X^3-(R^3+27/4S^2)$ ($R,S\in K$), $w=(R,3/2\sqrt{-3}S)$, and
  $\bar\psi((x,y))=(-9RS-9xS+2xy-2Ry)/(6R^2+6Rx+6x^2)$.
\item[(c)] If $\cT=(2,4,4)$, then $E$ is given by $Y^2=X^3+AX$,
  $w=(0,0)$, and $\bar\psi((x,y))=(x^2-A)/(4x)$.
\end{enumerate}
\end{Theorem}

A consequence of the theorems is

\begin{Theorem}\label{ellmain} Let $K$ be a field of characteristic
  $0$, $f\in K(X)$ be an indecomposable rational function of degree
  $n$ with Euclidean ramification type $\abs{\cT}=(2,2,2,2)$,
  $(2,3,6)$, $(3,3,3)$, or $(2,4,4)$. Then $n=p$ or $p^2$ for an odd
  prime $p$. If $n=p$, then the following holds.\begin{itemize}
  \item[(a)] If $\cT=(2,2,2,2)$ and $K=\QQQ$, then $n\in\{3,$ $
    5,$ $ 7,$ $ 11,$ $ 13,$ $ 17,$ $ 19,$ $ 37,$ $ 43,$ $ 67,$ $
    163\}$.
  \item[(b)] If $\cT=(2,3,6)$ or $(3,3,3)$, then
    $\QQQ(\sqrt{-3})\subseteq K$ and $n\equiv1\pmod{6}$.
  \item[(c)] If $\cT=(2,4,4)$, then $\QQQ(\sqrt{-1})\subseteq K$ and
    $n\equiv1\pmod{4}$.
\end{itemize}\end{Theorem}

\begin{proof} From Lemma \ref{sph_euc_cor} it follows that $n=p$ or
  $p^2$ for an odd prime $p$.
  
  Suppose that $\abs{\cT}=(2,2,2,2)$ and $n=p$. From Theorem
  \ref{T:2222main} we obtain a rational isogeny of elliptic curves of
  degree $p$, with everything defined over $\QQQ$. A deep result of
  Mazur \cite[Theorem 1]{Mazur} implies that $p$ is as in the theorem.
  
  Next suppose that $\cT=(2,3,6)$ and $n=p$. We obtain a $K$-rational
  endomorphisms of degree $p$ of an elliptic curves which has complex
  multiplication by the maximal order of $\QQQ(\sqrt{-3})$. As this
  endomorphism is not the multiplication map, we get $\sqrt{-3}\in K$.
  Also, $p\equiv1\pmod{6}$, because $p$ is a norm in the maximal order
  of $\QQQ(\sqrt{-3})$.
  
  Analogously argue for the remaining cases.
\end{proof}

In Section \ref{S:ECex} we show that such functions $f$ indeed exist
in those cases which are not excluded by the theorem. However, in
order to keep the discussion of the field of constants reasonably
simple, we will assume $w=0$.

\begin{proof}[Proof of Theorems \ref{T:2222main}, \ref{T:CMmain}]
  We need a simple observation.

\begin{Lemma}\label{L:qe} Let $K$ be a field of characteristic
  $\ne2$, and $q(X)=X^4+aX^3+bX^2+cX+d\in K[X]$ be separable. Let $C$
  be the projective completion of the curve $Y^2=q(X)$. Then $C$ is
  $K$--birationally equivalent to an elliptic curve over
  $K$.\end{Lemma}

\begin{proof} Without loss assume that $a=0$. A birational
  correspondence is given by
\begin{align*}
(U,V) &\mapsto (X,Y)=(\frac{V-c}{2(U+b)},\frac{U-2X^2}{2}),\\
(X,Y) &\mapsto (U,V)=(2(Y+X^2),(4(Y+X^2)+2b)X+c),
\end{align*}
where $U,V$ are the coordinate functions of an elliptic curve in
Weierstra\ss\ normal form.\end{proof}

\paragraph{The $(2,2,2,2)$-case.}
Suppose that $f$ has ramification type $(2,2,2,2)$. Let
$\lambda_i\in\PP(\bK)$, $i=1,2,3,4$, be the branch points of $f$.  A
linear fractional change over $K$ allows us to assume the following:
\begin{itemize}
\item[(a)] The branch points of $f$ are finite.
\item[(b)] $f(X)=Q(X)/R(X)$ with $Q,R\in K[X]$ monic and relatively
  prime, and $n:=\deg(f)=\deg(Q)=1+\deg(R)$.
\end{itemize}

Let $G_K$ be the absolute Galois group of $K$. As $G_K$
permutes the branch points among themselves, we have
\[
q_\lambda(X):=(X-\lambda_1)(X-\lambda_2)
             (X-\lambda_3)(X-\lambda_4)\in K[X].
\]
Let $\mu_i$ be the unique non--critical point in $f^{-1}(\lambda_i)$.
Again, $G_K$ permutes the $\mu_i$ among themselves, thus
\[
q_\mu(X):=(X-\mu_1)(X-\mu_2)(X-\mu_3)(X-\mu_4)\in K[X].
\]
The ramification data imply the existence of monic polynomials
$Q_i\in\bK[X]$ such that
\[
f(X)-\lambda_i=(X-\mu_i)\frac{Q_i(X)^2}{R(X)}.
\]
Multiply these equations for $i=1,2,3,4$ to obtain
\begin{equation}\label{dgl}q_\lambda(f(X))=q_\mu(X)(S(X))^2,\end{equation}
where $S(X)=Q_1(X)Q_2(X)Q_3(X)Q_4(X)/R^2(X)\in K(X)$. (Actually, one
can easily show that $S(X)=f'(X)$, but we do not need that here.)

Consider the projective curves $C_\lambda, C_\mu$, where the affine
parts are given by
\begin{align*}
C_\lambda:\;  Y^2 &= q_\lambda(X)\\
C_\mu:\; Y^2 &= q_\mu(X).
\end{align*}
From \eqref{dgl} we infer that the map
\[
\tilde\Phi:(X,Y)\mapsto(f(X),YS(X))
\]
induces a rational map from $C_\mu$ to $C_\lambda$.  By Lemma
\ref{L:qe}, there are elliptic curves $E_\lambda, E_\mu$ over $K$,
which are $K$--birationally equivalent to $C_\lambda$ and $C_\mu$,
respectively. So $\tilde\Phi$ induces an isogeny $\Phi$ from $E_\mu$
to $E_\lambda$. As $\tilde\Phi^{-1}(P)$ has $n=\deg f$ elements in
$C_\lambda(\bK)$ for all $P$ in an open subset of $C_\mu(\bK)$, we get
that the isogeny $\Phi$ has degree $n$. Compose the isomorphism of
$E_\mu$ with $C_\mu$ with the projection to the $X$-coordinates of
$C_\mu$.  This gives a degree $2$ cover $\bar\psi:E_\mu\mapsto\PP$
which is defined over $K$. Thus $E_\mu$ has an isomorphism $\bar\beta$
of order $2$ to itself with $\bar\psi(P)=\bar\psi(\bar\beta(P))$ for
all $P\in E_\mu(\bK)$. The claim follows from Corollary \ref{C:EPP}.

\paragraph{The $(2,3,6)$--case.}\label{236case} As the absolute
Galois group $G_K$ of $K$ permutes the branch points among themselves,
but also preserves the cycle type above them, the branch points need
to be $K$--rational.  Without loss assume that the branch points
belonging to the inertia generators of orders $6$, $3$, and $2$ are
$\infty$, $0$, and $1$, and that the non-critical point in the fiber
$f^{-1}(\lambda)$ of a branch point $\lambda$ (which has to be
$K$--rational) is $\lambda$. The ramification information yields the
existence of polynomials $R,Q_0,Q_1\in K[X]$ and a non-zero constant
$c\in K$ such that
\begin{align*} f(X) &= cX\frac{Q_0(X)^3}{R(X)^6}\\
\intertext{and}
  f(X)-1 &= c(X-1)\frac{Q_1(X)^2}{R(X)^6}.
\end{align*}
Note that $f(1)=1$, so substituting $1$ in the first equation shows
that $c$ is a $3$rd power in $K$. Similarly, $c$ is a square in $K$,
as we see from substituting $0$ in the second equation. Thus $c$ is a
$6$th power in $K$. Therefore we may assume that $c=1$ by suitably
changing $R$.

These two equations show that the mapping
\[
(U,V)\mapsto(U',V'):=(U\frac{Q_0(U^3)}{R(U^3)^2},V\frac{Q_1(U^3)}{R(U^3)^3})
\]
induces a rational map of the elliptic curve
\[
E:\;V^2=U^3-1.
\]
to itself, sending $0_E$ to $O_E$. Namely
\begin{align*}
{U'}^3-1 &= U^3\frac{Q_0(U^3)^3}{R(U^3)^6}-1\\
       &= f(U^3)-1\\
       &= (U^3-1)\frac{Q_1(U^3)^2}{R(U^3)^6}\\
       &= V^2\frac{Q_1(U^3)^2}{R(U^3)^6}\\
       &= {V'}^2.
\end{align*}
Comparing degrees we see that this rational map has degree
$n=\deg(f)$. So we get an endomorphism of degree $n$ of the elliptic
curve $E$ which is defined over $K$, and the claim follows.

\paragraph{The $(3,3,3)$--case.}\label{333case} Let $\lambda_i$,
$i=1,2,3$ be the branch points of $f$, and $\mu_i$ be the simple
element in the fiber $f^{-1}(\lambda)$. Let $L(X)\in\bK(X)$ be the
linear fractional function with $L(\lambda_i)=\mu_i$ for $i=1,2,3$. As
$G_K$ permutes the $\lambda_i$, and the action on the $\lambda_i$ is
compatible with the action on the $\mu_i$, we obtain that
$L^\gamma(\lambda_i)=\mu_i$ for all $i$ and $\gamma\in G_K$. As $L$ is
unique, we have $L=L^\gamma$, so $L(X)\in K(X)$. Thus, upon replacing
$f(X)$ with $f(L(X))$, we may assume that $\lambda_i=\mu_i$. Also, as
before, we may assume that the $\lambda_i$ are finite and
$f(\infty)=\infty$. Set $q(X)=\prod(X-\lambda_i)$.

The ramification information gives polynomial $Q_i(X)\in\bK[X]$,
$R(X)\in K[X]$ such that
\[
f(X)-\lambda_i=(X-\lambda_i)\frac{Q_i(X)^3}{R(X)}.
\]
Multiply for $i=1,2,3$ to get
\[
q(f(X))=q(X)S(X)^3,
\]
where $S(X)\in K(X)$. Thus
\[
(X,Y)\mapsto(f(X),YS(X))
\]
gives a degree $n$ map of the cubic curve $Y^3=q(X)$ to itself. This
non-singular cubic curve is isomorphic over $K$ to an elliptic curve
of the form $V^2=U^3-\ell$, because one of the three points at
infinity is $K$-rational, and the cubic has the degree $3$
automorphism $(X,Y)\mapsto (X,\omega Y)$. From that we obtain the
claim.

\paragraph{The $(2,4,4)$--case.}\label{244case} The action of
$G_K$ on the branch points shows that the point corresponding to the
inertia generator of order $2$ is $K$--rational. So a
linear--fractional change can move it to $\infty$. Additionally, we
assume that the non--critical point in $f^{-1}(\infty)$, which must be
$K$--rational, is $\infty$. An argument as in the previous case shows
that there is $\ell\in K$ and $S\in K(X)$ with
\[
f(X)^2-\ell=(X^2-\ell)S(X)^4.
\]
Let $C$ be the curve defined by $Y^4=X^2-\ell$. This curve is
$K$--birationally equivalent to the elliptic curve $V^2=U^3-4\ell U$
by Lemma \ref{L:qe}. (The correspondence is given by
$(U,V)=(2(X+Y^2),4Y(X+Y^2))$, $(X,Y)=(U^2+4\ell)/(4U),V/(2U)$.)
Furthermore, we have the rational map $(X,Y)\mapsto (f(X),YS(X))$ of
$C$ to itself, which induces a $K$--rational endomorphism of degree
$n$ of the elliptic curve. Conclude as in the previous
case.\end{proof}

\subsection{Existence results}\label{S:ECex}
The existence results are, as suggested by the arguments for the
non-existence results, based on isogenies or endomorphisms of elliptic 
curves.

\subsubsection{Rational functions from isogenies}\label{S:Isog}
We fix a field $K$ of characteristic $0$ and elliptic curves $E$ and
$E'$ defined over $K$. Let $\Phi:E\to E'$ be an isogeny which is
defined over $K$. In particular, $\Phi$ is a homomorphism of groups.

Let $N:=E[\Phi]<E(\bK)$ be the kernel of $\Phi$. Furthermore, let
$\beta$ denote a non-trivial automorphism of $E$, which is defined
over $K$, and fulfills $\beta(N)=N$. Note that $E'(\bK)\isom
E(\bK)/N$, so $\beta$ induces an automorphism $\beta'$ of $E'$ such
that the following diagram commutes:
\[
\begin{array}{ccccccccc}
0 & \longrightarrow & N & \longrightarrow & E(\bK) &
\stackrel{\Phi}\longrightarrow & E'(\bK) & \longrightarrow & 0\\
 & & \;\;\Big\downarrow\beta & & \Big\downarrow\beta &
& \Big\downarrow\beta' & & \\
0 & \longrightarrow & N & \longrightarrow & E(\bK) &
\stackrel{\Phi}\longrightarrow & E'(\bK) & \longrightarrow & 0
\end{array}
\]
Set $\tK:=K(N)$, the field generated by the finite coordinates of
$N$. Note that $N$ is invariant under the absolute Galois group of $K$ 
because $\Phi$ is defined over $K$, so $\tK$ is a normal extension of
$K$.

Let $T$ be the group of translations on $E$ by elements of $N$. For
$\eta\in N$ and $P\in E(\bK)$, let $t_\eta(P):=P+\eta$. The action of
$T$ is obviously defined over $\tK$. Let $K(E)$ denote the function
field of $E$ over $K$, and let $\tK(E)$ be the compositum of $\tK$ and
$K(E)$. As usual, if $\rho$ is an endomorphism (in the sense of
curves) of $E$ defined over $\tK$, then the comorphism $\rho^\star$ is
defined by $\rho^\star(f)(P):=f(\rho(P))$ for $f\in\tK(E)$ and $P\in
E(\bK)$. Similarly, $\Phi^\star$ is an injective homomorphism
$K(E')\to K(E)$ of fields. By replacing $E'$ with a $K$-isomorphic
copy, we may assume that $\Phi^\star$ is the inclusion $K(E')\subseteq
K(E)$.

We get $\Gal(\tK(E)|\tK(E'))=T^\star$.

Next set $H:=\Gal(\tK|K)$. Let $H$ act on $\tK(E)$ by fixing
elementwise $K(E)$, and acting naturally on $\tK$. To stay confirm
with the above $\star$-notation, we denote this group by $H^\star$.
Finally, let $\beta^\star$ be the comorphisms of $\beta$.

As $\beta$ is defined over $K$, we obtain that $\beta$ centralizes
$H$, so $\beta^\star$ centralizes $H^\star$.

Let $A\le\Aut(\tK(E))$ be the group generated by $T^\star$,
$H^\star$, and $\beta^\star$.

\begin{Lemma}\label{L:act} $T^\star$ is a normal subgroup of
  $A$. More precisely, if $\eta\in N$ and $h\in H$, then
\begin{align*}
{\beta^\star}t_\eta^\star{\beta^\star}^{-1} &= t_{\beta(\eta)}^\star,\\
{h^\star}t_\eta^\star{h^\star}^{-1} &= t_{h(\eta)}^\star.
\end{align*}
\end{Lemma}

\begin{proof} For $f\in\tK(E)$ and $P\in E(\tK)$ we compute
\begin{align*}
{\beta^\star}t_\eta^\star{\beta^\star}^{-1}(f)(P) &=
                    f(\beta(\beta^{-1}(P)+\eta)) \\
  &= f(P+\beta(\eta)) \\
  &= t_{\beta(\eta)}^\star(f)(P).
\end{align*}
The other case is not quite analogous, because $h^\star$ is not the
comorphism of an automorphism of $E$. Instead, we have the relation
\[
h^\star(f)(P)=h(f(h^{-1}(P))).
\]
Using this, we obtain
\begin{align*}
{h^\star}t_\eta^\star{h^\star}^{-1}(f)(P) &=
                    h(t^\star{h^\star}^{-1})(f)(h^{-1}(P)))\\
  &= h({h^\star}^{-1}(f)(h^{-1}(P)+\eta))\\
  &= h(h^{-1}(f(h(h^{-1}(P)+h(\eta))))\\
  &= f(P+h(\eta))\\
  &= t_{h(\eta)}^\star(f)(P).
\end{align*}
\end{proof}

\begin{Lemma} There are $x,z\in\tK(E)$ such that $K(x)$ and $K(z)$
  are the fixed fields of $H^\star\times\gen{\beta}$ and
  $A=T^\star\rtimes(H^\star\times\gen{\beta})$, respectively.
\end{Lemma}

\begin{proof} Let $F$ be the fixed field of
  $H^\star\times\gen{\beta}$. As $H^\star$ is trivial on $K(E)$ and
  $\beta^\star$ acts on $K(E)$, we obtain that $F$ is the fixed field
  under $\beta^\star$ on $K(E)$. Let $\frakP$ be the rational place of
  $K(E)$ corresponding to the neutral element of $E(\bK)$. As $\beta$
  fixes this neutral element, the place $\frakP$ is ramified. Denote
  by $g(F)$ the genus of $F$. Of course $g(F)\le1$. But the
  Riemann-Hurwitz formula shows that there is no proper finite
  unramified extension of genus $1$ function fields. Thus $g(F)=0$.
  Also, the restriction of $\frakP$ to $F$ gives a rational place, so
  $F$ is a rational field. Rationality of the other fixed field
  follows either analogously, or by L\"uroth's Theorem.
\end{proof}

Write $z=R(x)$ with $R(X)\in K(X)$ a rational function. For suitable
choices of the above setting these functions $R$ will be
arithmetically exceptional. Thus we need a description of the
geometric and arithmetic monodromy groups.

\begin{Lemma}\label{L:ag} Denote by $\bar{\ }$ the natural
  homomorphism
\[
H^\star\times\gen{\beta^\star}\to (H^\star\times\gen{\beta^\star})/
\Cen_{H^\star\times\gen{\beta^\star}}(T^\star).
\]
(This amounts to passing from $H^\star\times\gen{\beta^\star}$ to the
induced faithful action on $T^\star$).
\begin{itemize}
\item[(a)] The geometric monodromy group of $R(X)$ is
  $T^\star\rtimes\overline{\gen{\beta^\star}}$.
\item[(b)] The arithmetic monodromy group of $R(X)$ is
  $T^\star\rtimes\overline{(H^\star\times\gen{\beta^\star})}$.
\end{itemize}
\end{Lemma}

\begin{proof} By construction, the normal hull $L$ of $K(x)$ over $K(z)$
  is contained in $\tK(E)$. The fix group of $L$ is the core of
  $H^\star\times\gen{\beta^\star}$ in
  $T^\star\rtimes(H^\star\times\gen{\beta^\star})$, which is
  obviously the centralizer of $T^\star$ in
  $H^\star\times\gen{\beta^\star}$. From that the assertion about the 
  arithmetic monodromy group follows. Similarly, the geometric
  monodromy group is the Galois group of $\tK L|\tK(z)$, and the
  Galois correspondence gives the claim.
\end{proof}

Next we give an explicit description of the field of constants in $L$
in terms of $N$.

\begin{Lemma}\label{L:foc} Let $\hat K=\bar K\cap L$ be the field of
  constants in the normal hull $L$ of $K(x)|K(z)$. Then $\hat
  K=K(\psi(N))$, where $\psi:E\to\PP$ is the morphism corresponding
  to the inclusion $K(x)\subset K(E)$.
\end{Lemma}

In the proof of this lemma, we need a purely group theoretic
result.

\begin{Lemma}\label{L:Uab} Let $U$ be a finite nilpotent group, and
  $\alpha,\beta\in\Aut(U)$ with $\alpha\beta=\beta\alpha$. Suppose
  that for each $u\in U$ there is an integer $i$ such
  $\alpha(u)=\beta^i(u)$. Then there is an integer $j$ such that
  $\alpha=\beta^j$.
\end{Lemma}

\begin{proof} The proof is by induction on $\abs{U}+\abs{\beta}$.
  We reduce to the case that $\abs{\beta}$ is a power of a prime. For
  let $r,s>1$ be two relatively prime divisors of $\abs{\beta}$. By
  induction the assertion is true if we replace the automorphisms by
  their $r$th or $s$th powers. So $\alpha^r=\beta^{rj_1}$ and
  $\alpha^s=\beta^{sj_2}$. But $\alpha=(\alpha^r)^m(\alpha^s)^n$ for
  suitable integers $m$ and $n$.  Let $\abs{\beta}$ be a power of
  the prime $p$. Similarly as above, considering the $p$th powers shows
  that $\alpha^p=\beta^{pj}$. We may replace $\alpha$ by
  $\alpha\beta^{-j}$, which amounts to assuming that $\alpha$ has
  order $p$.
  
  Next suppose that $U$ is a direct product of two proper
  $\beta$-invariant normal subgroups $U_1$ and $U_2$. Because
  $\abs{\beta}$ is a power of a prime, the group $\gen{\beta}$ acts
  faithfully on one of the $U_i$, say on $U_1$. By induction, there is
  a $j$ with $\alpha=\beta^j$ on $U_1$. Choose $u_1\in U_1$ such that
  the $\gen{\beta}$-orbit is regular through $u_1$. Let $u_2\in U_2$
  be arbitrary. The assumption gives an $i$ (depending on $u_2$) with
  $\alpha\beta^{-i}(u_1u_2)=1$. So
  $\alpha\beta^{-i}(u_1),\alpha\beta^{-i}(u_2)\in U_1\cap U_2=\{1\}$.
  Then $\beta^i=\beta^j$ by the assumption on $u_1$. We obtain
  $\alpha\beta^{-j}(u_2)=1$, so $\alpha=\beta^j$ on all of $U$.
  
  So we are reduced to the situation that $U$ is a $q$-group for some
  prime $q$. First suppose that $q\ne p$. Let $\Phi(U)$ be the
  Frattini subgroup of $U$ (i.e.\ the intersection of the maximal
  subgroups of $U$). If $\Phi(U)>1$, then by induction there is a $j$
  with $\alpha\beta^{-j}$ the trivial automorphism on $U/\Phi(U)$. But
  then $\alpha\beta^{-j}$ is trivial on $U$ as well, because the order
  of $\alpha\beta^{-j}$ is relatively prime to $q$ (see \cite[Theorem
  5.1.4]{Gorenstein}). If however $\Phi(U)=1$, then $U$ is completely
  reducible with respect to $\beta$ (by Maschke), so $U$ is
  irreducible by the previous paragraph. Then Schur's Lemma yields the
  claim.
  
  So we are left with the case $q=p$. Let $p^m$ be the order of
  $\beta$, and let $U_0<U$ be the fix group of $\beta^{p^{m-1}}$. Then
  $\gen{\beta}$ has only faithful orbits on $U\setminus U_0$. Recall
  that $\alpha$ has order $p$. We obtain that
  $\alpha(u)=\beta^{ip^{m-1}}(u)$ for each $u\in U\setminus U_0$ and
  $i\in\{1,2,\dots,p-1\}$ depending on $u$. Thus $U\setminus U_0$
  contains a subset $M$ of size
  $\ge(\abs{U}-\abs{U_0})/(p-1)\ge\abs{U}/p$ such that
  $\alpha(u)=\beta^j(u)$ for all $u\in M$ and $j$ independent of $u$.
  Of course, $\alpha=\beta^j$ holds on the group generated by $M$. As
  $\abs{M}\ge\abs{U}/p$ and $1\not\in M$, we obtain $U=\gen{M}$ and
  the claim follows.
\end{proof}

\begin{proof}[Proof of Lemma \ref{L:foc}.] Set $\hat K=\bar K\cap
  L$. The field $\hat K$ is the fixed field in $\tK$ under
  $\Cen_{H^\star\times\gen{\beta^\star}}(T^\star)$, the fix group of
  $L$.
  
  We first show that $\psi(N)\subset L$. Let $\eta\in N$ and
  $h^\star\beta^\star$ be in the centralizer of $T^\star$. We need to
  show that $(h^\star\beta^\star)(\psi(\eta))=\psi(\eta)$. Use the
  fact that $\psi$ is defined over $K$ to compute the left hand side:
\begin{align*}
  (h^\star\beta^\star)(\psi(\eta)) &= h^\star(\psi(\beta(\eta))) \\
  &= \psi(h(\beta(\eta))).
\end{align*}
Thus we are done once we know that $h(\beta(\eta))=\eta$. Recall that 
$h^\star\beta^\star$ commutes with $t_\eta$, so
\begin{align*}
t_\eta^\star &= h^\star\beta^\star
t_\eta^\star{\beta^\star}^{-1}{h^\star}^{-1} \\
  &=  h^\star t_{\beta(\eta)}^\star{h^\star}^{-1} \\
  &=  t_{h(\beta(\eta))}^\star,
\end{align*}
hence $h(\beta(\eta))=\eta$ and the claim follows.

In order to complete the proof, we need to show that $\tK\cap
L\subseteq K(\psi(N))$. By Galois theory, this amounts to showing that
the fix group of $\psi(N)$ is contained in
$\Cen_{H^\star\times\gen{\beta^\star}}(T^\star)$. Thus suppose that
$h^\star{\beta^\star}^r\in H^\star\times\gen{\beta^\star}$ fixes
$\psi(\eta)$ for each $\eta\in N$. By a computation as above, this
implies
\[
\psi(h(\eta))=\psi(\eta).
\]
The definition of $\psi$ then shows that for each $\eta\in N$, there
is an exponent $i$, depending on $\eta$, such
$h(\eta)=\beta^i(\eta)$. By Lemma \ref{L:Uab}, however, this exponent 
$i$ does not depend on $\eta$, so $h\beta^{-i}$ is the identity map
on $N$. As above, this is equivalent to $h^\star{\beta^\star}^{-i}$
to commute with each $t^\star\in T^\star$.
\end{proof}

We need a lower bound in terms of $H$ and $\beta$ for the field of
constants.

\begin{Lemma}\label{L:beta} There is a canonical injection
  $\Gal(K(N)|K(\psi(N)))\to\overline{\gen{\beta}}$. In particular,
  $\abs{\Gal(K(\psi(N))|K)}\ge\abs{\Gal(K(N)|K)}/\abs{\beta}$.
\end{Lemma}

\begin{proof} Recall that
  $\Cen_{H^\star\times\gen{\beta^\star}}(T^\star)$ is the fix group of
  $L$. The fix group of $\tK L$ is
  $\Cen_{\gen{\beta^\star}}(T^\star)$. Clearly,
  $\Gal(K(N)|K(\psi(N)))\isom\Gal(\tK
  L|L)=\Cen_{H^\star\times\gen{\beta^\star}}(T^\star)/
  \Cen_{\gen{\beta^\star}}(T^\star)$. But $\Cen_{H^\star}(T^\star)=1$
  by the definition of $H^\star$ and Lemma \ref{L:act}, so the map
  coming from projection to the second component
  $\Cen_{H^\star\times\gen{\beta^\star}}(T^\star)/
  \Cen_{\gen{\beta^\star}}(T^\star)\to
  \gen{\beta^\star}/\Cen_{\gen{\beta^\star}}(T^\star)$ is injective.
\end{proof}

The following result is more or less well known, for an analytic proof
(using Weierstra\ss\ $\wp$-functions) see \cite{Fried:CM}.

\begin{Lemma}\label{L:ram} Suppose that $\abs{N}$ is relatively prime
  to $\abs{\beta}$. Then the ramification of the extension
  $\bK(x)|\bK(R(x))$ is Euclidean of type, $(2,2,2,2)$ if
  $\abs{\beta}=2$, $(2,3,6)$ if $\abs{\beta}=6$, $(3,3,3)$ if
  $\abs{\beta}=3$, and $(2,4,4)$ if $\abs{\beta}=4$.
\end{Lemma}

\begin{proof} Let $\gamma$ be a power of $\beta$ of order $m\ge1$. We
  claim that a fixed point in $E(\bK)$ under $\gamma$ is an
  $m$-torsion point. By Lemma \ref{L:betaE}
  $1+\gamma+\gamma^2+\dots+\gamma^{m-1}=0$, so if $P\in E(\bK)$ is
  fixed under $\gamma$, then $0=mP$, the claim follows.
  
  The branched places of $\bK(x)$ in $\bK(E)$ correspond to the images
  under $\psi$ of the non-regular orbits of $\gen{\beta}$ on $\bK(E)$.
  As $\abs{N}$ is relatively prime to $\abs{\beta}$, the non-regular
  orbits are mapped bijectively to the non-regular orbits under
  $\gen{\beta'}$ on $E'(\bK)$. Note that $\Phi$ is unramified.
  Therefore the place $z\mapsto\lambda$ of $\bK(z)$ is ramified in
  $\bK(x)$ if and only if this place is branched in $\bK(E')$. Suppose
  this place is branched in $\bK(E')$, with ramification index $e$
  ($e$ is well-defined, because this extension is Galois). The above
  considerations then imply that $R^{-1}(\lambda)$ has one single
  point, and the remaining ones have multiplicity $e$. (This is
  geometrically clear, algebraically one may use Abyhankar's lemma
  \cite[Prop.\ III.8.9]{Stich}.) Furthermore, if $\mu$ is the simple
  point in the fiber $R^{-1}(\lambda)$, then the ramification indices
  of $z\mapsto\lambda$ in $\bK(E')$ and $x\mapsto\mu$ in $\bK(E)$ are
  equal. Thus the assertion follows easily if we exhibit the functions
  $\beta$ and $\psi$ in the respective cases:
  
  $\abs{\beta}=2$. Here $K(E)$ is the function field of a curve given
  by $y^2=x^3+Ax+B$, where $A,B\in K$. Then $\beta((x,y))=(x,-y)$, and
  $\psi((x,y))=x$ defines a map $\psi$ as required.
  
  $\abs{\beta}=3$. Then $K(E)$ is the function field of $y^2=x^3+B$,
  $B\in K$. Then $\beta((x,y))=(\omega x,y)$, where $\omega$ is a
  primitive $3$rd root of unity, and $\psi((x,y))=y$. Thus the places
  of $\bK(y)$ which are branched in $\bK(E)$ are $y\mapsto\infty$,
  $y\mapsto\sqrt{B}$, and $y\mapsto-\sqrt{B}$, all with ramification
  index $3$.
  
  $\abs{\beta}=6$. We have the same curve as in the previous case,
  with $\beta((x,y))=(\omega x,-y)$ and $\psi((x,y))=y^2$. The places
  of $\bK(y^2)$ which are branched in $\bK(E)$ are $y^2\mapsto\infty$,
  $y^2\mapsto B$, $y^2\mapsto0$ of ramification indexes $6$, $3$, and
  $2$, respectively.
  
  $\abs{\beta}=4$. Then $K(E)$ is the function field of $y^2=x^3+Ax$,
  $A\in K$. Then $\beta((x,y))=(-x,iy)$, where $i^2=-1$, and
  $\psi((x,y))=x^2$. Thus the places of $\bK(x^2)$ which are branched
  in $\bK(E)$ are $x^2\mapsto\infty$, $x^2\mapsto0$, and
  $x^2\mapsto-A$ of ramification indexes $4$, $4$, and $2$,
  respectively.
\end{proof}

\subsubsection{Fields of constants}
In order to get arithmetical exceptionality of $R$ we need to assure
that the field of constants $\hat K$ is sufficiently big, in
particular it need to be bigger than $K$. The following proposition
shows that there are actually many instances for this to happen. For
an elliptic curve $E$ over $K$ and a prime $n$ denote by $E[n]$ the
$n$--division points and by $K(E[n])$ the field generated by the
nonzero coordinates of the elements in $E[n]$ over $K$. Note that
$\Gal(K(E[n])|K)$ acts faithfully as a subgroup of $\GL_2(\ZZZ/n\ZZZ)$
on $E[n]\cong C_n\times C_n$ (see \cite[III, \S7]{Silverman1}).

\begin{Proposition}\label{P:gl2} Let $K$ be a number field, $t$ a
  transcendental, $E_t$ be the elliptic curve defined by
  $Y^2=X^3-t(18X-1)$, and $p$ a prime. Then for infinitely many
  $t_0\in K$ the following holds: $E_{t_0}$ is an elliptic curve over
  $K$ with $\SL_2(p)\le\Gal(K(E_{t_0}[p])|K)\le\GL_2(p)$. If $K=\QQQ$,
  then $\Gal(K(E_{t_0}[p])|K)=\GL_2(p)$ for infinitely many
  $t_0\in\QQQ$.\end{Proposition}

\begin{proof} The $j$-invariant of $E(t)$ is easily computed to be
  $j=3456t/(2t-1)$, so $K(j)=K(t)$. It follows that
  $\Gal(K(t)(E_t[p])|K(t))$ contains $\SL_2(p)$ (or equals $\GL_2(p)$
  if $K=\QQQ$), see \cite[6, \S3, Corollary 1]{Lang:EF}. The claim
  then follows from Hilbert's irreducibility theorem.\end{proof}

For the construction of arithmetically exceptional rational functions
of prime degree or of types different from $(2,2,2,2)$, we need some
consequences from the theory of complex multiplication
\cite{Silverman2} and class field theory \cite{Lang:AN}.

Let $E$ be an elliptic curve. In order to state the first result we
need to introduce the Weber function which maps $E$ to $E/\Aut(E)$.
More explicitly, assume that $E$ is given by a Weierstra\ss\ equation
form $Y^2=X^3+aX+b$. Then $h((x,y)):=x^{\abs{\Aut(E)}/2}$.

\begin{Lemma}
\label{L:cm1}
Let $K$ be an imaginary quadratic field. Let $U$ be the group of units
in $\oK$, and $E$ be an elliptic curve defined over $K$, which we
assume to have complex multiplication by $\oK$. So $U\cong\Aut(E)$.
\begin{itemize}
\item[(a)] Let $\Phi$ be the multiplication by $p$ map for an odd
  prime $p$, which is not ramified in $\oK$. Identify $E[\Phi]$ with
  the vector space $\FFF_p^2$. Then
\[
\Gal(K(E[\Phi])|K)\le
\begin{cases} C_{p^2-1} & \text{if $p$ is inert in $\oK$},\\
  C_{p-1}\times C_{p-1} & \text{if $p$ splits in $\oK$},
\end{cases}
\]
where $C_{p^2-1}<\GL_2(p)$ is a conjugate of a
Singer group, and $C_{p-1}\times C_{p-1}$ is a conjugate of the group
of diagonal matrices.

Furthermore
\[
\Gal(K(h(E[\Phi]))|K)\isom
\begin{cases} C_{p^2-1}/C_{\abs{U}} & \text{if $p$ is inert in $\oK$},\\
  (C_{p-1}\times C_{p-1})/C_{\abs{U}} & \text{if $p$ splits in $\oK$}.
\end{cases}.
\]

\item[(b)] If $\Phi$ is an endomorphism of odd prime degree $p$, then
\begin{align*}
\Gal(K(E[\Phi])|K) &\le C_{p-1}\\
\Gal(K(h(E[\Phi]))|K) &\isom C_{p-1}/C_{\abs{U}}.
\end{align*}
\end{itemize}
\end{Lemma}

\begin{proof} Note that $K$ has class number $1$ by \cite[Theorem
  4.3]{Silverman2}. We will use that below. Let $\kappa\in\oK$ induce
  $\Phi$, and let $\frakm=\kappa\oK$ be the ideal in $\oK$ generated
  by $\kappa$. Then $E[\Phi]$ is a free $\oK/\frakm$-module of rank
  $1$ (see \cite[Proposition II.1.3]{Silverman2}), where
  $\Gal(K(E[\Phi])|K)$ respects this structure. Thus
  $\Gal(K(E[\Phi])|K)$ is a subgroup of the group of units of
  $\oK/\frakm$, and the assertion about $\Gal(K(E[\Phi])|K)$ follows.
  
  Next we compute $\Gal(K(h(E[\Phi]))|K)$. By \cite[Theorem
  II.5.6]{Silverman2}, $K(h(E[\Phi]))$ is the ray class field of $K$
  with respect to the ideal $\frakm$.
  
  So we are left to compute the Galois group of this ray class field.
  For this we use global class field theory. Let
  $J=\sideset{}{'}\prod_\nu K_\nu^\star$ be the idele group of $K$,
  where $\nu$ is running over the finite places of $K$, and $K_\nu$
  are the corresponding completions. (In our context, we may omit the
  infinite places of $K$.) Note that the product is restricted, that
  is almost all components have to lie in the valuation rings of
  $K_\nu$.  Thus we still have $K^\star$ diagonally embedded in $J$.
  Let $J^\frakm$ be the subgroup of $J$ such that the elements at the
  $\nu$--th positions are restricted to be in $1+\frakp^{m_\frakp}$,
  where $\frakp$ is the associated valuation ideal, and $m_\frakp$ the
  exponent of $\frakp$ in $\frakm$. As usual, let $1+\frakp^0$ be the
  group of elements in $K_\nu$ of $\nu$-adic norm $1$. Passing to
  class groups, set $C=J/K^\star$, $C^\frakm=J^\frakm
  K^\star/K^\star$. Now the ray class field of $K$ modulo $\frakm$ has
  Galois group $C/C^\frakm$ over $K$.  But $C/C^1=1$, as $K$ has class
  number $1$. So the the Galois group is isomorphic to
\begin{align*}
  C^1/C^\frakm &\cong J^1K^\star/J^\frakm K^\star\\
  &\cong J^1/(J^1\cap J^\frakm K^\star)\\
  &=     J^1/(J^\frakm U)\\
  &\cong (J^1/J^\frakm)/((J^\frakm U)/J^\frakm).
\end{align*}

But $(J^\frakm U)/J^\frakm\cong U$ as $J^\frakm\cap U=1$ under our
assumptions, so the claim follows from noting that
$J^1/J^\frakm\isom(\oK/\frakm)^\star$.\end{proof}

The previous result gave only upper (and lower) bounds for the Galois
group of $K(E[\Phi])$ over $K$. The upper bound is not necessarily
achieved, however it is in most of the cases. In order to formulate
the next result, we say that an elliptic curve with a given
Weierstra\ss\ equation over $K$ has good reduction at $p\in\PPP$, if
the reduction modulo each prime ideal dividing $p\oK$ is a
non-singular curve. The following result, though it looks very similar
to the previous one, is much deeper, its proof uses the main theorem
of complex multiplication. See also \cite[5.20(ii)]{Rubin:CM} for a
similar statement.

\begin{Lemma}\label{L:cm2}
  Let $K$ be an imaginary quadratic field, and $E$ be an elliptic
  curve defined over $K$, which we assume to have complex
  multiplication by $\oK$. Let $p$ be a prime such that $E$ has good
  reduction at $p$.
\begin{itemize}
\item[(a)] Let $\Phi$ be the multiplication by $p$, and suppose that
  $p$ is not ramified in $\oK$. Then
\[
\Gal(K(E[\Phi])|K)=
\begin{cases} C_{p^2-1} & \text{if $p$ is inert in $\oK$},\\
  C_{p-1}\times C_{p-1} & \text{if $p$ splits in $\oK$}.
\end{cases}
\]
\item[(b)] If $\Phi$ is an endomorphism of degree $p$, then
\begin{align*}
\Gal(K(E[\Phi])|K) &= C_{p-1}\\
\end{align*}
\end{itemize}
\end{Lemma}

\begin{proof} The uniformization theorem (see e.g.\
  \cite[\S5]{Silverman1}) gives a lattice $\Lambda\subset\CCC$ and an
  analytic isomorphism $f:\CCC/\Lambda\to E(\CCC)$. Without loss
  assume that $1\in\Lambda$. The fact that the endomorphism ring of
  $E$ is the maximal order $\oK$, and that this endomorphism ring is
  isomorphic to $\alpha\in K$ with $\alpha\Lambda\subseteq\Lambda$,
  implies that $\Lambda$ is a fractional ideal of $\oK$. Under this
  identification, the points of finite order in $E(\CCC)$ are given by
  $K/\Lambda$. The torsion points of $E(\CCC)$ generate an abelian
  extension of $K$, so $f$ maps $K/\Lambda$ to $E(K^{ab})$. The main
  theorem of complex multiplication describes the action of
  $\Gal(K^{ab}|K)$ on $E(K^{ab})$ in terms of the action of the idele
  group $J$ on $K/\Lambda$. As in the previous proof, let $J$ be the
  the idele group of $K$, where we again drop the factor $\CCC^\star$
  coming from the archimedian place. In order to describe the action
  of $J$ on $K/\Lambda$, we need to fix a $\frakp$-primary
  decomposition of $K/\Lambda$: Let $\nu$ run through the
  non-archimedian valuations of $K$, let $K_\nu$ be the completion of
  $K$ with respect to $\nu$, and let $\Lambda_\nu:=\Lambda R_\nu$ be
  the fractional ideal of $K_\nu$, where $R_\nu$ is the valuation ring
  of $\nu$. Then (see \cite[\S8]{Silverman2})
\[
K/\Lambda\isom\bigoplus_{\nu}K_\nu/\Lambda_\nu.
\]
Let $\psi:J\to\CCC^\star$ be the Hecke character corresponding to $E$
and $K$, and $\nu$ be a valuation of $K$ lying above the $p$-adic
valuation of $\ZZZ$. As $E$ has good reduction at $p$, \cite[Theorem
9.2]{Silverman2} gives $\psi(R_\nu^\star)=\{1\}$. Let $x\in J$ be an
idele with $x_{\nu'}=1$ for $\nu'\ne\nu$ and $x_\nu\in R_\nu^\star$.
Note that $x_\nu\Lambda_\nu\subseteq\Lambda_\nu$, so $\gen{x_\nu}$
acts on the quotient $K_\nu/\Lambda_\nu$ by multiplication with
$x_\nu$. Now use the results from \cite[\S9]{Silverman2} to see the
following:

Use the $\frakp$-primary decomposition of $K/\Lambda$ from above to
define an action of $x$ on $K/\Lambda$: Let $x_{\nu'}$ act trivially
on $K_{\nu'}/\Lambda_{\nu'}$ for $\nu'\ne\nu$, and let $x_{\nu}$ act
by multiplication with $x_\nu$ on $K_\nu/\Lambda_\nu$. Denote by
$[x,K]\in\Gal(K^{ab}|K)$ the Artin symbol of $x$. Then we have a
commutative diagram
\[
\begin{array}{ccc}
K/\Lambda & \stackrel{x^{-1}}\longrightarrow & K/\Lambda\\
\Big\downarrow f & & \Big\downarrow f\\
E(K^{ab}) & \stackrel{[x,K]}\longrightarrow & E(K^{ab})
\end{array}
\]
From that we can read off the claim.
\end{proof}

\subsubsection{The $(2,2,2,2)$-case}
We are now prepared to prove the existence of arithmetically
exceptional rational functions of elliptic type $(2,2,2,2)$ in those
cases which are not ruled out by Theorem \ref{ellmain}.

\begin{Theorem}\label{ellmainexi} There exist indecomposable
arithmetically exceptional rational functions of type
$(2,2,2,2)$
\begin{itemize}
\item[(a)] of degree $p^2$ for each odd prime $p$ over any number
field $K$. For $K=\QQQ$, there are examples with arithmetic monodromy
group as big as $\AGL_2(p)$.
\item[(b)] over $\QQQ$ of degrees $5,7,11,13,17,19,37,43,67,163$.
\item[(c)] of each prime degree $p\ge5$ for some number field $K$.
Hereby, $K$ can be chosen to be the Hilbert class field of
$\QQQ(\sqrt{-p})$, giving the arithmetic monodromy group to be
$\AGL_1(p)$.
\end{itemize}
\end{Theorem}

\begin{proof} The construction is based on isogenies $\Phi:E\to E'$ and
  $\beta\in\Aut(E)$ the canonical involution. As
  $\Phi(\beta(P))=\Phi(-P)=-\Phi(P)=\beta'(\Phi(P))$, we obtain that
  the kernel of $\Phi$ is invariant under $\beta$, so the results from
  Section \ref{S:Isog} apply. By Lemma \ref{L:ram}, the ramification
  type of the associated rational function $R(X)$ is correct. Let $A$
  and $G$ denote the arithmetic and geometric monodromy group of
  $R(X)$, respectively.
  
  To (a). Let $K=\QQQ$ and $E=E'$ be an elliptic curve from the
  conclusion of Proposition \ref{P:gl2}, and let $\Phi$ be the
  multiplication by $p$ map. By Lemma \ref{L:ag}, we have
  $A\le\FFF_p^2\rtimes\GL_2(p)$ and $G=\FFF_p^2\rtimes C_2$. On the
  other hand, $A/G=\Gal(K(\psi(N))|K)$, and the claim follows from
  Lemma \ref{L:beta} together with $\Gal(K(N)|K)=\GL_2(p)$. Similarly
  argue for $K$ a number field to obtain $\ASL_2(p)\le A\le\AGL_2(p)$.
  That clearly gives arithmetical exceptionality.
  
  To (b). First suppose that $p\in\{7,11,19,43,67,163\}$, and set
  $K=\QQQ(\sqrt{-p})$. Then $\oK=\ZZZ+\omega\ZZZ$ with
  $\omega=(1+\sqrt{-p})/2$ is the ring of integers of $K$.
  Furthermore, $K$ has class number one in these cases. Thus there
  exists an elliptic curve $E$ over $\QQQ$ with complex multiplication
  by $\oK$, see \cite[11.3.1]{Silverman1}. The multiplication by
  $\kappa=\sqrt{-p}$ yields an endomorphism $\Phi$ of degree $p$. We
  claim that the kernel is defined over $\QQQ$ (even though $\Phi$
  obviously is not). Let $\Phi'$ be a Galois conjugate of $\Phi$
  different from $\Phi$. Then $\Phi'$ is induced by multiplication
  with $-\kappa$. Passing to the lattice $\CCC/\oK$, an easy
  computation shows that the kernels of $z\mapsto\kappa z$ and
  $z\mapsto-\kappa z$ are the same. Thus $E[\Phi]$ is defined over
  $\QQQ$. By \cite[III.4.12, III.4.13.1]{Silverman1}, there exists a
  rational isogeny $E\to E/E[\Phi]$ with kernel $E[\Phi]$. As $E$ has
  no nontrivial automorphisms, $\psi$ is the Weber function on $E$, so
  $\Gal(K(\psi(N))|K)=C_{(p-1)/2}$ by Lemma \ref{L:cm1}(b). Thus the
  arithmetic monodromy group taken over $K$ is $\AGL_1(p)$, so this
  is even more true over $\QQQ$.
  
  The remaining cases are $p\in\{5,13,17,37\}$. First suppose $p=13$.
  Let $X_0(p)$ and $X_1(p)$ be the modular curves as defined in
  \cite[App.\ C, \S13]{Silverman1}.  These curves are defined over
  $\QQQ$.  Roughly, the non-cuspidal points of $X_0$ classify the
  pairs of elliptic curves $E$ and a subgroup of order $p$, with both
  defined over a field $K$ if the corresponding point on $X_0(p)$ is
  defined over $K$. Similarly for $X_1(p)$, with the group of order
  $p$ replaced by a point of order $p$. There is a natural covering
  map $\phi:X_1(p)\to X_0(p)$, which amounts to sending the pair
  $(E,P)$ with $P\in E$ of order $p$ to $(E,\gen{P})$. This covering
  map has degree $(p-1)/2$. (We need to divide by $2$, as $(E,P)$ and
  $(E,-P)$ correspond to the same point on $X_1(p)$.) Now $X_0(p)$ is
  a rational curve for $p=13$. Thus Hilbert's irreducibility theorem
  implies that there are infinitely many rational points $Q$ on
  $X_0(p)$, such that $\QQQ(Q')|\QQQ$ has degree $(p-1)/2$ where
  $Q'\in X_1(p)$ with $\phi(Q')=Q$. Such a point $Q$ thus corresponds
  to an elliptic curve over $\QQQ$ with an isogeny $\Phi$ of degree
  $p$, and with $[\QQQ(N):\QQQ]=(p-1)/2$. So
  $[\QQQ(\psi(N)):\QQQ]\ge(p-1)/4\ge3$ by Lemma \ref{L:beta}, and
  arithmetical exceptionality follows.
  
  We finish also the cases $p=17$ and $37$ once we know that there is
  an elliptic curve over $\QQQ$ with a rational isogeny $\Phi$ of
  degree $p$ and $\psi(N)\not\subset\QQQ$. For $p\in\{17,37\}$ it is
  known that there exists an elliptic curve over $\QQQ$ with a
  rational isogeny of degree $p$, see \cite{Mazur} and the references
  given there. Suppose that for such a curve $\psi(N)\subset\QQQ$. But
  then there would be a point of order $p$ defined over a quadratic
  number field. But this cannot happen by results of Kamienny, see the
  references in \cite{KamMaz}.
  
  So the only case left is $p=5$. Here we can even give an
  arithmetically exceptional function quite explicitly. Namely let
  $R(X)=X(11X^4+40X^3+10X^2-40X-5)/(5X^2-1)^2$. Set
  $q_1(X)=11X^3-5X^2-3X-11/12$ and $q_2(X)=X^3-5X^2+7X-1$. Then one
  verifies that $q_2(f(X))=q_1(X)(f'(X)/5)^2$, thus
  $(X,Y)\mapsto(f(X),1/5Yf'(X))$ gives an isogeny of the elliptic
  curve $Y^2=q_1(X)$ to $Y^2=q_2(X)$ of degree $5$. Furthermore, no
  $X$--coordinate of a nonzero element in the kernel is rational, as
  $5X^2-1$ is irreducible. This proves the claim. (We thank Ian
  Connell for pointing out on how to use his elliptic curves software
  {\em apecs}, available from http://euclid.math.mcgill.ca/~connell,
  to compute isogenies of elliptic curves.)
  
  To (c). The argument is a slight extension of the proof of Lemma
  \ref{L:cm1}, again using the theory of complex multiplication
  \cite[Chapter II]{Silverman2}. Let $K_h$ be the Hilbert class field
  of $K=\QQQ(\sqrt{-p})$. Then there is an elliptic curve defined over
  $K_h$, which admits complex multiplication by $\oK$. Let $\Phi$ be
  the endomorphism (of degree $p$) induced by $\kappa=\sqrt{-p}$. Set
  $\frakm=\kappa\oK$. The field $K_h':=K_h(\psi(N))$ is the ray class
  field of $K$ modulo $\frakm$. With the notation from the proof of
  Lemma \ref{L:cm1} we have $\Gal(K_h|K)\cong C/C^1$. So the proof there
  gives $\Gal(K_h'|K_h)\cong C_{(p-1)/2}$. This implies that the
  arithmetic monodromy group of $R$ over $K_h$ is $\AGL_1(p)$, and
  everything follows.\end{proof}

\begin{Remark*} The claim of part (c) already appears in the work of
  Fried, see \cite{Fried:CM}. His argument is completely different
  (using modular curves rather than complex multiplication) and does
  not give a hint about possible fields of definition.
  
  There is yet another proof of part (c), which does not use complex
  multiplication. Take $R(X)$ as in (a) of degree $p^2$, ($p\ge5$),
  such that the arithmetic monodromy group is $\AGL_2(p)$. So the
  field of constants $\hat\QQQ$ has Galois group $\GL_2(p)/C_2$ over
  $\QQQ$.  Now let $B$ be the group of upper triangular matrices in
  $\GL_2(p)$, and $K$ be the fixed field of $B/C_2$ in $\hat\QQQ$. So
  $f$ has arithmetic monodromy group $C_p^2\rtimes B$ over $K$, so
  $R(X)$ decomposes over $K$ into two rational functions of degree
  $p$. It is immediate to verify that the two composition factors have
  the right ramification, and arithmetic monodromy groups
  $\AGL_1(p)$.\end{Remark*}

\subsubsection{The cases $(2,3,6)$, $(3,3,3)$, and $(2,4,4)$}
The next theorems treat the ramification types $(2,3,6)$, $(3,3,3)$,
and $(2,4,4)$. The construction is based on the setting from Section
\ref{S:Isog}. Specifically, let $K$ be a number field (to be specified
in the theorems), $E$ an elliptic curve with an automorphism $\beta$,
and $\Phi:E\to E$ an endomorphism. As the endomorphism ring of an
elliptic curve is abelian \cite[Theorem 6.1]{Silverman1}, we obtain
that the kernel of $\Phi$ is invariant under $\beta$. The results from
Section \ref{S:Isog} may be interpreted that we get a rational
function $R$ such that the following diagram is commutative:
\[
\begin{array}{ccc}
E & \stackrel{\Phi}\longrightarrow & E\\
\Big\downarrow \psi & & \Big\downarrow \psi\\
E/\gen{\beta}\isom\PP & \stackrel{R}\longrightarrow & E/\gen{\beta}\isom\PP
\end{array}
\]
This is the motivation for talking about the rational function
$E/\gen{\beta}\to E/\gen{\beta}$ in the theorems below. For the sake
of easier reading, we separate the prime degree from the prime square
degree cases.

\begin{Theorem}\label{236p} Let $p$ be a prime with
  $p\equiv1\pmod{3}$. Set $K=\QQQ(\sqrt{-3})$, and let $E$ be an
  elliptic curve over $K$ with complex multiplication by $\oK$. So
  there is an automorphism $\beta$ of order $6$.  There exists an
  endomorphism $\Phi$ of $E$ of degree $p$. The associated rational
  function $R:E/\gen{\beta}\to E/\gen{\beta}$ has ramification type
  $(2,3,6)$, is defined over $K$, and has arithmetic monodromy group
  $\AGL_1(p)$.
\end{Theorem}

\begin{Theorem}\label{236p^2} Let $p>3$ be a prime, and $E$ be an
  elliptic curve over $\QQQ$ with complex multiplication by the
  integers of $\QQQ(\sqrt{-3})$. Let $\beta$ be an automorphism of
  order $6$. The rational function $R:E/\gen{\beta}\to E/\gen{\beta}$
  associated to the multiplication by $p$ map has ramification type
  $(2,3,6)$, is defined over $\QQQ$, and its arithmetic monodromy
  group over $\QQQ$ is
\begin{itemize}
\item $(C_p\times C_p)\rtimes(C_{p^2-1}\rtimes C_2)$
  if $p\equiv-1\pmod{3}$, or
\item $(C_p\times C_p)\rtimes((C_{p-1}\times
  C_{p-1})\rtimes C_2)$ if $p\equiv1\pmod{3}$.
\end{itemize}
\end{Theorem}

\begin{proof}[Proof of Theorem \ref{236p}.] By assumption $p$ is a
  norm in $\oK$, so there is an endomorphism $\Phi$ of $E$ of degree
  $p$ which is defined over $K$. Lemma \ref{L:ram} gives the correct
  ramification type of the associated rational function $R(X)$. Let
  $A,G\le\AGL_1(p)$ be the arithmetic and geometric monodromy group of
  $R(X)$ over $K$. Of course $G=C_p\rtimes C_6$. In our case, the
  function $\psi:E\to\PP$ coincides with the Weber function $h$ on
  $E$. So, by Lemma \ref{L:cm1}(b), the field of constants of $R(X)$
  has Galois group $C_{p-1}/C_6$. As this Galois group is isomorphic
  to $A/G$, and $A\le\AGL_1(p)$, the claim follows.
\end{proof}
  
\begin{proof}[Proof of Theorem \ref{236p^2}.] Now let $\Phi$ be the
  multiplication by $p$ map, and $N=E[\Phi]$ the kernel. As $\Phi$ and
  $\psi$ are defined over $\QQQ$, we obtain that $R(X)\in\QQQ(X)$,
  even though $\beta$ is not defined over $\QQQ$.  We first work over
  $K=\QQQ(\sqrt{-3})$, and afterwards descend to $\QQQ$. Let
  $(C_p\times C_p)\rtimes C$ with $C\le\GL_2(p)$ be the arithmetic
  monodromy group of $R(X)$ taken over $K$. The field of constants is
  $K(\psi(N))$. First consider the case $p\equiv-1\pmod{3}$. Note that
  $C_6\le C$, and $C$ is abelian by Lemma \ref{L:ag}(b) and Lemma
  \ref{L:cm1}(a). The Galois group of $K(\psi(N))|K$ is given by Lemma
  \ref{L:cm1}(a). We obtain $C/C_6\isom C_{p^2-1}/C_6$. So $C$ is
  abelian of order $p^2-1$, and irreducible (otherwise by Maschke $C$
  were a subgroup of $C_{p-1}\times C_{p-1}$, so $p+1$ would divide
  $p-1$).  Thus by Schur's Lemma $C$ is cyclic, indeed $C$ is a Singer
  group in $GL_2(p)$.
  
  If $p\equiv1\pmod{3}$, then we obtain similarly
  $C/C_6\isom(C_{p-1}\times C_{p-1})/C_6$. If the abelian group $C$
  were irreducible, then $\abs{C}=(p-1)^2$ would divide $p^2-1$, a
  contradiction. Thus $C\isom C_{p-1}\times C_{p-1}$, the group of
  diagonal matrices in $GL_2(p)$.
  
  The assertion about the arithmetic monodromy group over $\QQQ$
  follows now from these two observations: $\psi$ is defined over
  $\QQQ$, thus the extension $K(\psi(N))|\QQQ$ is normal. Secondly,
  $K$ is contained in the field of constants of $R(X)$, for instance
  by the branch cycle argument Proposition \ref{BCA}, noting that a
  generator of $C_6\le C$ is not rational in the geometric monodromy
  group $(C_p\times C_p)\rtimes C_6$. Thus, if $(C_p\times C_p)\rtimes
  \hat C$ is the arithmetic monodromy group over $\QQQ$, then $C$ has
  index two in $\hat C$. Is is easy to see that this extension split.
  Indeed, $\hat C$ may be identified with
  $\FFF_{p^2}^\star\rtimes\Aut(\FFF_{p^2})$ in the first case, or with
  the group of monomial matrices in $GL_2(p)$ in the second case.
\end{proof}

Though one might expect that the situation for the ramification type
$(3,3,3)$ is completely analogous to the previous case, there is a
difficulty with bad reduction. Indeed, the arithmetic monodromy group
of $R$ depends on the chosen elliptic curve, see Example \ref{E:cm7}
and Remark \ref{R:333}. Possible choices for $E$ in the following two
theorems are $Y^2=X^3+B$, where $B\in\ZZZ$ is not divisible by $p$.

\begin{Theorem}\label{333p} Let $p$ be a prime with
  $p\equiv1\pmod{3}$. Set $K=\QQQ(\sqrt{-3})$, and let $E$ be an
  elliptic curve over $K$ with complex multiplication by $\oK$ and
  good reduction at $p$. So there is an automorphism $\beta$ of order
  $3$. Let $\Phi$ be an endomorphism of $E$ of degree $p$. The
  associated rational function $R:E/\gen{\beta}\to E/\gen{\beta}$ has
  ramification type $(3,3,3)$, is defined over $K$, and has arithmetic
  monodromy group $\AGL_1(p)$.
\end{Theorem}

\begin{Theorem}\label{333p^2} Let $p>3$ be a prime, and $E$ be an
  elliptic curve over $\QQQ$ with complex multiplication by the
  integers of $\QQQ(\sqrt{-3})$ and good reduction at $p$. Let $\beta$
  be an automorphism of order $3$. The rational function
  $R:E/\gen{\beta}\to E/\gen{\beta}$ associated to the multiplication
  by $p$ map has ramification type $(3,3,3)$, is defined over $\QQQ$,
  and its arithmetic monodromy group over $\QQQ$ is
\begin{itemize}
\item $(C_p\times C_p)\rtimes(C_{p^2-1}\rtimes C_2)$
  if $p\equiv-1\pmod{3}$, or
\item $(C_p\times C_p)\rtimes((C_{p-1}\times
  C_{p-1})\rtimes C_2)$ if $p\equiv1\pmod{3}$.
\end{itemize}
\end{Theorem}

\begin{proof}[Proof of Theorem \ref{333p}.] Use the notation from
  the proof of Theorem \ref{236p}, with the only modification that
  $\beta$ has order $3$. Now $G=C_p\rtimes C_3$, and $A\le\AGL_1(p)$.
  So we are done once we show that the degree of the field of
  constants $K(\psi(N))$ over $K$ is $\ge(p-1)/3$, where $N$ is the
  kernel of $\Phi$. But $[K(N):K]=p-1$ by Lemma \ref{L:cm2}, and
  $[K(\psi(N)):K]\ge [K(N):K]/\abs{\beta}$ by Lemma \ref{L:beta}, so
  the claim follows.
\end{proof}

\begin{proof}[Proof of Theorem \ref{333p^2}] Combine the ideas of
    the proofs of Theorems \ref{333p} and \ref{236p^2}.
\end{proof}

\begin{Theorem}\label{244p} Let $p$ be a prime with
  $p\equiv1\pmod{4}$. Set $K=\QQQ(\sqrt{-1})$, and let $E$ be an
  elliptic curve over $K$ with complex multiplication by $\oK$. So
  there is an automorphism $\beta$ of order $4$.  There exists an
  endomorphism $\Phi$ of $E$ of degree $p$. The associated rational
  function $R:E/\gen{\beta}\to E/\gen{\beta}$ has ramification type
  $(2,4,4)$, is defined over $K$, and has arithmetic monodromy group
  $\AGL_1(p)$.
\end{Theorem}

\begin{Theorem}\label{244p^2} Let $p$ be an odd prime, and $E$ be an
  elliptic curve over $\QQQ$ with complex multiplication by the
  integers of $\QQQ(\sqrt{-1})$. Let $\beta$ be an automorphism of
  order $4$. The rational function $R:E/\gen{\beta}\to E/\gen{\beta}$
  associated to the multiplication by $p$ map has ramification type
  $(2,4,4)$, is defined over $\QQQ$, and its arithmetic monodromy
  group over $\QQQ$ is
\begin{itemize}
\item $(C_p\times C_p)\rtimes(C_{p^2-1}\rtimes C_2)$
  if $p\equiv-1\pmod{4}$, or
\item $(C_p\times C_p)\rtimes((C_{p-1}\times C_{p-1})\rtimes C_2)$ if
  $p\equiv1\pmod{4}$.
\end{itemize}
\end{Theorem}

\begin{proof}[Proof of Theorems \ref{244p} and \ref{244p^2}.] The
  proofs are completely analogous to the proofs of the $(2,3,6)$
  cases. In order to prove Theorem \ref{244p^2} we need to know that
  $\sqrt{-1}$ is in the field of constants. As previously, that
  follows from the branch cycle argument. Note that the branch points
  of the function $R(X)$ are all rational by the proof of Lemma
  \ref{L:ram}, even though the two points with inertia index $4$ could
  a priori be algebraically conjugate.
\end{proof}

\begin{Example}\label{E:cm7} Let $\omega$ be a primitive $3$rd root of
  unity, and $K=\QQQ(\omega)=\QQQ(\sqrt{-3})$. Let $E$ be the elliptic 
  curve $Y^2=X^3+B$ with $B\in\oK=\ZZZ[\omega]$. Set $p=7$. Clearly
  $E$ fulfills the assumptions of Theorem \ref{333p} except possibly
  the one on good reduction at $7$. The norm of $3+\omega$ is $7$, so
  the endomorphism $[3]+\beta$ of $E$ has degree $7$. Using the
  addition formula, one can explicitly compute this endomorphism. It
  sends a point $(x,y)\in E$ to $(*,R(y))$ with
\[
R(Y)=\frac{(1-18\omega)(y^6+(9+108\omega)By^4+(459+216B^2\omega)B^2y^2
-(405+324\omega)B^3)y}{(7y^2-(3-12\omega)B)^3}.
\]
The field of constants is generated by the roots of the denominator of 
$R(Y)$. Thus the arithmetic monodromy group of $R$ over $K$ is
$\AGL_1(7)$, unless this denominator splits, in which case we obtain
the index $2$ subgroup of $\AGL_1(7)$. This latter case happens
whenever $B=5r^2-8rs-s^2+(4r^2+2rs-5s^2)\omega$ with $r,s\in\ZZZ$.
\end{Example}

\begin{Remark}\label{R:333} We do not know whether the good reduction
  assumption is really needed in Theorem \ref{333p^2}. One can show
  that this assumption is superfluous if $p\equiv{-1}\pmod{3}$. Namely
  let $E$ be the curve $Y^2=X^3+B$. Then the multiplication by $p$ map
  sends $(x,y)$ to $(*,\frac{U(y)}{V(y)})$, where $V$ is essentially
  the $p$-division polynomial, expressed in terms of $y$. Assume that
  $V$ is monic. We have $V(Y)=h(Y^2)$, where $h\in\QQQ[Y]$ has degree
  $(p^2-1)/6$. Note that $Y^2$ is the Weber function on $E$, so Lemma
  \ref{L:cm1} shows that $h(Y)$ is irreducible over $K=\QQQ(\omega)$.
  On the other hand, the proof of the theorem fails exactly if
  $h(Y^2)$ becomes reducible over $K$. If that happens, then
  $h(Y^2)=u(Y)u(-Y)$, where $u\in K[Y]$. Setting $Y=0$ shows that
  $V(0)=h(0)$ is a square in $K$. On the other hand, the recursion
  formulae for division polynomials (see \cite[Chapter II]{Lang:EC},
  \cite[Exercise 3.7]{Silverman1}) allow to explicitly compute the
  constant term of $V(Y)$: Is is $\pm3^{(p^2-1)/4}B^{(p^2-1)/6}/p$.
  Note that the exponents $(p^2-1)/6$ and $(p^2-1)/4$ are even, so
  this coefficient is never a square in $K$.
  
  Unfortunately, this kind of reasoning does not seem to extend to the
  case $p\equiv{1}\pmod{3}$. Here $h(Y)$ is reducible and the absolute
  term of the factor which should stay irreducible upon replacing $Y$
  with $Y^2$ turns out to be a square if $p\equiv1\pmod{12}$!
\end{Remark}

\section{Sporadic cases of arithmetic exceptionality}\label{g>1}

We are now going through the list of sporadic group theoretic
candidates for arithmetically exceptional functions suggested by
Theorems \ref{AS} and \ref{Affine} which do not belong to the series
we dealt with, decide in each case whether they occur over some number
field, and if they do, determine the minimal field of definition.

\subsection{$G=C_2\times C_2$ (Theorem \ref{Affine}(a)(iii))}\label{A4}

This is a very rich case, as every number field can be a field of
definition, and also each $C_3$ or $S_3$ extension can be the field of
constants. This is in big contrast to all the other cases. More
precisely, we have the following. We thank Mike Zieve for an idea
which helped constructing the function $f$ below.

\begin{Theorem}\label{A4S4} Let $K$ be any field of characteristic
  $0$, and $E$ a Galois extension of $K$ with group $C_3$ or $S_3$.
  Then there exists a rational function $f(X)\in K(X)$ of degree $4$,
  such that the geometric monodromy group $G$ is $C_2\times C_2$, the
  arithmetic monodromy group is $A=A_4$ or $S_4$, and $E$ is the field
  of constants.\end{Theorem}

\begin{proof} Let $E_0$ be a cubic extension of $K$ inside $E$, and
  $Z^3+pZ+q\in K[Z]$ be a minimal polynomial of a primitive element of 
  $E_0|K$. Set
\[
f(X) := \frac{X^4-2pX^2-8qX+p^2}{4(X^3+pX+q)}.
\]
We claim that $f$ has the required properties. Let $\lambda$ be a root 
of $X^3+pX+q$. Then one verifies that
\[
f(X)-\lambda = \frac{(X^2-2\lambda X-2\lambda^2-p)^2}{4(X^3+pX+q)}.
\]
This shows that the roots of $X^3+pX+q$ are branch points of $f$, of
cycle type $(2^2)$. Thus the geometric monodromy group of $f$
is $C_2\times C_2$. Looking at the fiber $f^{-1}(\infty)$ shows that
the field of constants is the splitting field of $X^3+pX+q$, thus
$E$.\end{proof}

The above theorem has an interesting consequence about compositions of
arithmetically exceptional functions, see also Corollary
\ref{RedeiCor}.

\begin{Corollary} Let $K$ be any number field. Then there are $4$
arithmetically exceptional functions of degree $4$ over $K$, such that
their composition is not arithmetically exceptional.\end{Corollary}

\begin{proof} Let $L$ be a Galois extension of $K$ with group
$C_3\times C_3$. (It is well known that over a number field any finite
abelian group is a Galois group). For $i=1,2,3,4$ let $E_i$ be the
fixed fields of the $4$ subgroups of order $3$ in $\Gal(L|K)$. For
each extension $E_i|K$ there is, by the previous theorem, an
arithmetically exceptional function over $K$ of degree $4$ with $E_i$
being the field of constants. To conclude the assertion, we only need
to show (confer Theorem \ref{togroups}) that every rational prime $p$
splits completely in at least one of the fields $E_i$. Let $\frakP$ be
a prime of $L$ above $p$, and $\sigma\in\Gal(L|K)$ be a generator of
the corresponding decomposition group. Of course, $\sigma$ has either
order $1$ (and we are done), or lies in at least one of the groups
$\Gal(L|E_i)$. But then $p$ splits in the corresponding extension
$E_i|K$.\end{proof}

\subsection{$G=(C_{11}^2)\rtimes\GL_2(3)$ (Theorem
  \ref{Affine}(c)(i))} Let $\cT=(\s1,\s2,\s3)$ be a genus $0$ system
of $G$, with the orders of $\s1$, $\s2$, $\s3$ being $2$, $3$, and
$8$, respectively. Let $\cC_i$ be the conjugacy class $\s i^G$. One
verifies that $\cC_1$ and $\cC_2$ are rational in $G$, whereas
$\cC_3=\cC_3^3\ne\cC_3^5=\cC_3^7$. Let $\zeta$ be a primitive $8$th
root of unity. Suppose that we have a realization of these data over
the number field $K$.

\begin{Claim} $\QQQ(\sqrt{-2})=\QQQ(\zeta+\zeta^3)\subseteq
K$.\end{Claim}

\begin{proof} Suppose not. Then $K$ and $\QQQ(\zeta+\zeta^3)$ are
  linearly disjoint over $\QQQ$, so there is an automorphism of
  $\bQ|\QQQ$ which is trivial on $K$, but moves $\beta:=\zeta+\zeta^3$
  to $\zeta^5+\zeta^7$. So the number $m$ in Proposition \ref{BCA} is
  congruent $5$ or $7$ modulo $8$. However, $A$ is the full
  automorphism group of $G$, but $\cC_3\ne\cC_3^m$, a
  contradiction.\end{proof}

Now set $K:=\QQQ(\sqrt{-2})$. The above proof actually shows that
$\cH:=(\cC_1,\cC_2,\cC_3)$ is $K$--rational. One verifies that $\cH$
is weakly rigid. Also, $G$ contains only $1$ conjugacy class of
subgroups of index $121$ (for instance by the Schur--Zassenhaus
Theorem), and $\s2$ has a unique fixed point, so Proposition
\ref{WeakRigCrit} applies. So there is a rational function $f(X)\in
K(X)$ with geometric monodromy group $G$. The symmetric normalizer of
$G$ is $A$, with $A/G=C_5$. Hence the arithmetic monodromy group of
$f$ is either $G$ or $A$. We show that the latter holds.

We start by assuming that the arithmetic monodromy group of $f$ is
$G$, and derive a contradiction. Let $L$ be a splitting field of
$f(X)-t$ over $K(t)$, hence $G=\Gal(L|K(t))$. There are subgroups
$V<M<G$ (unique up to conjugacy), such that the following holds:
$[G:M]=12$, $[M:V]=11$, and $\cT$ induces a genus $0$ system on $G/V$.
Without loss assume that the branch points belonging to $\s1$, $\s2$,
and $\s3$ are $-27/256$, $0$, and $\infty$. We now identify the genus
$0$ fixed fields of $M$ and $V$. First note that there is a subgroup
$W$ with $M<W<G$, $[G:W]=4$, and the action of $G$ on $G/W$ is the
symmetric group $S_4$. The cycle types of the $\s i$ are $(1^22^1)$,
$(1^13^1)$, and $(4^1)$, respectively. Thus, there is $x\in L_W$, the
fixed field of $W$, such that $t=x^3(x-1)$. Then $M$ is the stabilizer
of $K(x,y)$ for a conjugate $y$ of $x$. Of course, also $t=y^3(y-1)$.
Note that
\begin{align*}
x &= \frac{r(r^3-1)}{r^4-1}\\
y &= \frac{r^3-1}{r^4-1}
\end{align*}
gives a parametrization of the non--diagonal part of the curve
$x^3(x-1)-y^3(y-1)=0$. Also $r=x/y$ so $K(x,y)=K(r)$. We get
\[
t=\frac{r^3(1-r)(r^3-1)^3}{(r^4-1)^4}.
\]
The linear fractional change
\[
r=\frac{\sqrt{-2}-s}{\sqrt{-2}+s}
\]
yields
\[
t=\frac{1}{1024}\,\frac{(s^2-6)^3(s^2+2)^3}{(s^2-2)^4}=:g(s).
\]
For later use we record
\begin{equation}\label{27}
t+\frac{27}{256}=\frac{1}{1024}\,\frac{(s-2)(s+2)
[s(s^2-4s+6)(s^2+4s+6)]^2}{(s^2-2)^4}.
\end{equation}
Let $L_V$ be the fixed field of $V$. The cycle types of the $\s i$ on
$G/M$ are $(1^22^5)$, $(3^4)$, and $(4^3)$, whereas they are
$(1^{12}2^{60})$, $(3^{44})$, and $(4^38^{15})$ on $G/V$. From that we
see that the branch points of $L_V$ over $L_M=K(s)$ are the $3$
quadruple points in the fiber $g^{-1}(\infty)$, and one of the two
simple points in the fiber $g^{-1}(-27/256)$, and that the cycle types
in each of these points is $(1^12^5)$. Thus we get a $K$--rational
isogeny of degree $11$ of an elliptic curve $E_1$ defined over $K$ to
the elliptic curve $E_2:\,Y^2=(X\pm2)(X^2-2)$, see Section \ref{g=1}.
Upon possibly replacing $s$ with $-s$ (which does not change $g(s)$),
we may assume that the elliptic curve in question is
\[
E_2:\,Y^2=(X-2)(X^2-2).
\]
\begin{Lemma} Let $\phi:E_1\to E_2$ be a $K$--rational isogeny of
  degree $11$. Then $E_1$ and $E_2$ are $K$--isomorphic, and this
  isogeny is induced by complex multiplication with an integer in
  $K$.\end{Lemma}

\begin{proof}
  Let $\hat\phi:\,E_2\to E_1$ be the dual isogeny. We first show that
  there are at most two isogenies of degree $11$, which are
  interchanged by the action of $\Gal(K|\QQQ)$. Let $\psi_{11}(X)$ be
  the $11$th division polynomial. This can easily be computed using
  \cite{apecs}. First $\psi$ factors modulo $5$ in irreducible factors
  of degree $10$, so in particular no irreducible factor of $\psi$
  over $\QQQ$ has degree $<10$. Furthermore, $-2$ is a square modulo
  $17$, and $\psi$ has two factors of degree $5$ and the remaining of
  degree $10$ modulo $17$, so $\psi$ has at most two factors of degree
  $5$ over $K$, which are furthermore interchanged by $\Gal(K|\QQQ)$.
  The assertion follows if we indeed exhibit an isogeny of degree
  $11$. First note that
\[\rho:\;(X,Y)\to(-\frac{1}{2}\,\frac{X^2-4X+6}{X-2},
                  -\frac{\sqrt{-2}}{4}Y\frac{X^2-4X+2}{(X-2)^2})
\]
defines an endomorphism of degree $2$ of $E_2$, which does not come
from dividing by a point of order $2$ (for instance because then this
isogeny were defined over $\sqrt{2}$). So if we identify $E_2$ with
the torus $\CCC/\Lambda$ with a lattice $\Lambda$ in the usual way,
then this endomorphism is induced by multiplication with $\sqrt{-2}$.
As $\abs{3+\sqrt{-2}}^2=11$, we get that $[3]+\rho$ is a $K$--rational
endomorphism (where $[3]$ denotes multiplication with $3$) of degree
$11$ of $E_2$. The dual isogeny is the requested isogeny. (Using the
addition formulae it would be possible to write out this degree $11$
isogeny explicitly.)
\end{proof}

Using this and the factorization results from above, we see that the
associated map $h(X)\in K(X)$ on the $X$--coordinates (see Section
\ref{S:Isog}) is arithmetically exceptional, with the field of
constants a cyclic extension of degree $5$ over $K$. Hence also the
function $f$ we started with has this field of constants, hence $A>G$.

\subsection{$G=(C_5^2)\rtimes S_3$ (Theorem \ref{Affine}(c)(ii))}
This case occurs over the rationals, see \cite{PM:AriEx}.

\subsection{$G=(C_5^2)\rtimes((C_4\times C_2)\rtimes C_2)$ (Theorem
\ref{Affine}(c)(iii))}

Similarly as in case (i), we see that a minimal field of definition
has to contain $K:=\QQQ(\sqrt{-1})$. Let $\cT$ be a branch cycle
description. Recall that the type is $(2,2,2,4)$. The involutions of
course are rational, and the element of order $4$ is $K$--rational.
Furthermore, $\cT$ is a weakly rigid tuple, and $G$ contains only one
conjugacy class of subgroups of index $25$. Let $\infty$ be a branch
point corresponding to the element of order $4$, and let
$\frakp_1,\frakp_2,\frakp_3$ be the roots of a cubic polynomial which
has Galois group $S_3$ and splitting field $\hK$ over $K$. Apply
Proposition \ref{WeakRigCrit} to this configuration to obtain $f(X)\in
K(X)$ with ramification data as above and geometric monodromy group
$G$.  (Note that the element of order $4$ has a unique fixed point.)
The branch cycle argument Proposition \ref{BCA} shows that $\hK$ is
the field of constants of $f$, so the arithmetic monodromy group is
$A$.

Similarly, we realize the arithmetic monodromy group $B$. Let $\hK$ be
any cyclic cubic extension of $K$, $\frakp_1$ be a primitive element
of this extension, and $\frakp_2$, $\frakp_3$ the conjugates.
Associate the three $G$--conjugacy classes $\cC_i$ with the
involutions in such a way, that the permutation of these classes by
$B$ is compatible with the permutation of the $\frakp_i$ by
$\Gal(\hK|K)$. Now apply Propositions \ref{WeakRigCrit} and \ref{BCA}
to get the conclusion.

\subsection{$G=(C_5^2)\rtimes D_{12}$ (Theorem \ref{Affine}(c)(iv))}
This is by far the most difficult case. In \cite{PM:AriEx} we show that
$B/G=C_2$ occurs (up to linear fractional transformations) two times
over the rationals. There is strong evidence that $A/G=C_4$ does also
occur over the rationals.

\subsection{$G=(C_3^2)\rtimes D_{8}$ (Theorem \ref{Affine}(c)(v))}
Here we have the possible ramification types $(2,4,6)$, $(2,2,2,6)$,
$(2,2,2,4)$, and $(2,2,2,2,2)$. The first two ramification types do
not occur over any field $K$ of characteristic $0$. Namely the place
with ramification index $6$ must be rational, so the normalizer of the
corresponding inertia group in $A$ must not be contained in $G$ by
Lemma \ref{AG}. However, the normalizer of a cyclic group of order $6$
of $G$ has order $12$ and is a subgroup of $G$, a contradiction.

The other two ramification types occur over the rationals, see
\cite{PM:AriEx}.

\subsection{$G=(C_2^4)\rtimes(C_5\rtimes C_2)$ (Theorem
  \ref{Affine}(c)(vi))}\label{deg16} Suppose that this case occurs
over some number field $K$. If $B<A$, then the field of constants
$\hat K$ has Galois group $A/G=S_3$ over $K$. Let $K'$ be the
quadratic extension of $K$ in $\hK$. Then $f\in K(X)$ has, over $K'$,
arithmetic and geometric monodromy groups $B$ and $G$ respectively. So
upon replacing $K$ by $K'$, we may assume that $A=B$. The ramification
type of $G$ is either $(2,4,5)$ or $(2,2,2,4)$. Let $I$ be the inertia
group of order $4$. One verifies that the normalizer in $A$ of each
cyclic subgroup of order $4$ in $G$ is contained in $G$. So
$\Nor_A(I)\le G$, contrary to Lemma \ref{AG}. This shows that this
case does not occur over any number field.

\subsection{$G=\PSL_2(8)$ (Theorem \ref{AS}(a))}\label{PGL28} Here
$A=\PgL_2(8)$, $G=\PSL_2(8)$, and $n=28$. There are three kinds of
genus $0$ systems, namely of types $(2,2,2,3)$, $(2,3,7)$, and
$(2,3,9)$, where the latter two come from coalescing branch points of
the former type.  We show that all cases indeed yield arithmetically
exceptional functions over the rationals.

Let $\s1,\s2,\s3$ be a generating triple of $G$, with $\s1\s2\s3=1$,
$\abs{\s1}=2$, $\abs{\s2}=3$, and $\abs{\s3}=m$ with $m=7$ or $9$. Let
$\cC_i$ be the conjugacy class of $\s i$ in $G$. One immediately
verifies that the triple $\cH=(\cC_1,\cC_2,\cC_3)$ is weakly rigid.
Furthermore, the character values of the elements $\s i$ lie in the
real field $K=\QQQ(\zeta_m+1/\zeta_m)$, which is a Galois extension of
$\QQQ$ of degree $3$. The classes $\cC_1$ and $\cC_2$ are rational in
$G$ and fixed under $A=\Aut(G)$ (simply because there is only one
class of elements of order $2$ or $3$).  Furthermore, $A$ permutes the
three conjugacy classes in $G$ with elements of order $m$
transitively. Associate a triple $\cB$ of $3$ rational points to the
$\cC_i$. Thus the pair $(\cB,\cH)$ is weakly $\QQQ$--rational. Let $H$
be a subgroup of $G$ of index $28$. Then $H$ is determined up to
conjugacy, because these subgroups are the normalizers of the Sylow
$3$--subgroups of $G$.  Also, $\s2$ has a unique fixed point. Use
Proposition \ref{WeakRigCrit} to obtain $f(X)\in\QQQ(X)$ with
geometric monodromy group $G$. Let $t$ be a transcendental, and $L$ be
a splitting field of $f(X)-t$ over $\QQQ(t)$.  We need to determine
$\tilde A:=\Gal(L|\QQQ(t))$. Let $\hat\QQQ$ be the algebraic closure
of $\QQQ$ in $L$. Then $\Gal(L|\hat Q(t))=G$. We cannot have
$\hat\QQQ=\QQQ$, for then $\tilde A=G$, but $\cC_3$ is not rational in
$\tilde A$, contrary to Corollary \ref{BCAC}. Thus $\tilde A$ properly
contains $G$, but also acts on the $28$ roots of $f(X)-t$ and
normalizes $G$.  So $\tilde A$ must be $\PgL_2(8)$.

We now construct an example with ramification type $(2,2,2,3)$. For
this let $f$ be as above of type $(2,3,9)$. After making suitable
choices of the necessarily rational branch points and some of the
necessarily rational preimages, we may assume that
\[
f(X)=\frac{X\cdot a(X)^3}{b(X)^9},
\]
with polynomials $a,b\in\QQQ[X]$. So $f(X^3)=g(X)^3$ for
$g(X)\in\QQQ(X)$ of degree $28$. One easily verifies that $g$ has
ramification type $(2,2,2,3)$. We show that $g$ has the same pair of
arithmetic and geometric monodromy group as $f$. We give the argument
only for the arithmetic monodromy group, a slight variation handles
the geometric group too.

Let $t$ be a transcendental, $x$ a root of $f(X^3)-t$, $L$ a Galois
closure of $\QQQ(x)|\QQQ(t)$ with group $H=\Gal(L|\QQQ(t))$. Set
$U_f:=\Gal(L|\QQQ(x^3))$, $M:=\Gal(L|g(x))$, $U_g:=\Gal(L|\QQQ(x))$,
$N_f:=\core_H(U_f)$, and $N_g:=\core_H(U_g)$. So the arithmetic
monodromy group of $f$ and $g$ is $H/N_f=\PgL_2(8)$ and $M/N_g$,
respectively.

We first claim that $H=MN_f$. Suppose false. Then $N_f\le U_f\cap
M=U_g$, but the core of $U_g$ in $H$ is trivial, so $N_f=1$. So
$H=\PgL_2(8)$ has a non--normal subgroup $M$ of index $3$ (note that
$M<U$ corresponds to the non--normal cubic field extension
$\QQQ(g(x))|\QQQ(g(x)^3)$), a contradiction.

From $H=MN_f$ and $N_f\le U_f$ we obtain
\begin{align*}
N_g &= \bigcap_{m\in M}U_g^m         \\
    &= \bigcap_{m\in M}(U_f\cap M)^m \\
    &= (\bigcap_{m\in M}U_f^m)\cap M \\
    &= (\bigcap_{h\in H}U_f^h)\cap M \\
    &= N_f \cap M,
\end{align*}
so the claim follows from
\[
H/N_f=MN_f/N_f=M/(N_f\cap M)=M/N_g.
\]

\subsection{$G=\PSL_2(9)$ (Theorem \ref{AS}(b))}\label{PSL29} Here
$\M10\le A\le\PgL_2(9)$, $G=\PSL_2(9)$, and $n=45$. The genus $0$
systems of $G$ are of type $(2,4,5)$. They are weakly rigid. Let
$\cC_1$ and $\cC_2$ be the conjugacy classes of elements of order $2$
and $4$, and $\cC_3$ and $\bar\cC_3$ be the two conjugacy classes of
elements of order $5$.

Now let $(\s1,\s2,\s3)$ be a generating triple of $G$ with $\s
i\in\cC_i$. The classes $\cC_1$ and $\cC_2$ are rational in $G$ and
fixed under $\M10<\Aut(G)$, whereas $\M10$ switches $\cC_3$ and
$\bar\cC_3$. Therefore $(\cC_1,\cC_2,\cC_3)$ is weakly
$\QQQ$--rational.  The group $G$ has only one conjugacy class of
subgroups $H$ of index $45$, as these subgroups are the Sylow
$2$--subgroups of $G$. As $\s2$ has a unique fixed point we can apply
Proposition \ref{WeakRigCrit} to obtain $f(X)\in\QQQ(X)$ with
geometric monodromy group $G$ and given ramification. Let $t$ be a
transcendental, and $L$ be a splitting field of $f(X)-t$ over
$\QQQ(t)$.  As in the previous case, we need to determine the
arithmetic monodromy group $\tilde A=\Gal(L|\QQQ(t))$. The group $G$
in its given action is normalized by $\PgL_2(9)$. First note that
Corollary \ref{BCAC} shows that $\tilde A$ is not contained in
$\PGL_2(9)$, for $\tilde A$ must switch the conjugacy classes $\cC_3$
and $\bar\cC_3$. Suppose that $f$ is not arithmetically exceptional.
Then $\tilde A=\PsL_2(9)$. So $\tilde A$ acts imprimitively, moving
$15$ blocks of size $3$. Thus we have $f(X)=a(b(X))$ with
$a(X)\in\QQQ(X)$ of degree $15$.  Therefore $a$ has geometric
monodromy group $A_6\cong\PSL_2(9)$. But the $(2,4,5)$ system also has
genus $0$ with respect to the natural action of $A_6$. Let $H$ be a
point stabilizer in $A_6$ of this action, and $J$ the normalizer of
$H$ in the arithmetic monodromy group (which is $A_6$ or $S_6$) of
$a$. Then the fixed field of $J$ gives a rational function
$h(X)\in\QQQ(X)$ of degree $6$ and ramification type $(2,4,5)$.
Suppose that $\infty$ and $0$ correspond to the inertia generators of
order $5$ and $4$ respectively. Without loss let the single point in
$h^{-1}(\infty)$ be $0$, and the quintuple point be $\infty$.
Furthermore, assume that the quadruple point in $h^{-1}(0)$ is $1$,
and the other one is $\alpha\in\QQQ$. So
$$h(X)=\frac{(X-1)^4(X-\alpha)^2}{X}.$$
The third branch point of $h$
has to be rational, and is the root $\ne0$ of the discriminant with
respect to $X$ of the numerator of $h(X)-t$.  Dividing by $t^4$ we get
$3125t^2+(-256+15060\alpha+33360\alpha^2-2120\alpha^3+
720\alpha^4-108\alpha^5)t-1024\alpha+
6144\alpha^2-15360\alpha^3+20480\alpha^4-15360\alpha^5+
6144\alpha^6-1024\alpha^7$. So this polynomial has only one root,
therefore its discriminant
$16(\alpha^2-11\alpha+64)^2(9\alpha^2+26\alpha+1)^3$ vanishes.
However, this polynomial does not have a rational solution $\alpha$.
This contradiction proves that $\M10\le\tilde A$.

\subsection{A remark about one of the sporadic cases}
The degree $28$ case in Theorem \ref{main}, which corresponds to a
$(2,3,7)$ generating system of $G=\PSL_2(8)$ is somehow classically
known. Note that $A=\PgL_2(8)$ contains a subgroup $U$ of index $9$
with $A=UG$, such that the $(2,3,7)$ system induces genus $0$ on
$A/U$. (The induced genera are $>0$ for the two other genus $0$
systems though.) Thus there should be a rational function
$f(X)\in\QQQ(X)$ of degree $9$ and $(A,G)$ as arithmetic and geometric
monodromy group. Such a function has been computed frequently in the
past by Goursat, J-F.~Mestre, and others, see \cite{Serre:Abhy}. In
\cite{Serre:Abhy}, Serre sketches a computation of $f$. Noting that
the coefficients of $f$ happen to lie in $\QQQ$, he asks ``Could this
be proved a priori? I suspect it could.'' Proposition
\ref{WeakRigCrit}, and the way we proved existence without
computation, gives an affirmative answer. Later on Serre speculates
about a ``modular'' interpretation of this case and a connection to
Shimura curves. A recent preprint \cite{Elkies:Shimura} by Elkies
contains a lot of material in this regard. A similar interpretation
holds for the degree $45$ case (Elkies, personal communication).

It would be interesting to find nice interpretations for the other
sporadic cases too.

\cleardoublepage
\addcontentsline{toc}{section}{Bibliography}
\newcommand{\etalchar}[1]{$^{#1}$}

\vskip.5cm

\noindent{\sc Department of Mathematics, University of Southern California,
Los Angeles, CA 90089}, USA\par
\noindent{\sl E-mail: }{\tt guralnic@math.usc.edu}\vskip.3cm\par
\noindent{\sc IWR, Universit\"at Heidelberg, Im Neuenheimer Feld 368,
  69120 Heidelberg, Germany}\par
\noindent{\sl E-mail: }{\tt
Peter.Mueller@iwr.uni-heidelberg.de}\vskip.3cm\par
\noindent{\sc DPMMS, University of Cambridge, 16 Mill Lane,
  Cambridge CB2 1SB, England}\par
\noindent{\sl E-mail: }{\tt saxl@dpmms.cam.ac.uk}

\begin{thebibliography}{CCN{\etalchar{+}}85}

\bibitem[Asc86]{Asch:Book}
M.~Aschbacher, \textit{Finite Group Theory}, Cambridge University Press,
  Cambridge (1986).

\bibitem[Asc92]{Asch:GT}
M.~Aschbacher, \textit{On conjectures of {G}uralnick and {T}hompson}, J.
  Algebra (1992), \textbf{135}, 277--341.

\bibitem[AS85]{AschScott}
M.~Aschbacher, L.~Scott, \textit{Maximal subgroups of finite groups}, J.
  Algebra (1985), \textbf{92}, 44--88.

\bibitem[BGL77]{BGL}
N.~Burgoyne, R.~Griess, R.~Lyons, \textit{Maximal subgroups and automorphisms
  of {C}hevalley groups}, Pacific J. Math. (1977), \textbf{71}(2), 365--403.

\bibitem[Che51]{Chevalley}
C.~Chevalley, \textit{Algebraic Functions of one Variable}, Mathematical
  Surveys VI, Amer. Math. Soc., Providence (1951).

\bibitem[CM94]{CohenMatthews}
S.~D. Cohen, R.~W. Matthews, \textit{A class of exceptional polynomials},
  Trans. Amer. Math. Soc. (1994), \textbf{345}, 897--909.

\bibitem[Con97]{apecs}
I.~Connell, \textit{Apecs (arithmetic of plane elliptic curves), a program
  written in maple}, available via anonymous ftp from math.mcgill.ca in
  /pub/apecs (1997).

\bibitem[At85]{ATLAS}
J.~H. Conway, R.~T. Curtis, S.~P. Norton, R.~A. Parker, R.~A. Wilson,
  \textit{Atlas of Finite Groups}, Clarendon Press, Oxford (1985).

\bibitem[Elk98]{Elkies:Shimura}
N.~Elkies, \textit{Shimura curve computations}, in J.~P. Buhler, ed.,
  \textit{Algorithmic Number Theory}, vol. 1423 of \textit{Lecture Notes in
  Computer Science}. Springer (1998), 1998 pp. 1--47.

\bibitem[Eno72]{Enomoto}
H.~Enomoto, \textit{The characters of the finite symplectic group {Sp}$_4(q)$,
  $q=2^f$}, Osaka J. Math. (1972), \textbf{9}, 75--94.

\bibitem[Fri70]{Fried:Schur}
M.~Fried, \textit{On a conjecture of {S}chur}, Michigan Math. J. (1970),
  \textbf{17}, 41--55.

\bibitem[Fri74]{Fried:HIT}
M.~Fried, \textit{On {H}ilbert's irreducibility theorem}, J. Number Theory
  (1974), \textbf{6}, 211--231.

\bibitem[Fri78]{Fried:CM}
M.~Fried, \textit{Galois groups and complex multiplication}, Trans. Amer. Math.
  Soc. (1978), \textbf{235}, 141--163.

\bibitem[FGS93]{FGS}
M.~Fried, R.~M. Guralnick, J.~Saxl, \textit{{S}chur covers and {C}arlitz's
  conjecture}, Israel J. Math. (1993), \textbf{82}, 157--225.

\bibitem[GAP95]{GAP}
M.~Sch{\"o}nert, et~al., \textit{GAP -- Groups, Algorithms, and Programming},
  {Lehrstuhl D f{\"u}r Mathematik}, RWTH Aachen, Germany (1995).

\bibitem[Gor68]{Gorenstein}
D.~Gorenstein, \textit{Finite Groups}, Harper and Row, New
  York--Evanston--London (1968).

\bibitem[GL83]{GorLy}
D.~Gorenstein, R.~Lyons, \textit{The local structure of finite groups of
  characteristic $2$\ type}, Mem. Amer. Math. Soc. (1983), \textbf{42}(276),
  vii+731.

\bibitem[GLS98]{GLS3}
D.~Gorenstein, R.~Lyons, R.~Solomon, \textit{The classification of the finite
  simple groups. {N}umber 3. {P}art {I}. {C}hapter {A}}, American Mathematical
  Society, Providence, RI (1998), Almost simple $K$-groups.

\bibitem[Gro71]{SGA1}
A.~Grothendieck, \textit{Rev{\^e}tement {\'e}tales et groupe fondamental
  ({SGA}1)}, vol. 224 of \textit{Lecture Notes in Math.}, Springer--Verlag
  (1971).
  
\bibitem[GLPS]{GLPS} R.~M. Guralnick, C.~H. Li, C.~E. Praeger,
  J.~Saxl, \textit{Exceptional primitive group actions and partitions
    of orbitals}, preprint.

\bibitem[GM97]{GM:ExPol}
R.~M. Guralnick, P.~M{\"u}ller, \textit{Exceptional polynomials of affine
  type}, J. Algebra (1997), \textbf{194}, 429--454.

\bibitem[GN95]{GurNeu}
R.~M. Guralnick, M.~Neubauer, \textit{Monodromy groups of branched coverings:
  The generic case}, in M.~Fried, ed., \textit{Recent developments in the
  inverse {G}alois problem}, vol. 186 of \textit{Contemp. Maths.} Amer. Math.
  Soc. (1995), 1995 pp. 325--352.

\bibitem[GN01]{GurNeu:Affine}
R.~M. Guralnick, M.~Neubauer, \textit{Arithmetically indecomposable rational
  functions of affine type}, preprint.

\bibitem[GRZ01]{GRZ}
R.~Guralnick, J.~Rosenberg, M.~Zieve, \textit{A new class of exceptional
  polynomials in characteristic $2$}, preprint.

\bibitem[GS01]{GSRam}
R.~M. Guralnick, J.~Saxl, \textit{Exceptional polynomials over arbitrary
  fields}, preprint.

\bibitem[GT90]{GT}
R.~M. Guralnick, J.~G. Thompson, \textit{Finite groups of genus zero}, J.
  Algebra (1990), \textbf{131}, 303--341.

\bibitem[GW97]{GW}
R.~Guralnick, D.~Wan, \textit{Bounds for fixed point free elements in a
  transitive group and applications to curves over finite fields}, Israel J.
  Math. (1997), \textbf{101}, 255--287.

\bibitem[GZ]{GZ}
R.~M. Guralnick, M.~Zieve, \textit{Exceptional rational functions of small
  genus}, in preparation.

\bibitem[HY79]{HarYam}
K.~Harada, H.~Yamaki, \textit{The irreducible subgroups of {GL}$_n(2)$ with
  $n\le6$}, C. R. Math. Rep. Acad. Sci. Canada (1979), \textbf{1}, 75--78.

\bibitem[Hup67]{Huppert1}
B.~Huppert, \textit{Endliche Gruppen I}, Springer--Verlag, Berlin Heidelberg
  (1967).

\bibitem[KM95]{KamMaz}
S.~Kamienny, B.~Mazur, \textit{Rational torsion of prime order in elliptic
  curves over number fields}, Ast{\'e}risque (1995), \textbf{228}, 81--98.

\bibitem[Kan79]{Kantor:RootElements}
W.~M. Kantor, \textit{Subgroups of classical groups generated by long root
  elements}, Trans. Amer. Math. Soc. (1979), \textbf{248}, 347--379.

\bibitem[Lan78]{Lang:EC}
S.~Lang, \textit{Elliptic Curves Diophantine Analysis}, Springer--Verlag, New
  York (1978).

\bibitem[Lan86]{Lang:AN}
S.~Lang, \textit{Algebraic Number Theory}, Springer--Verlag, New York (1986).

\bibitem[Lan87]{Lang:EF}
S.~Lang, \textit{Elliptic Functions}, Springer--Verlag, New York (1987).

\bibitem[LZ96]{LZ}
H.~W. {Lenstra, Jr.}, M.~Zieve, \textit{A family of exceptional polynomials in
  characteristic three}, in S.~Cohen, H.~Niederreiter, eds., \textit{Finite
  Fields and Applications}, vol. 233 of \textit{Lecture Note Series}, Cambridge
  (1996), Cambridge University Press, 1996 pp. 209--218.

\bibitem[LP]{LiPr}
C.~H. Li, C.~E. Praeger, \textit{On partitioning the orbitals of a transitive
  permutation group}, preprint.

\bibitem[LPS88]{LPS}
M.~W. Liebeck, C.~E. Praeger, J.~Saxl, \textit{On the {O'Nan--Scott Theorem}
  for finite primitive permutation groups}, J. Austral. Math. Soc. Ser. A
  (1988), \textbf{44}, 389--396.

\bibitem[LS91]{LS:MinDeg}
M.~W. Liebeck, J.~Saxl, \textit{Minimal degrees of primitive permutation
  groups, with an application to monodromy groups of covers of {R}iemann
  surfaces}, Proc. London Math. Soc. (3) (1991), \textbf{63}, 266--314.

\bibitem[Mag93]{Magaard:Sporadic}
K.~Magaard, \textit{Monodromy and sporadic groups}, Comm. Algebra (1993),
  \textbf{21}(12), 4271--4297.

\bibitem[Mag74]{Magnus:Tess}
W.~Magnus, \textit{Noneuclidean Tesselations and their Groups}, Academic Press,
  New York (1974).

\bibitem[MM99]{MM}
G.~Malle, B.~H. Matzat, \textit{{Inverse Galois Theory}}, Springer Verlag,
  Berlin (1999).

\bibitem[Mat84]{Matthews}
R.~Matthews, \textit{Permutation polynomials over algebraic number fields}, J.
  Number Theory (1984), \textbf{18}, 249--260.

\bibitem[Maz78]{Mazur}
B.~Mazur, \textit{Rational isogenies of prime degree}, Invent. Math. (1978),
  \textbf{44}, 129--162.

\bibitem[McL69]{McL:Some}
J.~McLaughlin, \textit{Some subgroups of {SL}$_n(\mathbb {F}_2)$}, Illinois J.
  Math. (1969), \textbf{13}, 108--115.

\bibitem[M{\"u}l94]{PM:Ex28}
P.~M{\"u}ller, \textit{New examples of exceptional polynomials}, in G.~L.
  Mullen, P.~J. Shiue, eds., \textit{Finite Fields: Theory, Applications and
  Algorithms}, vol. 168 of \textit{Contemp. Maths.} Amer. Math. Soc. (1994),
  1994 pp. 245--249.

\bibitem[M{\"u}l99]{PM:AriEx}
P.~M{\"u}ller, \textit{Arithmetically exceptional functions and elliptic
  curves}, in {H. V\"olklein, P. M\"uller, D. Harbater, J. G. Thompson}, ed.,
  \textit{Aspects of Galois Theory}, vol. 256 of \textit{Lecture Note Series},
  Cambridge (1999), Cambridge University Press, 1999 pp. 180--201.

\bibitem[Neu93]{Neubauer:CA1}
M.~Neubauer, \textit{On primitive monodromy groups of genus zero and one, {I}},
  Comm. Algebra (1993), \textbf{21}(3), 711--746.

\bibitem[N{\"o}b89]{Nobauer}
R.~N{\"o}bauer, \textit{R{\'e}dei--{F}unktionen und das {S}chur'sche
  {P}roblem}, Arch. Math. (Basel) (1989), \textbf{52}, 61--65.

\bibitem[Row95]{Rowley}
P.~Rowley, \textit{Finite groups admitting a fixed-point-free automorphism
  group}, J. Algebra (1995), \textbf{174}, 724--727.

\bibitem[Rub99]{Rubin:CM}
K.~Rubin, \textit{Elliptic curves with complex multiplication and the
  conjecture of {B}irch and {S}winnerton-{D}yer}, in \textit{Arithmetic theory
  of elliptic curves (Cetraro, 1997)}, Springer, Berlin, 1999 pp. 167--234.

\bibitem[Sch23]{Schur}
I.~Schur, \textit{{\"U}ber den {Z}usammenhang zwischen einem {P}roblem der
  {Z}ahlentheorie und einem {S}atz {\"u}ber algebraische {F}unktionen}, S.-B.
  Preuss. Akad. Wiss., Phys.--Math. Klasse (1923), pp. 123--134.

\bibitem[Ser79]{Serre:LF}
J.-P. Serre, \textit{Local Fields}, Springer--Verlag, New York (1979).

\bibitem[Ser94]{Serre:Abhy}
J.-P. Serre, \textit{A letter as an appendix to the square--root
  parameterization paper of {A}bhyankar}, in C.~L. Bajaj, ed.,
  \textit{Algebraic Geometry and its Applications}, Springer--Verlag, 1994 pp.
  85--88.

\bibitem[Shi74]{Shih:PSL}
K.~Shih, \textit{On the construction of {G}alois extensions of function fields
  and number fields}, Math. Ann. (1974), \textbf{207}, 99--120.

\bibitem[Shi91]{Shih}
T.~Shih, \textit{A note on groups of genus zero}, Comm. Algebra (1991),
  \textbf{19}, 2813--2826.

\bibitem[Sil86]{Silverman1}
J.~H. Silverman, \textit{The Arithmetic of Elliptic Curves}, Springer--Verlag,
  New York (1986).

\bibitem[Sil94]{Silverman2}
J.~H. Silverman, \textit{Advanced Topics in the Arithmetic of Elliptic Curves},
  Springer--Verlag, New York (1994).

\bibitem[Sti93]{Stich}
H.~Stichtenoth, \textit{Algebraic Function Fields and Codes}, Springer--Verlag,
  Berlin Heidelberg (1993).

\bibitem[Suz82]{Suzuki1}
M.~Suzuki, \textit{Group Theory I}, Springer--Verlag, Berlin (1982).

\bibitem[Tur95]{Turnwald:Schur}
G.~Turnwald, \textit{On {S}chur's conjecture}, J. Austral. Math. Soc. Ser. A
  (1995), \textbf{58}, 312--357.

\bibitem[V{\"o}l96]{H:Buch}
H.~V{\"o}lklein, \textit{Groups as Galois Groups -- an Introduction}, Cambridge
  University Press, New York (1996).

\end{thebibliography}
\end{document}